\PassOptionsToPackage{unicode}{hyperref}
\PassOptionsToPackage{noend}{algorithmic}
\documentclass[USenglish,journal]{IEEEtran}
\usepackage[%
	noload={authSort},
]{myPreamble}
\ifieee
	\crefname{appsec}{Appendix}{Appendices}
	\Crefformat{appsec}{#2Appendix #1#3}
\fi
\newcommand\added[1]{#1}
\newcommand\removed[1]{}

	\newcommand\weak{\mathcal W}
	\newcommandx\linop[2][{1=\HH},{2={}}]{%
		\ifstrempty{#2}{%
			\ifstrempty{#1}{\mathcal B}{\mathcal B(#1)}%
		}{%
			\mathcal B(#1;#2)%
		}%
	}
	\newcommandx\HS[2][{1=\HH},{2={}}]{%
		\ifstrempty{#2}{%
			\ifstrempty{#1}{\mathcal S}{\mathcal S(#1)}%
		}{%
			\mathcal S(#1;#2)%
		}%
	}

	\newcommand\avg[1][\alpha]{%
		\ifstrempty{#1}{averaged}{\ensuremath{#1}-av\-er\-aged}%
	}
	\newcommand\KM{\texorpdfstring{Krasnosel'ski\v{\i}-Mann}{Krasnosel'skii-Mann}}
	\newcommand\refSM[1][{ (\cref{alg:SuperMann})}]{%
		\textit{%
			\hyperref[alg:SuperMann]{Su}\-
			\hyperref[alg:SuperMann]{per}\-%
			\hyperref[alg:SuperMann]{Mann}
			\hyperref[alg:SuperMann]{scheme}#1%
		}%
	}

	\newcommandx\boxedalgorithmic[3][{1={}},{3={draw,black,thin,fill=black!10,rounded corners}}]{{%
		\tikzexternaldisable%
		\begin{tikzpicture}%
			\node[#3] {%
				\begin{minipage}[t]{0.75\linewidth}%
					#1\begin{algorithmic}[1]#2\end{algorithmic}%
				\end{minipage}%
			};%
		\end{tikzpicture}%
	}}%

	\usepackage{pgfplots}
	\pgfplotsset{compat=1.10}
	\usepgfplotslibrary{fillbetween}
	\usetikzlibrary{arrows,arrows.meta,calc,decorations.markings,patterns,positioning,through,shapes}
	\newcounter{axiscounter}
	\definecolor{mycolor1}{rgb}{0.00000,0.44700,0.74100}%
	\definecolor{mycolor2}{rgb}{0.85000,0.32500,0.09800}%
	\definecolor{mycolor3}{rgb}{0.92900,0.69400,0.12500}%
	\definecolor{mycolor4}{rgb}{0.49400,0.18400,0.55600}%

	\pgfplotscreateplotcyclelist{plotcycle}{
		{mycolor1},
		{mycolor2},
		{mycolor3},
		{mycolor4},
	}

	\pgfplotsset{%
		commonaxis/.style={%
			cycle list name=plotcycle,
			scale only axis,
			axis background/.style={fill=white},
			title style={font=\bfseries\vphantom{Ag}},
			ylabel style={font=\sffamily\vphantom{Ag}},
			xlabel style={font=\sffamily\vphantom{Ag}},
			legend style={%
				at={(1,0)},
				anchor=south east,
				outer sep=2pt,
				legend cell align=left,
				align=left,
				draw=white!15!black,
				font={\sffamily},
				scale=2,
				opacity=0.8,
			},
		},
		graphaxis/.style={%
			commonaxis,
			unit vector ratio*=1 1 1,
			hide axis,
			legend style={%
				nodes={scale=1.5, transform shape},
			},
		},
		plotaxis/.style={%
			commonaxis,
			/pgfplots/ymode=log,
			yminorticks=true,
			xmajorgrids,
			ymajorgrids,
			yminorgrids,
			grid=major,
			solid,
			xmin=1,
			width=3.0in,
			height=4.5in,
			at={(\theaxiscounter in,0)},
		},
		plotaxis/.belongs to family=/pgfplots/scale,
		plotstyle/.style = {%
			line width=2.0pt,
			solid,
			opacity=0.8,
		},
	}
	\AtBeginEnvironment{tikzpicture}{\setcounter{axiscounter}{0}}

	\pgfplotsset{%
		every axis/.append code = {\addtocounter{axiscounter}{4}},
	}
	\newcommand\getDistance[2]{let \p1=(#1), \p2=(#2), \n1={veclen((\x2-\x1),(\y2-\y1))/2} in}

	\newenvironment{customlegend}[1][]{%
		\begingroup
		\csname pgfplots@init@cleared@structures\endcsname
		\pgfplotsset{#1}%
	}{%
		\csname pgfplots@createlegend\endcsname
		\endgroup
	}%

	\def\addlegendimage{\csname pgfplots@addlegendimage\endcsname}

	\pgfkeys{/pgfplots/number in legend/.style={%
			/pgfplots/legend image code/.code={%
				\node at (0.125,-0.0225){#1}; 
			},%
		},%
	}%
	\pgfplotsset{%
		every legend to name picture/.style={west}
	}%

\begin{document}
	\author{%
		Andreas Themelis and Panagiotis Patrinos%
		\thanks{%
			\TheAddressKU\ 
			\{%
				\href{mailto:andreas.themelis@esat.kuleuven.be}{andreas.themelis},
				\href{mailto:panos.patrinos@esat.kuleuven.be}{panos.patrinos}%
			\}%
			\href{mailto:andreas.themelis@esat.kuleuven.be,panos.patrinos@esat.kuleuven.be}{@esat.kuleuven.be}.%
		}%
	}%
	\title{SuperMann: a superlinearly convergent algorithm for finding fixed points of nonexpansive operators}%
	\maketitle

	\begin{abstract}
		Operator splitting techniques have recently gained popularity in convex optimization problems arising in various control fields.
		Being fixed-point iterations of nonexpansive operators, such methods suffer many well known downsides, which include high sensitivity to ill conditioning and parameter selection, and consequent low accuracy and robustness.
		As universal solution we propose \emph{SuperMann}, a Newton-type algorithm for finding fixed points of nonexpansive operators.
		It generalizes the classical \KM{} scheme, enjoys its favorable global convergence properties and requires exactly the same oracle.
		It is based on a novel separating hyperplane projection tailored for nonexpansive mappings which makes it possible to include steps along any direction.
		In particular, when the directions satisfy a Dennis-Moré condition we show that \emph{SuperMann} converges superlinearly under mild assumptions, which, surprisingly, do not entail nonsingularity of the Jacobian at the solution but merely metric subregularity.
		As a result, \emph{SuperMann} enhances and robustifies all operator splitting schemes for structured convex optimization, overcoming their well known sensitivity to ill conditioning.
	\end{abstract}


	\section{Introduction}
		\label{sec:Intro}
		Operator splitting techniques (also known as proximal algorithms), introduced in the 50's for solving PDEs and optimal control problems, have been successfully used to reduce complex problems into a series of simpler subproblems.
The most well known operator splitting methods are the alternating direction method of multipliers (ADMM), forward-backward splitting (FBS) also known as proximal-gradient method in composite convex minimization, Douglas-Rachford splitting (DRS) and the alternating minimization method (AMM) \cite{parikh2014proximal}.
Operator splitting techniques pose several advantages over traditional optimization methods such as sequential quadratic programming and interior point methods:
(1) they can easily handle nonsmooth terms and abstract linear operators, (2) each iteration requires only simple arithmetic operations, (3) the algorithms scale gracefully as the dimension of the problem increases, and (4) they naturally lead to parallel and distributed implementation.
Therefore, operator splitting methods cope well with limited amount of hardware resources making them particularly attractive for (embedded) control \cite{stathopoulos2016operator}, signal processing \cite{combettes2011proximal}, and distributed optimization \cite{benameur2016robust,iutzeler2016explicit}.

The key idea behind these techniques when applied to convex optimization is to reformulate the optimality conditions of the problem at hand into a problem of finding a fixed point of a nonexpansive operator and then apply relaxed fixed-point iterations.
Although sometimes a fast convergence rate can be observed, the norm of the fixed-point residual decreases, at best, with \(Q\)-linear rate, and due to an inherent sensitivity to ill conditioning oftentimes the \(Q\)-factor is close to one.
Moreover, all operator splitting methods are basically ``open-loop'', since the tuning parameters, such as stepsizes and preconditioning, must be set before their execution.
In fact, such methods are very sensitive to the choice of parameters and sometimes there is not even a concrete way of selecting them, as it is the case of ADMM.
All these are serious obstacles when it comes to using such types of algorithms for real-time applications such as embedded MPC, or to reliably solve cone programs.

As an attempt to solve the issue, people have considered the employment of variable metrics to reshape the geometry of the problem and enhance convergence rate \cite{combettes2014variable}.
However, unless such metrics have a very specific structure, even for simple problems the cost of operating in the new geometry outweights the benefits.

Another interesting approach that is gaining more and more popularity tries to exploit possible sparsity patterns by means of chordal decomposition techniques \cite{zheng2017fast}.
These methods can improve scalability and reduce memory usage, but unless the problem comes with an inherent sparse structure they yield no tangible benefit.

Alternatively, the task of searching fixed points of an operator \(T\) can be translated to that of finding zeros of the corresponding residual \(R=\id-T\).
Many methods with fast asymptotic convergence rates such as Newton-type exist that can be employed for efficiently solving nonlinear equations, see, \eg \cite[§7]{facchinei2003finite} and \cite{izmailov2014newton}.
However, such methods converge only when close enough to the solution, and in order to globalize the convergence there comes the need of a merit function to perform a line search along candidate directions of descent.
The typical choice of the square residual \(\|Rx\|^2\) unfortunately is of no use, as in meaningful applications \(R\) is nonsmooth.

		\subsection{Proposed methodology}
			In response to these issues, in this paper we propose a universal scheme that globalizes Newton-type methods for finding fixed points of any nonexpansive operator on real Hilbert spaces.
Admittedly with an intended pun, since it exhibits \emph{super}linear convergence rates and generalizes the Krasnosel'ski\v{\i}-\emph{Mann} iterations we name our algorithm \emph{SuperMann}.
The method is based on a novel hyperplane projection step tailored for nonexpansive mappings.

Furthermore, we consider a modified Broyden's scheme which was first introduced in \cite{powell1970numerical} and show how it fits into our framework enabling superlinear asymptotic convergence rates.
One of the most appealing properties of \emph{SuperMann} is that thanks to its quasi-Fejérian behavior, achieving superlinear convergence does not necessitate nonsingularity of the Jacobian at the solution, which is the usual requirement of quasi-Newton schemes, but merely metric subregularity.
This property considerably widens the range of problems which can be solved efficiently, in that, for instance, the solutions need not be isolated for superlinear convergence to take place.

		\subsection{Contributions}
Our contributions can be summarized as follows:
\begin{enumerate}[%
	label=(\arabic*),
]
\item\label{item:contributions:General}
	In \Cref{sec:General} we design a universal algorithmic framework (\Cref{alg:General}) for finding fixed points of nonexpansive operators, which generalizes the classical \KM{} scheme and possessess its same global and local convergence properties.
\item
	In \Cref{sec:GenMann} we introduce a novel separating hyperplane projection	tailored for nonexpansive mappings;
	based on this, in \Cref{defin:GKM} we then propose a generalized KM iteration (GKM).
\item
	We define a line search based on the novel projection, suited for any nonexpansive operator and update direction (\Cref{thm:LS}).
\item
	In \Cref{sec:SuperMann} we combine these ideas and derive the \refSM, an algorithm that
	\begin{itemize}
	\item
		globalizes the convergence of Newton-type methods for finding fixed points of nonexpansive operators (\Cref{thm:SuperMann});
	\item
		reduces to the local method \(x_{k+1}=x_k+d_k\) when the directions \(d_k\) are \emph{superlinear}, as it is the case for a modified Broyden's scheme (\Cref{thm:dSuperlinear,thm:SuperMann:Broyden});
	\item
		has superlinear convergence guarantees even without the usual requirement of nonsingularity of the Jacobian at the limit point, but simply under metric subregularity; in particular, the solution need not be unique!
	\end{itemize}
\end{enumerate}

		\subsection{Paper organization}
			The paper is organized as follows.
\Cref{sec:Motivations} serves as an informal introduction to highlight the known limitations of fixed-point iterations and to motivate our interest in Newton-type methods with some toy examples.
The formal presentation begins in \Cref{sec:Definitions} with the introduction of some basic notation and known facts.
In \Cref{sec:General} we define the problem at hand and propose a general abstract algorithmic framework for solving it.
In \Cref{sec:GenMann} we provide a generalization of the classical KM iterations that is key for the global convergence and performance of \emph{SuperMann}, an algorithm which is presented and analyzed in \Cref{sec:SuperMann}.
Finally, in \Cref{sec:Simulations} we show how the theoretical findings are backed up by promising numerical simulations, where \emph{SuperMann} dramatically improves classical splitting schemes.
For the sake of readability some of the proofs are referred to the Appendix.

	\section{Motivating examples}
		\label{sec:Motivations}%
		Given a nonexpansive operator \(\func T{\R^n}{\R^n}\), consider the problem of finding a fixed point, \ie, a point \(x_\star\in\R^n\) such that \(x_\star=Tx_\star\).
The independent works of Krasnosel'ski\v{\i} and Mann \cite{krasnoselskii1955two,mann1953mean} provided a very elegant solution which is simply based on recursive iterations \(x^+=(1-\alpha)x+\alpha Tx\) with \(\alpha\in(0,\bar\alpha)\) for some \(\bar\alpha\geq 1\).
The method, known as \KM{} scheme or KM scheme for short, has been studied intensively ever since, also because it generalizes a plethora of optimization algorithms.
It is well known that the scheme is globally convergent with square-summable and monotonically decreasing residual \(R=\id-T\) (in norm), and also locally \(Q\)-linearly convergent if \(R\) is \DEF{metrically subregular} at the limit point \(x_\star\).
Metric subregularity basically amounts to requiring the distance from the set of solutions to be upper bounded by a multiple of the norm of \(R\) for all points sufficiently close to \(x_\star\); it is quite mild a requirement --- for instance, it does not entail \(x_\star\) to be an isolated solution --- and as such linear convergence is 
quite frequent in practice.
However, the major drawback of the KM scheme is its high sensitivity to ill conditioning of the problem, and cases where convergence is prohibitively slow in practice despite the theoretical (sub)linear rate are also abundant.
\ifieee
	\begin{figure*}
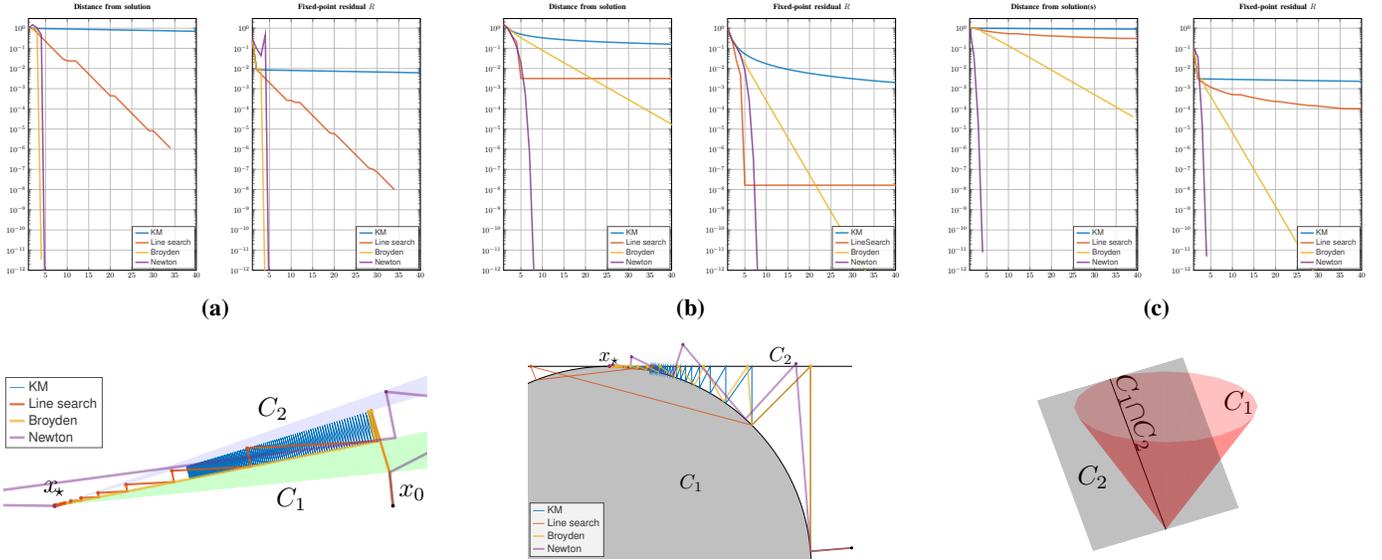

		\def\myHeight{7\baselineskip}%
		\begin{subfigure}[t]{0.31\textwidth}%
			{{%
		\pgfkeys{/pgf/images/include external/.code={\includegraphics[width=\linewidth]{#width=\linewidth}}}%
		\tikzsetnextfilename{Cones}%
%
%
\begin{tikzpicture}
\pgfplotsset{%
	every axis/.append style={
		plotaxis,
		xmax=40,
		ymin=1e-12,
		ymax=3,
	},
	every axis plot/.append style={
		plotstyle,
	},
}

\begin{axis}[%
	title={Distance from solution},
]
\addplot
  table[row sep=crcr]{%
1	1\\
2	0.976793142354244\\
3	0.96817640680066\\
4	0.959635683381458\\
5	0.951170301559109\\
6	0.942779596711199\\
7	0.934462910078249\\
8	0.926219588711997\\
9	0.918048985424136\\
10	0.909950458735496\\
11	0.901923372825691\\
12	0.893967097483192\\
13	0.886081008055853\\
14	0.878264485401867\\
15	0.870516915841159\\
16	0.862837691107203\\
17	0.855226208299271\\
18	0.847681869835092\\
19	0.840204083403943\\
20	0.832792261920139\\
21	0.825445823476948\\
22	0.818164191300899\\
23	0.810946793706501\\
24	0.803793064051363\\
25	0.7967024406917\\
26	0.789674366938244\\
27	0.782708291012536\\
28	0.775803666003604\\
29	0.768959949825026\\
30	0.762176605172371\\
31	0.755453099481012\\
32	0.74878890488432\\
33	0.742183498172215\\
34	0.735636360750088\\
35	0.729146978598094\\
36	0.722714842230785\\
37	0.716339446657119\\
38	0.710020291340807\\
39	0.703756880161019\\
40	0.69754872137343\\
41	0.691395327571618\\
42	0.685296215648791\\
43	0.679250906759864\\
44	0.673258926283859\\
45	0.667319803786648\\
46	0.661433072984014\\
47	0.655598271705044\\
48	0.649814941855846\\
49	0.644082629383582\\
50	0.638400884240819\\
51	0.632769260350198\\
52	0.627187315569409\\
53	0.621654611656482\\
54	0.616170714235377\\
55	0.610735192761882\\
56	0.605347620489812\\
57	0.600007574437502\\
58	0.594714635354603\\
59	0.58946838768916\\
60	0.584268419554993\\
61	0.579114322699356\\
62	0.574005692470886\\
63	0.568942127787832\\
64	0.563923231106571\\
65	0.558948608390387\\
66	0.554017869078545\\
67	0.549130626055622\\
68	0.544286495621116\\
69	0.53948509745932\\
70	0.534726054609468\\
71	0.530008993436131\\
72	0.525333543599891\\
73	0.520699338028261\\
74	0.516106012886869\\
75	0.511553207550887\\
76	0.507040564576726\\
77	0.502567729673967\\
78	0.498134351677549\\
79	0.493740082520196\\
80	0.489384577205092\\
81	0.485067493778794\\
82	0.480788493304387\\
83	0.47654723983487\\
84	0.472343400386786\\
85	0.468176644914072\\
86	0.464046646282155\\
87	0.459953080242264\\
88	0.45589562540597\\
89	0.451873963219961\\
90	0.447887777941027\\
91	0.443936756611272\\
92	0.440020589033543\\
93	0.436138967747079\\
94	0.432291588003368\\
95	0.428478147742224\\
96	0.424698347568068\\
97	0.42095189072643\\
98	0.417238483080643\\
99	0.413557833088753\\
100	0.409909651780631\\
101	0.406293652735283\\
};
\addlegendentry{KM}

\addplot
  table[row sep=crcr]{%
1	1\\
2	0.976793142354244\\
3	0.594133214006442\\
4	0.361380788500269\\
5	0.219809415158464\\
6	0.133698803401965\\
7	0.0813221305765452\\
8	0.049464098347676\\
9	0.0300858881313636\\
10	0.0240648196015739\\
11	0.0237834784800043\\
12	0.0235767779992842\\
13	0.0143405458957974\\
14	0.0087226192058862\\
15	0.00530552228371978\\
16	0.00322707733074233\\
17	0.00196286576819047\\
18	0.00119391086810326\\
19	0.000726185873055654\\
20	0.000442108035117088\\
21	0.000437810511288956\\
22	0.000266297698989762\\
23	0.000161975244218013\\
24	9.85212408481344e-05\\
25	5.99254222174732e-05\\
26	3.64495634592027e-05\\
27	2.21704049308539e-05\\
28	1.34849785442791e-05\\
29	8.2086742613249e-06\\
30	8.12995074224795e-06\\
31	4.94503242781172e-06\\
32	3.00780982411566e-06\\
33	1.82949658472784e-06\\
34	1.11278902230023e-06\\
};
\addlegendentry{Line search}

\addplot
  table[row sep=crcr]{%
1	1\\
2	0.976793142354244\\
3	0.96819335493068\\
4	3.34510993414479e-12\\
};
\addlegendentry{Broyden}

\addplot
  table[row sep=crcr]{%
1	1\\
2	1.4830618466841\\
3	1.03623442509679\\
4	0.461632689240763\\
5	2.22044604925031e-16\\
};
\addlegendentry{Newton}

\end{axis}

\begin{axis}[%
	title={Fixed-point residual \(R\)},
]
\addplot
  table[row sep=crcr]{%
1	0.287973174113417\\
2	0.00861673555358389\\
3	0.00854072341920133\\
4	0.00846538182234893\\
5	0.00839070484791034\\
6	0.00831668663295008\\
7	0.00824332136625145\\
8	0.00817060328786176\\
9	0.00809852668863924\\
10	0.00802708590980503\\
11	0.00795627534249908\\
12	0.00788608942733953\\
13	0.00781652265398612\\
14	0.00774756956070816\\
15	0.00767922473395523\\
16	0.00761148280793234\\
17	0.00754433846417856\\
18	0.00747778643114948\\
19	0.00741182148380345\\
20	0.00734643844319119\\
21	0.00728163217604937\\
22	0.00721739759439759\\
23	0.00715372965513855\\
24	0.00709062335966281\\
25	0.00702807375345558\\
26	0.00696607592570794\\
27	0.006904625008932\\
28	0.00684371617857806\\
29	0.00678334465265558\\
30	0.00672350569135843\\
31	0.00666419459669201\\
32	0.00660540671210586\\
33	0.00654713742212625\\
34	0.00648938215199448\\
35	0.00643213636730861\\
36	0.00637539557366616\\
37	0.00631915531631196\\
38	0.00626341117978836\\
39	0.00620815878758853\\
40	0.00615339380181219\\
41	0.00609911192282672\\
42	0.00604530888892731\\
43	0.00599198047600446\\
44	0.00593912249721128\\
45	0.00588673080263455\\
46	0.00583480127896979\\
47	0.00578332984919783\\
48	0.005732312472264\\
49	0.00568174514276278\\
50	0.00563162389062114\\
51	0.00558194478078866\\
52	0.0055327039129271\\
53	0.0054838974211052\\
54	0.00543552147349497\\
55	0.00538757227207021\\
56	0.00534004605230968\\
57	0.0052929390828997\\
58	0.00524624766544295\\
59	0.00519996813416696\\
60	0.00515409685563706\\
61	0.00510863022847007\\
62	0.00506356468305301\\
63	0.0050188966812619\\
64	0.00497462271618356\\
65	0.00493073931184182\\
66	0.004887243022923\\
67	0.00484413043450628\\
68	0.0048013981617954\\
69	0.004759042849853\\
70	0.00471706117333699\\
71	0.00467544983623981\\
72	0.00463420557162932\\
73	0.00459332514139258\\
74	0.00455280533598163\\
75	0.00451264297416114\\
76	0.00447283490275874\\
77	0.00443337799641831\\
78	0.00439426915735325\\
79	0.00435550531510415\\
80	0.00431708342629764\\
81	0.0042790004744072\\
82	0.00424125346951652\\
83	0.0042038394480846\\
84	0.00416675547271346\\
85	0.0041299986319167\\
86	0.00409356603989202\\
87	0.00405745483629384\\
88	0.00402166218600897\\
89	0.00398618527893411\\
90	0.003951021329755\\
91	0.0039161675777283\\
92	0.00388162128646388\\
93	0.00384737974371103\\
94	0.0038134402611448\\
95	0.00377980017415519\\
96	0.00374645684163793\\
97	0.00371340764578719\\
98	0.00368064999188992\\
99	0.00364818130812241\\
100	0.00361599904534786\\
101	0.00358410067691677\\
};
\addlegendentry{KM}

\addplot
  table[row sep=crcr]{%
1	0.287973174113417\\
2	0.00861673555358389\\
3	0.00524111868389591\\
4	0.00318790392114241\\
5	0.00193903859509607\\
6	0.00117941782109733\\
7	0.000717379110385476\\
8	0.000436358538756519\\
9	0.000264897242019169\\
10	0.000261597628716851\\
11	0.000209366119955775\\
12	0.000207981457297854\\
13	0.000126504462736394\\
14	7.6946182126788e-05\\
15	4.68024195786965e-05\\
16	2.84675136101302e-05\\
17	1.7315325007602e-05\\
18	1.05322309219357e-05\\
19	6.3979642146503e-06\\
20	5.89432779119163e-06\\
21	3.86212518779957e-06\\
22	2.34913284218214e-06\\
23	1.42885712965875e-06\\
24	8.69100572015525e-07\\
25	5.28629343779747e-07\\
26	3.21538065138756e-07\\
27	1.95578420071368e-07\\
28	1.18826936725607e-07\\
29	1.01754859883765e-07\\
30	7.17179846705025e-08\\
31	4.36223749806963e-08\\
32	2.65332553291637e-08\\
33	1.61388195554128e-08\\
34	9.81641689396417e-09\\
};
\addlegendentry{Line search}

\addplot
  table[row sep=crcr]{%
1	0.287973174113417\\
2	0.00861673555358389\\
3	0.00854087292634684\\
4	3.00649896771539e-14\\
};
\addlegendentry{Broyden}

\addplot
  table[row sep=crcr]{%
1	0.287973174113417\\
2	0.0977579646317743\\
3	0.0435502537019617\\
4	0.461632689240763\\
5	2.22044604925031e-16\\
};
\addlegendentry{Newton}

\end{axis}
\end{tikzpicture}%
%
	}}
			\caption[Alternating projections on polyhedral cones]{}%
			\label{fig:Cones}%
			\centering%
			\vspace*{15pt}%
			{{%
		\pgfkeys{/pgf/images/include external/.code={\includegraphics[width=\linewidth,height=\myHeight,keepaspectratio]{#width=\linewidth,height=\myHeight,keepaspectratio}}}%
		\tikzsetnextfilename{Cones_graph}%
		\input{./TeX/Tikz/Cones_graph.tex}%
	}}
		\end{subfigure}
		\hfill
		\begin{subfigure}[t]{0.31\textwidth}%
			{{%
		\pgfkeys{/pgf/images/include external/.code={\includegraphics[width=\linewidth]{#width=\linewidth}}}%
		\tikzsetnextfilename{LineCircleH}%
		\input{./TeX/Tikz/LineCircleH.tex}%
	}}
			\caption[Alternating projections on ball and tangent line]{}%
			\label{fig:LineCircle}%
			\centering%
			{{%
		\pgfkeys{/pgf/images/include external/.code={\includegraphics[width=\linewidth,height=\myHeight,keepaspectratio]{#width=\linewidth,height=\myHeight,keepaspectratio}}}%
		\tikzsetnextfilename{LineCircleH_graph}%
		\input{./TeX/Tikz/LineCircleH_graph.tex}%
	}}
		\end{subfigure}%
		\hfill
		\begin{subfigure}[t]{0.31\textwidth}%
			{{%
		\pgfkeys{/pgf/images/include external/.code={\includegraphics[width=\linewidth]{#width=\linewidth}}}%
		\tikzsetnextfilename{IceCream}%
		\input{./TeX/Tikz/IceCream.tex}%
	}}
			\caption[Alternating projections on polyhedral cones]{}%
			\label{fig:IceCream}%
			\centering%
			\vspace*{8pt}%
			~~{{%
		\pgfkeys{/pgf/images/include external/.code={\includegraphics[width=\linewidth,height=\myHeight,keepaspectratio]{#width=\linewidth,height=\myHeight,keepaspectratio}}}%
		\tikzsetnextfilename{IceCream_graph}%
%
%
\begin{tikzpicture}

\begin{axis}[%
	hide axis,
	scale only axis,
	xmin=-10,
	xmax=10,
	ymin=-10,
	ymax=10,
	zmin=0,
	zmax=1,
	view={-37.5}{30},
	axis background/.style={fill=white},
]

\node[scale=2.5] at (axis cs:7,-4,1) {\(C_1\)};
\node[scale=2.5] at (axis cs:-10,0,0.6) {\(C_2\)};
\path[draw, black, thick] (axis cs:0,0,0) -- node [above,pos=0.75,sloped,inner sep=0pt,scale=2.5] {\(C_1{\cap}C_2\)} (axis cs:0,10,1);

\addplot3[%
surf,faceted color=red,opacity=0,
fill opacity=0.25, fill=red, z buffer=sort, colormap={mymap}{[1pt] rgb(0pt)=(0.2081,0.1663,0.5292); rgb(1pt)=(0.211624,0.189781,0.577676); rgb(2pt)=(0.212252,0.213771,0.626971); rgb(3pt)=(0.2081,0.2386,0.677086); rgb(4pt)=(0.195905,0.264457,0.7279); rgb(5pt)=(0.170729,0.291938,0.779248); rgb(6pt)=(0.125271,0.324243,0.830271); rgb(7pt)=(0.0591333,0.359833,0.868333); rgb(8pt)=(0.0116952,0.38751,0.881957); rgb(9pt)=(0.00595714,0.408614,0.882843); rgb(10pt)=(0.0165143,0.4266,0.878633); rgb(11pt)=(0.0328524,0.443043,0.871957); rgb(12pt)=(0.0498143,0.458571,0.864057); rgb(13pt)=(0.0629333,0.47369,0.855438); rgb(14pt)=(0.0722667,0.488667,0.8467); rgb(15pt)=(0.0779429,0.503986,0.838371); rgb(16pt)=(0.0793476,0.520024,0.831181); rgb(17pt)=(0.0749429,0.537543,0.826271); rgb(18pt)=(0.0640571,0.556986,0.823957); rgb(19pt)=(0.0487714,0.577224,0.822829); rgb(20pt)=(0.0343429,0.596581,0.819852); rgb(21pt)=(0.0265,0.6137,0.8135); rgb(22pt)=(0.0238905,0.628662,0.803762); rgb(23pt)=(0.0230905,0.641786,0.791267); rgb(24pt)=(0.0227714,0.653486,0.776757); rgb(25pt)=(0.0266619,0.664195,0.760719); rgb(26pt)=(0.0383714,0.674271,0.743552); rgb(27pt)=(0.0589714,0.683757,0.725386); rgb(28pt)=(0.0843,0.692833,0.706167); rgb(29pt)=(0.113295,0.7015,0.685857); rgb(30pt)=(0.145271,0.709757,0.664629); rgb(31pt)=(0.180133,0.717657,0.642433); rgb(32pt)=(0.217829,0.725043,0.619262); rgb(33pt)=(0.258643,0.731714,0.595429); rgb(34pt)=(0.302171,0.737605,0.571186); rgb(35pt)=(0.348167,0.742433,0.547267); rgb(36pt)=(0.395257,0.7459,0.524443); rgb(37pt)=(0.44201,0.748081,0.503314); rgb(38pt)=(0.487124,0.749062,0.483976); rgb(39pt)=(0.530029,0.749114,0.466114); rgb(40pt)=(0.570857,0.748519,0.44939); rgb(41pt)=(0.609852,0.747314,0.433686); rgb(42pt)=(0.6473,0.7456,0.4188); rgb(43pt)=(0.683419,0.743476,0.404433); rgb(44pt)=(0.71841,0.741133,0.390476); rgb(45pt)=(0.752486,0.7384,0.376814); rgb(46pt)=(0.785843,0.735567,0.363271); rgb(47pt)=(0.818505,0.732733,0.34979); rgb(48pt)=(0.850657,0.7299,0.336029); rgb(49pt)=(0.882433,0.727433,0.3217); rgb(50pt)=(0.913933,0.725786,0.306276); rgb(51pt)=(0.944957,0.726114,0.288643); rgb(52pt)=(0.973895,0.731395,0.266648); rgb(53pt)=(0.993771,0.745457,0.240348); rgb(54pt)=(0.999043,0.765314,0.216414); rgb(55pt)=(0.995533,0.786057,0.196652); rgb(56pt)=(0.988,0.8066,0.179367); rgb(57pt)=(0.978857,0.827143,0.163314); rgb(58pt)=(0.9697,0.848138,0.147452); rgb(59pt)=(0.962586,0.870514,0.1309); rgb(60pt)=(0.958871,0.8949,0.113243); rgb(61pt)=(0.959824,0.921833,0.0948381); rgb(62pt)=(0.9661,0.951443,0.0755333); rgb(63pt)=(0.9763,0.9831,0.0538)}, mesh/rows=21]
table[row sep=crcr, point meta=\thisrow{c}] {%
x	y	z	c\\
0	0	0	0\\
10	0	1	1\\
0	0	0	0\\
9.51056516295153	3.09016994374947	1	1\\
0	0	0	0\\
8.09016994374947	5.87785252292473	1	1\\
0	0	0	0\\
5.87785252292473	8.09016994374947	1	1\\
0	0	0	0\\
3.09016994374947	9.51056516295153	1	1\\
0	0	0	0\\
6.12323399573677e-16	10	1	1\\
-0	0	0	0\\
-3.09016994374947	9.51056516295154	1	1\\
-0	0	0	0\\
-5.87785252292473	8.09016994374947	1	1\\
-0	0	0	0\\
-8.09016994374947	5.87785252292473	1	1\\
-0	0	0	0\\
-9.51056516295153	3.09016994374947	1	1\\
-0	0	0	0\\
-10	1.22464679914735e-15	1	1\\
-0	-0	0	0\\
-9.51056516295153	-3.09016994374948	1	1\\
-0	-0	0	0\\
-8.09016994374947	-5.87785252292473	1	1\\
-0	-0	0	0\\
-5.87785252292473	-8.09016994374947	1	1\\
-0	-0	0	0\\
-3.09016994374948	-9.51056516295153	1	1\\
-0	-0	0	0\\
-1.83697019872103e-15	-10	1	1\\
0	-0	0	0\\
3.09016994374947	-9.51056516295154	1	1\\
0	-0	0	0\\
5.87785252292473	-8.09016994374948	1	1\\
0	-0	0	0\\
8.09016994374947	-5.87785252292473	1	1\\
0	-0	0	0\\
9.51056516295153	-3.09016994374948	1	1\\
0	0	0	0\\
10	0	1	1\\
};

\addplot3[%
surf,opacity=0,
fill opacity=0.25, fill=black, z buffer=sort, colormap={mymap}{[1pt] rgb(0pt)=(0.2081,0.1663,0.5292); rgb(1pt)=(0.211624,0.189781,0.577676); rgb(2pt)=(0.212252,0.213771,0.626971); rgb(3pt)=(0.2081,0.2386,0.677086); rgb(4pt)=(0.195905,0.264457,0.7279); rgb(5pt)=(0.170729,0.291938,0.779248); rgb(6pt)=(0.125271,0.324243,0.830271); rgb(7pt)=(0.0591333,0.359833,0.868333); rgb(8pt)=(0.0116952,0.38751,0.881957); rgb(9pt)=(0.00595714,0.408614,0.882843); rgb(10pt)=(0.0165143,0.4266,0.878633); rgb(11pt)=(0.0328524,0.443043,0.871957); rgb(12pt)=(0.0498143,0.458571,0.864057); rgb(13pt)=(0.0629333,0.47369,0.855438); rgb(14pt)=(0.0722667,0.488667,0.8467); rgb(15pt)=(0.0779429,0.503986,0.838371); rgb(16pt)=(0.0793476,0.520024,0.831181); rgb(17pt)=(0.0749429,0.537543,0.826271); rgb(18pt)=(0.0640571,0.556986,0.823957); rgb(19pt)=(0.0487714,0.577224,0.822829); rgb(20pt)=(0.0343429,0.596581,0.819852); rgb(21pt)=(0.0265,0.6137,0.8135); rgb(22pt)=(0.0238905,0.628662,0.803762); rgb(23pt)=(0.0230905,0.641786,0.791267); rgb(24pt)=(0.0227714,0.653486,0.776757); rgb(25pt)=(0.0266619,0.664195,0.760719); rgb(26pt)=(0.0383714,0.674271,0.743552); rgb(27pt)=(0.0589714,0.683757,0.725386); rgb(28pt)=(0.0843,0.692833,0.706167); rgb(29pt)=(0.113295,0.7015,0.685857); rgb(30pt)=(0.145271,0.709757,0.664629); rgb(31pt)=(0.180133,0.717657,0.642433); rgb(32pt)=(0.217829,0.725043,0.619262); rgb(33pt)=(0.258643,0.731714,0.595429); rgb(34pt)=(0.302171,0.737605,0.571186); rgb(35pt)=(0.348167,0.742433,0.547267); rgb(36pt)=(0.395257,0.7459,0.524443); rgb(37pt)=(0.44201,0.748081,0.503314); rgb(38pt)=(0.487124,0.749062,0.483976); rgb(39pt)=(0.530029,0.749114,0.466114); rgb(40pt)=(0.570857,0.748519,0.44939); rgb(41pt)=(0.609852,0.747314,0.433686); rgb(42pt)=(0.6473,0.7456,0.4188); rgb(43pt)=(0.683419,0.743476,0.404433); rgb(44pt)=(0.71841,0.741133,0.390476); rgb(45pt)=(0.752486,0.7384,0.376814); rgb(46pt)=(0.785843,0.735567,0.363271); rgb(47pt)=(0.818505,0.732733,0.34979); rgb(48pt)=(0.850657,0.7299,0.336029); rgb(49pt)=(0.882433,0.727433,0.3217); rgb(50pt)=(0.913933,0.725786,0.306276); rgb(51pt)=(0.944957,0.726114,0.288643); rgb(52pt)=(0.973895,0.731395,0.266648); rgb(53pt)=(0.993771,0.745457,0.240348); rgb(54pt)=(0.999043,0.765314,0.216414); rgb(55pt)=(0.995533,0.786057,0.196652); rgb(56pt)=(0.988,0.8066,0.179367); rgb(57pt)=(0.978857,0.827143,0.163314); rgb(58pt)=(0.9697,0.848138,0.147452); rgb(59pt)=(0.962586,0.870514,0.1309); rgb(60pt)=(0.958871,0.8949,0.113243); rgb(61pt)=(0.959824,0.921833,0.0948381); rgb(62pt)=(0.9661,0.951443,0.0755333); rgb(63pt)=(0.9763,0.9831,0.0538)}, mesh/rows=2]
table[row sep=crcr, point meta=\thisrow{c}] {%
x	y	z	c\\
-10	0	0	0\\
-10	10	1	1\\
10	0	0	0\\
10	10	1	1\\
};

\end{axis}
\end{tikzpicture}%
%
	}}
		\end{subfigure}%
		\caption[]{%
			\emph{\textbf{(a)} Alternating projections on polyhedral cones.}
				\(R=\id-\proj_{C_2}\circ\proj_{C_1}\) is globally metrically subregular, however the \(Q\)-linear convergence of the KM scheme is very slow.
		
			\emph{\textbf{(b)} Alternating projections on ball and tangent line.}
			With or without line search the KM scheme is not linearly convergent due to the fact that the residual \(R\) is not metrically subregular at \(x_\star\).
		
			\emph{\textbf{(c)} Alternating projections on second-order cone and tangent plane}.
			In contrast with the slow sublinear rate of KM both with and without line search, and despite the non isolatedness of any solution, Broyden's scheme exhibits an appealing linear convergence rate.
		}%
	\end{figure*}%
\fi
Illustrative examples can be easily constructed for the problem of finding a point in the intersection of two closed convex sets \(C_1\) and \(C_2\) with \(C_1\cap C_2\neq\emptyset\).
The problem can be solved by means of fixed-point iterations of the (nonexpansive) \DEF{alternating projections} operator \(T=\proj_{C_2}\circ\proj_{C_1}\).

In \Cref{fig:Cones} we consider the case of two polyhedral cones, namely
\ifieee
	\(
		C_1
	{}={}
		\set{x\in\R^2}[0.1x_1\leq x_2\leq 0.2x_1]
	\)
	and
	\(
		C_2
	{}={}
		\set{x\in\R^2}[0.3x_1\leq x_2\leq 0.35x_1]
	\).
\else
	\[
		C_1
	{}={}
		\set{x\in\R^2}[0.1x_1\leq x_2\leq 0.2x_1]
	~~\text{and}~~
		C_2
	{}={}
		\set{x\in\R^2}[0.3x_1\leq x_2\leq 0.35x_1].
	\]
\fi
Alternating projections is then linearly convergent (to the unique intersection point \(0\)) due to the fact that \(R=\id-T\) is piecewise affine and hence globally metrically subregular.
However, the convergence is extremely slow due to the pathological small angle between the two cones, as it is apparent in \Cref{fig:Cones}.
\ifieee\else
	\begin{figure}
		\begin{minipage}[t][][t]{0.40\linewidth}\vspace*{5pt}
			{{%
		\pgfkeys{/pgf/images/include external/.code={\includegraphics[width=\linewidth]{#width=\linewidth}}}%
		\tikzsetnextfilename{Cones_graph}%
		\input{./TeX/Tikz/Cones_graph.tex}%
	}}
			\caption[Alternating projections on polyhedral cones]{%
				\emph{Alternating projections on polyhedral cones.}
			
				\(R=\id-\proj_{C_2}\circ\proj_{C_1}\) is globally metrically subregular, however the \(Q\)-linear convergence of the KM scheme is very slow.%
			}%
			\label{fig:Cones}
		\end{minipage}
		\hfill
		\begin{minipage}[t][][t]{0.59\linewidth}\vspace*{0pt}
			{{%
		\pgfkeys{/pgf/images/include external/.code={\includegraphics[width=\linewidth]{#width=\linewidth}}}%
		\tikzsetnextfilename{Cones}%
%
%
\begin{tikzpicture}
\pgfplotsset{%
	every axis/.append style={
		plotaxis,
		xmax=40,
		ymin=1e-12,
		ymax=3,
	},
	every axis plot/.append style={
		plotstyle,
	},
}

\begin{axis}[%
	title={Distance from solution},
]
\addplot
  table[row sep=crcr]{%
1	1\\
2	0.976793142354244\\
3	0.96817640680066\\
4	0.959635683381458\\
5	0.951170301559109\\
6	0.942779596711199\\
7	0.934462910078249\\
8	0.926219588711997\\
9	0.918048985424136\\
10	0.909950458735496\\
11	0.901923372825691\\
12	0.893967097483192\\
13	0.886081008055853\\
14	0.878264485401867\\
15	0.870516915841159\\
16	0.862837691107203\\
17	0.855226208299271\\
18	0.847681869835092\\
19	0.840204083403943\\
20	0.832792261920139\\
21	0.825445823476948\\
22	0.818164191300899\\
23	0.810946793706501\\
24	0.803793064051363\\
25	0.7967024406917\\
26	0.789674366938244\\
27	0.782708291012536\\
28	0.775803666003604\\
29	0.768959949825026\\
30	0.762176605172371\\
31	0.755453099481012\\
32	0.74878890488432\\
33	0.742183498172215\\
34	0.735636360750088\\
35	0.729146978598094\\
36	0.722714842230785\\
37	0.716339446657119\\
38	0.710020291340807\\
39	0.703756880161019\\
40	0.69754872137343\\
41	0.691395327571618\\
42	0.685296215648791\\
43	0.679250906759864\\
44	0.673258926283859\\
45	0.667319803786648\\
46	0.661433072984014\\
47	0.655598271705044\\
48	0.649814941855846\\
49	0.644082629383582\\
50	0.638400884240819\\
51	0.632769260350198\\
52	0.627187315569409\\
53	0.621654611656482\\
54	0.616170714235377\\
55	0.610735192761882\\
56	0.605347620489812\\
57	0.600007574437502\\
58	0.594714635354603\\
59	0.58946838768916\\
60	0.584268419554993\\
61	0.579114322699356\\
62	0.574005692470886\\
63	0.568942127787832\\
64	0.563923231106571\\
65	0.558948608390387\\
66	0.554017869078545\\
67	0.549130626055622\\
68	0.544286495621116\\
69	0.53948509745932\\
70	0.534726054609468\\
71	0.530008993436131\\
72	0.525333543599891\\
73	0.520699338028261\\
74	0.516106012886869\\
75	0.511553207550887\\
76	0.507040564576726\\
77	0.502567729673967\\
78	0.498134351677549\\
79	0.493740082520196\\
80	0.489384577205092\\
81	0.485067493778794\\
82	0.480788493304387\\
83	0.47654723983487\\
84	0.472343400386786\\
85	0.468176644914072\\
86	0.464046646282155\\
87	0.459953080242264\\
88	0.45589562540597\\
89	0.451873963219961\\
90	0.447887777941027\\
91	0.443936756611272\\
92	0.440020589033543\\
93	0.436138967747079\\
94	0.432291588003368\\
95	0.428478147742224\\
96	0.424698347568068\\
97	0.42095189072643\\
98	0.417238483080643\\
99	0.413557833088753\\
100	0.409909651780631\\
101	0.406293652735283\\
};
\addlegendentry{KM}

\addplot
  table[row sep=crcr]{%
1	1\\
2	0.976793142354244\\
3	0.594133214006442\\
4	0.361380788500269\\
5	0.219809415158464\\
6	0.133698803401965\\
7	0.0813221305765452\\
8	0.049464098347676\\
9	0.0300858881313636\\
10	0.0240648196015739\\
11	0.0237834784800043\\
12	0.0235767779992842\\
13	0.0143405458957974\\
14	0.0087226192058862\\
15	0.00530552228371978\\
16	0.00322707733074233\\
17	0.00196286576819047\\
18	0.00119391086810326\\
19	0.000726185873055654\\
20	0.000442108035117088\\
21	0.000437810511288956\\
22	0.000266297698989762\\
23	0.000161975244218013\\
24	9.85212408481344e-05\\
25	5.99254222174732e-05\\
26	3.64495634592027e-05\\
27	2.21704049308539e-05\\
28	1.34849785442791e-05\\
29	8.2086742613249e-06\\
30	8.12995074224795e-06\\
31	4.94503242781172e-06\\
32	3.00780982411566e-06\\
33	1.82949658472784e-06\\
34	1.11278902230023e-06\\
};
\addlegendentry{Line search}

\addplot
  table[row sep=crcr]{%
1	1\\
2	0.976793142354244\\
3	0.96819335493068\\
4	3.34510993414479e-12\\
};
\addlegendentry{Broyden}

\addplot
  table[row sep=crcr]{%
1	1\\
2	1.4830618466841\\
3	1.03623442509679\\
4	0.461632689240763\\
5	2.22044604925031e-16\\
};
\addlegendentry{Newton}

\end{axis}

\begin{axis}[%
	title={Fixed-point residual \(R\)},
]
\addplot
  table[row sep=crcr]{%
1	0.287973174113417\\
2	0.00861673555358389\\
3	0.00854072341920133\\
4	0.00846538182234893\\
5	0.00839070484791034\\
6	0.00831668663295008\\
7	0.00824332136625145\\
8	0.00817060328786176\\
9	0.00809852668863924\\
10	0.00802708590980503\\
11	0.00795627534249908\\
12	0.00788608942733953\\
13	0.00781652265398612\\
14	0.00774756956070816\\
15	0.00767922473395523\\
16	0.00761148280793234\\
17	0.00754433846417856\\
18	0.00747778643114948\\
19	0.00741182148380345\\
20	0.00734643844319119\\
21	0.00728163217604937\\
22	0.00721739759439759\\
23	0.00715372965513855\\
24	0.00709062335966281\\
25	0.00702807375345558\\
26	0.00696607592570794\\
27	0.006904625008932\\
28	0.00684371617857806\\
29	0.00678334465265558\\
30	0.00672350569135843\\
31	0.00666419459669201\\
32	0.00660540671210586\\
33	0.00654713742212625\\
34	0.00648938215199448\\
35	0.00643213636730861\\
36	0.00637539557366616\\
37	0.00631915531631196\\
38	0.00626341117978836\\
39	0.00620815878758853\\
40	0.00615339380181219\\
41	0.00609911192282672\\
42	0.00604530888892731\\
43	0.00599198047600446\\
44	0.00593912249721128\\
45	0.00588673080263455\\
46	0.00583480127896979\\
47	0.00578332984919783\\
48	0.005732312472264\\
49	0.00568174514276278\\
50	0.00563162389062114\\
51	0.00558194478078866\\
52	0.0055327039129271\\
53	0.0054838974211052\\
54	0.00543552147349497\\
55	0.00538757227207021\\
56	0.00534004605230968\\
57	0.0052929390828997\\
58	0.00524624766544295\\
59	0.00519996813416696\\
60	0.00515409685563706\\
61	0.00510863022847007\\
62	0.00506356468305301\\
63	0.0050188966812619\\
64	0.00497462271618356\\
65	0.00493073931184182\\
66	0.004887243022923\\
67	0.00484413043450628\\
68	0.0048013981617954\\
69	0.004759042849853\\
70	0.00471706117333699\\
71	0.00467544983623981\\
72	0.00463420557162932\\
73	0.00459332514139258\\
74	0.00455280533598163\\
75	0.00451264297416114\\
76	0.00447283490275874\\
77	0.00443337799641831\\
78	0.00439426915735325\\
79	0.00435550531510415\\
80	0.00431708342629764\\
81	0.0042790004744072\\
82	0.00424125346951652\\
83	0.0042038394480846\\
84	0.00416675547271346\\
85	0.0041299986319167\\
86	0.00409356603989202\\
87	0.00405745483629384\\
88	0.00402166218600897\\
89	0.00398618527893411\\
90	0.003951021329755\\
91	0.0039161675777283\\
92	0.00388162128646388\\
93	0.00384737974371103\\
94	0.0038134402611448\\
95	0.00377980017415519\\
96	0.00374645684163793\\
97	0.00371340764578719\\
98	0.00368064999188992\\
99	0.00364818130812241\\
100	0.00361599904534786\\
101	0.00358410067691677\\
};
\addlegendentry{KM}

\addplot
  table[row sep=crcr]{%
1	0.287973174113417\\
2	0.00861673555358389\\
3	0.00524111868389591\\
4	0.00318790392114241\\
5	0.00193903859509607\\
6	0.00117941782109733\\
7	0.000717379110385476\\
8	0.000436358538756519\\
9	0.000264897242019169\\
10	0.000261597628716851\\
11	0.000209366119955775\\
12	0.000207981457297854\\
13	0.000126504462736394\\
14	7.6946182126788e-05\\
15	4.68024195786965e-05\\
16	2.84675136101302e-05\\
17	1.7315325007602e-05\\
18	1.05322309219357e-05\\
19	6.3979642146503e-06\\
20	5.89432779119163e-06\\
21	3.86212518779957e-06\\
22	2.34913284218214e-06\\
23	1.42885712965875e-06\\
24	8.69100572015525e-07\\
25	5.28629343779747e-07\\
26	3.21538065138756e-07\\
27	1.95578420071368e-07\\
28	1.18826936725607e-07\\
29	1.01754859883765e-07\\
30	7.17179846705025e-08\\
31	4.36223749806963e-08\\
32	2.65332553291637e-08\\
33	1.61388195554128e-08\\
34	9.81641689396417e-09\\
};
\addlegendentry{Line search}

\addplot
  table[row sep=crcr]{%
1	0.287973174113417\\
2	0.00861673555358389\\
3	0.00854087292634684\\
4	3.00649896771539e-14\\
};
\addlegendentry{Broyden}

\addplot
  table[row sep=crcr]{%
1	0.287973174113417\\
2	0.0977579646317743\\
3	0.0435502537019617\\
4	0.461632689240763\\
5	2.22044604925031e-16\\
};
\addlegendentry{Newton}

\end{axis}
\end{tikzpicture}%
%
	}}%
		\end{minipage}
	\end{figure}
\fi

As an attempt to overcome this frequent phenomenon, \cite{giselsson2016line} proposes a foretracking line-search heuristic which is particularly effective when subsequent fixed-point iterations proceed along almost parallel directions.
Iteration-wise, in such instances the line search does yield a considerable improvement upon the plain KM scheme; however, each foretrack prescribes extra evaluations of \(T\) and unless \(T\) has a specific structure the computational overhead might outweight the advantages.
Moreover,
its asymptotic convergence rates do not improve upon the plain KM scheme.
\Cref{fig:LineCircle} illustrates this fact relative to
\ifieee
	\(
		C_1
	{}={}
		\set{x\in\R^2}[x_1^2+x_2^2\leq 1]
	\)
	and
	\(
		C_2
	{}={}
		\set{x\in\R^2}[x_1=1]
	\).
\else
	\[
		C_1
	{}={}
		\set{x\in\R^2}[x_1^2+x_2^2\leq 1]
	\quad\text{and}\quad
		C_2
	{}={}
		\set{x\in\R^2}[x_1=1].
	\]
\fi
\ifieee\else%
	\begin{figure}
		\begin{minipage}[t][][t]{0.40\linewidth}\vspace*{7pt}
			{{%
		\pgfkeys{/pgf/images/include external/.code={\includegraphics[width=\linewidth]{#width=\linewidth}}}%
		\tikzsetnextfilename{LineCircle_graph}%
		\input{./TeX/Tikz/LineCircle_graph.tex}%
	}}
		\end{minipage}
		\hfill
		\begin{minipage}[t][][t]{0.59\linewidth}\vspace*{0pt}
			{{%
		\pgfkeys{/pgf/images/include external/.code={\includegraphics[width=\linewidth]{#width=\linewidth}}}%
		\tikzsetnextfilename{LineCircle}%
		\input{./TeX/Tikz/LineCircle.tex}%
	}}%
			\caption[Alternating projections on ball and tangent line]{%
				\emph{Alternating projections on ball and tangent line.}
			
				With or without line search the KM scheme is not linearly convergent due to the fact that the residual \(R\) is not metrically subregular at \(x_\star\).%
			}%
			\label{fig:LineCircle}
		\end{minipage}
	\end{figure}%
\fi
Despite a good performance on early iterations, the line search cannot improve the asymptotic sublinear rate of the plain KM scheme due to the fact that the residual is not metrically subregular at the (unique) solution
\ifieee
	\(x_\star=(0,1)\).
\else
	\(x_\star=(1,0)\).
\fi
In particular, it is evident that medium-to-high accuracy cannot be achieved in a reasonable number of iterations with either methods.

In response to this limitation there comes the need to include some ``first-order-like information''.
Specifically, the problem of finding a fixed point of \(T\) can be rephrased in terms of solving the system of nonlinear (monotone) equations \(Rx=0\), which could \emph{possibly} be solved efficiently with Newton-type methods.
In the toy simulations of this section, the purple lines correspond to the semismooth Newton iterations
\[
	x^+=x-G^{-1}Rx
\quad\text{for some}\quad
	G\in\partial Rx,
\]
where \(\partial R\) is the \DEF{Clarke generalized Jacobian} of \(R\) \cite[Def. 7.1.1]{facchinei2003finite}.
Interestingly, in the proposed simulations this method exhibits fast convergence
even when the limit point is a non isolated solution, as in the case of the second-order cone
\(
	C_1
{}={}
	\set{x\in\R^3}[x_3\geq 0.1\sqrt{x_1^2+x_2^2}]
\)
and the tangent plane
\(
	C_2
{}={}
	\set{x\in\R^3}[x_3=0.1x_2]
\)
considered in \Cref{fig:IceCream}.
\ifieee\else
	\begin{figure}
		\begin{minipage}[t][][t]{0.46\linewidth}\vspace*{0pt}
			{{%
		\pgfkeys{/pgf/images/include external/.code={\includegraphics[width=\linewidth]{#width=\linewidth}}}%
		\tikzsetnextfilename{IceCream_graph}%
%
%
\begin{tikzpicture}

\begin{axis}[%
	hide axis,
	scale only axis,
	xmin=-10,
	xmax=10,
	ymin=-10,
	ymax=10,
	zmin=0,
	zmax=1,
	view={-37.5}{30},
	axis background/.style={fill=white},
]

\node[scale=2.5] at (axis cs:7,-4,1) {\(C_1\)};
\node[scale=2.5] at (axis cs:-10,0,0.6) {\(C_2\)};
\path[draw, black, thick] (axis cs:0,0,0) -- node [above,pos=0.75,sloped,inner sep=0pt,scale=2.5] {\(C_1{\cap}C_2\)} (axis cs:0,10,1);

\addplot3[%
surf,faceted color=red,opacity=0,
fill opacity=0.25, fill=red, z buffer=sort, colormap={mymap}{[1pt] rgb(0pt)=(0.2081,0.1663,0.5292); rgb(1pt)=(0.211624,0.189781,0.577676); rgb(2pt)=(0.212252,0.213771,0.626971); rgb(3pt)=(0.2081,0.2386,0.677086); rgb(4pt)=(0.195905,0.264457,0.7279); rgb(5pt)=(0.170729,0.291938,0.779248); rgb(6pt)=(0.125271,0.324243,0.830271); rgb(7pt)=(0.0591333,0.359833,0.868333); rgb(8pt)=(0.0116952,0.38751,0.881957); rgb(9pt)=(0.00595714,0.408614,0.882843); rgb(10pt)=(0.0165143,0.4266,0.878633); rgb(11pt)=(0.0328524,0.443043,0.871957); rgb(12pt)=(0.0498143,0.458571,0.864057); rgb(13pt)=(0.0629333,0.47369,0.855438); rgb(14pt)=(0.0722667,0.488667,0.8467); rgb(15pt)=(0.0779429,0.503986,0.838371); rgb(16pt)=(0.0793476,0.520024,0.831181); rgb(17pt)=(0.0749429,0.537543,0.826271); rgb(18pt)=(0.0640571,0.556986,0.823957); rgb(19pt)=(0.0487714,0.577224,0.822829); rgb(20pt)=(0.0343429,0.596581,0.819852); rgb(21pt)=(0.0265,0.6137,0.8135); rgb(22pt)=(0.0238905,0.628662,0.803762); rgb(23pt)=(0.0230905,0.641786,0.791267); rgb(24pt)=(0.0227714,0.653486,0.776757); rgb(25pt)=(0.0266619,0.664195,0.760719); rgb(26pt)=(0.0383714,0.674271,0.743552); rgb(27pt)=(0.0589714,0.683757,0.725386); rgb(28pt)=(0.0843,0.692833,0.706167); rgb(29pt)=(0.113295,0.7015,0.685857); rgb(30pt)=(0.145271,0.709757,0.664629); rgb(31pt)=(0.180133,0.717657,0.642433); rgb(32pt)=(0.217829,0.725043,0.619262); rgb(33pt)=(0.258643,0.731714,0.595429); rgb(34pt)=(0.302171,0.737605,0.571186); rgb(35pt)=(0.348167,0.742433,0.547267); rgb(36pt)=(0.395257,0.7459,0.524443); rgb(37pt)=(0.44201,0.748081,0.503314); rgb(38pt)=(0.487124,0.749062,0.483976); rgb(39pt)=(0.530029,0.749114,0.466114); rgb(40pt)=(0.570857,0.748519,0.44939); rgb(41pt)=(0.609852,0.747314,0.433686); rgb(42pt)=(0.6473,0.7456,0.4188); rgb(43pt)=(0.683419,0.743476,0.404433); rgb(44pt)=(0.71841,0.741133,0.390476); rgb(45pt)=(0.752486,0.7384,0.376814); rgb(46pt)=(0.785843,0.735567,0.363271); rgb(47pt)=(0.818505,0.732733,0.34979); rgb(48pt)=(0.850657,0.7299,0.336029); rgb(49pt)=(0.882433,0.727433,0.3217); rgb(50pt)=(0.913933,0.725786,0.306276); rgb(51pt)=(0.944957,0.726114,0.288643); rgb(52pt)=(0.973895,0.731395,0.266648); rgb(53pt)=(0.993771,0.745457,0.240348); rgb(54pt)=(0.999043,0.765314,0.216414); rgb(55pt)=(0.995533,0.786057,0.196652); rgb(56pt)=(0.988,0.8066,0.179367); rgb(57pt)=(0.978857,0.827143,0.163314); rgb(58pt)=(0.9697,0.848138,0.147452); rgb(59pt)=(0.962586,0.870514,0.1309); rgb(60pt)=(0.958871,0.8949,0.113243); rgb(61pt)=(0.959824,0.921833,0.0948381); rgb(62pt)=(0.9661,0.951443,0.0755333); rgb(63pt)=(0.9763,0.9831,0.0538)}, mesh/rows=21]
table[row sep=crcr, point meta=\thisrow{c}] {%
x	y	z	c\\
0	0	0	0\\
10	0	1	1\\
0	0	0	0\\
9.51056516295153	3.09016994374947	1	1\\
0	0	0	0\\
8.09016994374947	5.87785252292473	1	1\\
0	0	0	0\\
5.87785252292473	8.09016994374947	1	1\\
0	0	0	0\\
3.09016994374947	9.51056516295153	1	1\\
0	0	0	0\\
6.12323399573677e-16	10	1	1\\
-0	0	0	0\\
-3.09016994374947	9.51056516295154	1	1\\
-0	0	0	0\\
-5.87785252292473	8.09016994374947	1	1\\
-0	0	0	0\\
-8.09016994374947	5.87785252292473	1	1\\
-0	0	0	0\\
-9.51056516295153	3.09016994374947	1	1\\
-0	0	0	0\\
-10	1.22464679914735e-15	1	1\\
-0	-0	0	0\\
-9.51056516295153	-3.09016994374948	1	1\\
-0	-0	0	0\\
-8.09016994374947	-5.87785252292473	1	1\\
-0	-0	0	0\\
-5.87785252292473	-8.09016994374947	1	1\\
-0	-0	0	0\\
-3.09016994374948	-9.51056516295153	1	1\\
-0	-0	0	0\\
-1.83697019872103e-15	-10	1	1\\
0	-0	0	0\\
3.09016994374947	-9.51056516295154	1	1\\
0	-0	0	0\\
5.87785252292473	-8.09016994374948	1	1\\
0	-0	0	0\\
8.09016994374947	-5.87785252292473	1	1\\
0	-0	0	0\\
9.51056516295153	-3.09016994374948	1	1\\
0	0	0	0\\
10	0	1	1\\
};

\addplot3[%
surf,opacity=0,
fill opacity=0.25, fill=black, z buffer=sort, colormap={mymap}{[1pt] rgb(0pt)=(0.2081,0.1663,0.5292); rgb(1pt)=(0.211624,0.189781,0.577676); rgb(2pt)=(0.212252,0.213771,0.626971); rgb(3pt)=(0.2081,0.2386,0.677086); rgb(4pt)=(0.195905,0.264457,0.7279); rgb(5pt)=(0.170729,0.291938,0.779248); rgb(6pt)=(0.125271,0.324243,0.830271); rgb(7pt)=(0.0591333,0.359833,0.868333); rgb(8pt)=(0.0116952,0.38751,0.881957); rgb(9pt)=(0.00595714,0.408614,0.882843); rgb(10pt)=(0.0165143,0.4266,0.878633); rgb(11pt)=(0.0328524,0.443043,0.871957); rgb(12pt)=(0.0498143,0.458571,0.864057); rgb(13pt)=(0.0629333,0.47369,0.855438); rgb(14pt)=(0.0722667,0.488667,0.8467); rgb(15pt)=(0.0779429,0.503986,0.838371); rgb(16pt)=(0.0793476,0.520024,0.831181); rgb(17pt)=(0.0749429,0.537543,0.826271); rgb(18pt)=(0.0640571,0.556986,0.823957); rgb(19pt)=(0.0487714,0.577224,0.822829); rgb(20pt)=(0.0343429,0.596581,0.819852); rgb(21pt)=(0.0265,0.6137,0.8135); rgb(22pt)=(0.0238905,0.628662,0.803762); rgb(23pt)=(0.0230905,0.641786,0.791267); rgb(24pt)=(0.0227714,0.653486,0.776757); rgb(25pt)=(0.0266619,0.664195,0.760719); rgb(26pt)=(0.0383714,0.674271,0.743552); rgb(27pt)=(0.0589714,0.683757,0.725386); rgb(28pt)=(0.0843,0.692833,0.706167); rgb(29pt)=(0.113295,0.7015,0.685857); rgb(30pt)=(0.145271,0.709757,0.664629); rgb(31pt)=(0.180133,0.717657,0.642433); rgb(32pt)=(0.217829,0.725043,0.619262); rgb(33pt)=(0.258643,0.731714,0.595429); rgb(34pt)=(0.302171,0.737605,0.571186); rgb(35pt)=(0.348167,0.742433,0.547267); rgb(36pt)=(0.395257,0.7459,0.524443); rgb(37pt)=(0.44201,0.748081,0.503314); rgb(38pt)=(0.487124,0.749062,0.483976); rgb(39pt)=(0.530029,0.749114,0.466114); rgb(40pt)=(0.570857,0.748519,0.44939); rgb(41pt)=(0.609852,0.747314,0.433686); rgb(42pt)=(0.6473,0.7456,0.4188); rgb(43pt)=(0.683419,0.743476,0.404433); rgb(44pt)=(0.71841,0.741133,0.390476); rgb(45pt)=(0.752486,0.7384,0.376814); rgb(46pt)=(0.785843,0.735567,0.363271); rgb(47pt)=(0.818505,0.732733,0.34979); rgb(48pt)=(0.850657,0.7299,0.336029); rgb(49pt)=(0.882433,0.727433,0.3217); rgb(50pt)=(0.913933,0.725786,0.306276); rgb(51pt)=(0.944957,0.726114,0.288643); rgb(52pt)=(0.973895,0.731395,0.266648); rgb(53pt)=(0.993771,0.745457,0.240348); rgb(54pt)=(0.999043,0.765314,0.216414); rgb(55pt)=(0.995533,0.786057,0.196652); rgb(56pt)=(0.988,0.8066,0.179367); rgb(57pt)=(0.978857,0.827143,0.163314); rgb(58pt)=(0.9697,0.848138,0.147452); rgb(59pt)=(0.962586,0.870514,0.1309); rgb(60pt)=(0.958871,0.8949,0.113243); rgb(61pt)=(0.959824,0.921833,0.0948381); rgb(62pt)=(0.9661,0.951443,0.0755333); rgb(63pt)=(0.9763,0.9831,0.0538)}, mesh/rows=2]
table[row sep=crcr, point meta=\thisrow{c}] {%
x	y	z	c\\
-10	0	0	0\\
-10	10	1	1\\
10	0	0	0\\
10	10	1	1\\
};

\end{axis}
\end{tikzpicture}%
%
	}}%
		\end{minipage}
		\hfill
		\begin{minipage}[t][][t]{0.53\linewidth}\vspace*{0pt}
			{{%
		\pgfkeys{/pgf/images/include external/.code={\includegraphics[width=\linewidth]{#width=\linewidth}}}%
		\tikzsetnextfilename{IceCream}%
		\input{./TeX/Tikz/IceCream.tex}%
	}}%
		\end{minipage}
		\caption[Alternating projections on second-order cone and tangent plane]{%
			\emph{Alternating projections on second-order cone and tangent plane}.
		
			In contrast with the slow sublinear rate of KM both with and without line search, and despite the non isolatedness of any solution, Broyden's scheme exhibits an appealing linear convergence rate.%
		}%
		\label{fig:IceCream}
	\end{figure}%
\fi

However, computing the generalized Jacobian might be too demanding and require extra information not available in close form.
For this reason we focus on \emph{quasi-Newton} methods
\[
	x^+=x-HRx,
\]
where the linear operator \(H\) is progressively updated with only evaluations of \(R\) and direct linear algebra in such a way that the vector \(HRx\) is asymptotically a \emph{good approximation} of a Newton direction \(G^{-1}Rx\).
The yellow lines in the simulations of this section correspond to \(H\) being selected with Broyden's quasi-Newton method.

The crucial issue is \emph{convergence} itself.
Though in these trivial simulations it is not the case, it is well known that Newton-type methods in general converge only when close to a solution, and may even diverge otherwise.
\added{%
	In fact, globalizing the convergence of Newton-type methods is a key challenge in optimization, as the dedicated recent book \cite{izmailov2014newton} confirms.

	In this paper we provide the \refSM[], a globalization strategy for Newton-type methods (or any local scheme in general) that applies to any (nonsmooth) monotone equation deriving from fixed-point iterations of nonexpansive operators.
	Our method covers almost all splitting schemes in convex optimization, such as forward-backward splitting (FBS, also known as proximal gradient method), Douglas-Rachford splitting (DRS) and the alternating direction method of multipliers (ADMM), to name a few.
	We also provide sufficient conditions at the limit point under which the method reduces to the local scheme and converges superlinearly.
}%

	\section{Notation and known results}
		\label{sec:Definitions}
		With \(\boundary A\) we denote the boundary of the set \(A\), and given a sequence \((x_k)_{k\in\N}\) we write \((x_k)_{k\in\N}\subset A\) with the obvious meaning of \(x_k\in A\) for all \(k\in\N\).
For \(p>0\) we let
\[
	\ell^p
{}\coloneqq{}
	\set{\seq{x_k}\subset\R}[
		\textstyle
		\sum_{k\in\N}|x_k|^p<\infty
	]
\]
denote the set of real-valued sequences with summable \(p\)-th power, and with \(\ell_+^p\) the subset of the positive-valued ones.

The positive part of \(x\in\R\) is \([x]_+\coloneqq\max\set{x,0}\).

		\subsection{Hilbert spaces and bounded linear operators}
			Throughout the paper, \(\HH\) is a real separable Hilbert space endowed with an inner product \(\innprod\cdot\cdot\) and with induced norm \(\|{}\cdot{}\|\).
\added{%
	The Euclidean norm and scalar product are denoted as \(\|{}\cdot{}\|_{_2}\) and \(\innprod\cdot\cdot_{\!_2}\), respectively.
}%
For \(\bar x\in\HH\) and \(r>0\), the open ball centered at \(\bar x\) with radius \(r\) is indicated as
\(
	\ball{\bar x}{r}
{}\coloneqq{}
	\set{x\in\HH}[\|x-\bar x\|<r]
\).
For a closed and nonempty convex set \(C\subseteq\HH\) we let \(\proj_C\) denote the projection operator on \(C\).

Given \(\seq{x_k}\subset\HH\) and \(x\in\HH\) we write \(x_k\to x\) and \(x_k\rightharpoonup x\) to denote, respectively, strong and weak convergence of \(\seq{x_k}\) to \(x\).
The set of weak sequential cluster points of \(\seq{x_k}\) is indicated as \(\weak\seq{x_k}\).

The set of bounded linear operators \(\HH\to\HH\) is denoted as \(\linop\).
The adjoint operator of \(L\in\linop\) is indicated as \(\adj L\), \ie, the unique operator in \(\linop\) such that
\(
	\innprod{Lx}{y}
{}={}
	\innprod{x}{\adj Ly}
\)
for all \(x,y\in\HH\).

		\subsection{Nonexpansive operators and Fejér sequences}
			We now briefly recap some known definitions and results of nonexpansive operator theory that will be used in the paper.
\begin{defin}
An operator \(\func T\HH\HH\) is said to be
\begin{enumerate}
	\item
		\DEF{nonexpansive (NE)} if
		\(
			\|Tx-Ty\|\leq\|x-y\|
		\)
		for all \(x,y\in\HH\);
	\item
		\label{defin:avg}
		\DEF{\avg[]} if it is \DEF{\avg} for some \(\alpha\in(0,1)\), \ie if there exists a nonexpansive operator \(\func S\HH\HH\) such that
		\(
			T
		{}={}
			(1-\alpha)\id+\alpha S
		\);
	\item
		\label{defin:FNE}
		\DEF{firmly nonexpansive (FNE)} if it is \(\avg[\nicefrac 12]\).
\end{enumerate}
\end{defin}
Clearly, for any NE operator \(T\) the residual \(R=\id-T\) is monotone, in the sense that \(\innprod{Rx-Ry}{x-y}\geq 0\) for all \(x,y\in\HH\); if \(T\) is additionally FNE, then not only is \(R\) monotone, but it is FNE as well.
For notational convenience we extend the definition of \(\alpha\)-averagedness to the case \(\alpha=1\) which reduces to plain nonexpansiveness.

Given an operator \(\func T\HH\HH\) we let
	\[
		\zer T
	{}\coloneqq{}
		\set{z\in\HH}[Tz=0]
	~~\text{and}~~
		\fix T
	{}\coloneqq{}
		\set{z\in\HH}[Tz=z]
	\]
denote the set of its \DEF{zeros} and \DEF{fixed points}, respectively.
For \(\lambda\in\R\) we define the \(\lambda\)-averaging of \(T\) as
\[
	T_\lambda
{}\coloneqq{}
	(1-\lambda)\id+\lambda T.
\]
Notice that
\begin{equation}\label{eq:AvgRes}
	\id-T_\lambda
{}={}
	\lambda(\id-T)
\qquad
	\text{for all }
	\lambda\in\R,
\end{equation}
and therefore \(\fix T_\lambda=\fix T\) for all \(\lambda\neq 0\).
Moreover, if \(T\) is \avg{} and \(\lambda\in(0,\nicefrac 1\alpha]\), then
\begin{equation}\label{eq:AvgAvg}
	\text{\(T_\lambda\) is \avg[\alpha\lambda]}
\end{equation}
\cite[Cor. 4.28]{bauschke2011convex} and in particular \(T_{\nicefrac{1}{2\alpha}}\) is FNE.


			\begin{defin}\label{defin:Fejer}
Relative to a nonempty set \(S\subseteq\HH\), a sequence \(\seq{x_k}\subset\HH\) is
\begin{enumerate}
	\item
		\DEF{Fejér (-monotone)} if~
		\(
			\|x_{k+1}-s\|
		{}\leq{}
			\|x_k-s\|
		\)~
		for all \(k\in\N\) and \(s\in S\);
	\item\label{defin:QF}
		\DEF{quasi-Fejér (monotone)} if for all \(s\in S\) there exists a sequence
		\(
			\seq{\varepsilon_k(s)}
		{}\in{}
			\ell_+^1
		\)
		such that
		\[
			\|x_{k+1}-s\|^2
		{}\leq{}
			\|x_k-s\|^2
			{}+{}
			\varepsilon_k(s)
		\quad
			\forall k\in\N.
		\]
\end{enumerate}
\end{defin}
This definition of quasi-Fejér monotonicity is taken from \cite{combettes2001quasi} where it is referred to as \emph{of type III}, and generalizes the classical definition \cite{ermolev1973random}.
\begin{thm}\label{prop:WeakFix}
Let \(\func T\HH\HH\) be an NE operator with \(\fix T\neq\emptyset\), and suppose that \(\seq{x_k}\subset\HH\) is quasi-Fejér with respect to \(\fix T\).
If \(\seq{x_k-Tx_k}\to 0\), then there exists \(x_\star\in\fix T\) such that \(x_k\rightharpoonup x_\star\).
\begin{proof}
From \cite[Prop. 3.7(i)]{combettes2001quasi} we have \(\weak\seq{x_k}\neq\emptyset\); in turn, from \cite[Cor. 4.18]{bauschke2011convex} we infer that \(\weak\seq{x_k}\subseteq\fix T\).
The claim then follows from \cite[Thm. 3.8]{combettes2001quasi}.
\end{proof}
\end{thm}

	\section{General abstract framework}
		\label{sec:General}
		Unless differently specified, in the rest of the paper we work under the following assumption.
\begin{ass}\label{ass:T}
\(\func T\HH\HH\) is an \avg{} operator for some \(\alpha\in(0,1]\) and with \(\fix T\neq\emptyset\).
With \(R\coloneqq\id-T\) we denote its (\(2\alpha\)-Lipschitz continuous) fixed-point residual.
\end{ass}
We also stick to this notation, so that, whenever mentioned, \(T\), \(R\), and \(\alpha\) are as in \Cref{ass:T}.
Our goal is to find a fixed point of \(T\), or, equivalently, a zero of \(R\):
\begin{equation}\label{eq:Problem}
	\text{find }
	x_\star\in\fix T=\zer R.
\end{equation}

In this section we introduce \Cref{alg:General}, an abstract procedure to solve problem \eqref{eq:Problem}.
The scheme is not implementable in and of itself, as it gives no hint as to how to compute each of the iterates, but it rather serves as a comprehensive ground framework for a class of algorithms with global convergence guarantees.
In \Cref{sec:SuperMann} we will derive the \emph{SuperMann scheme}, an implementable instance which also enjoys appealing asymptotic properties.

\begin{algorithm*}[tb]
	\algcaption[General framework]{%
		\emph{General framework} for finding a fixed point of the \avg{} operator \(T\) with residual \(R=\id-T\)%
	}%
	\label{alg:General}%
	\newcommand\keyfont[1]{\textsc{#1}}
\begin{tabularx}{\linewidth}{@{}l@{~~}X@{}}
		\keyfont{Require}%
	&
		\(x_0\in\HH\),~
		\(c_0,c_1,q\in[0,1)\),~
		\(\sigma>0\)%
	\\
		\keyfont{Initialize}%
	&
		\(\eta_0=r_{\rm safe}=\|Rx_0\|\),~
		\(k=0\)%
\end{tabularx}
\begin{enumerate}[{label=\textbf{\arabic*.}},{ref=\textbf{\arabic*}}]
	\item\label{step:General:Initial}
		\keyfont{If} \(Rx_k=0\), \keyfont{then stop}.
	\item\label{step:General:K0}
		\keyfont{If} \(\|Rx_k\|\leq c_0\eta_k\),~
		\keyfont{then} set~
		\(
			\eta_{k+1}
		{}={}
			\|Rx_k\|
		\),~
		proceed with a \emph{blind update} \(x_{k+1}\) and go to step \ref{step:General:k++}.
	\item
		Set \(\eta_{k+1}=\eta_k\) and select \(x_{k+1}\) such that
		\begin{enumerate}[{label=\textbf{\theenumi(\alph*)}},{ref=\theenumi(\alph*)}]
			\item\label{step:General:K1}
				\keyfont{either} the \emph{safe condition} \(\|Rx_k\|\leq r_{\rm safe}\) holds, and \(x_{k+1}\) is \emph{educated}:
				\[
					\|Rx_{k+1}\|
				{}\leq{}
					c_1\hspace{0.5pt}\|Rx_k\|
				\]
				in which case update
				\(
					r_{\rm safe}
				{}={}
					\|Rx_{k+1}\|
					{}+{}
					q^k
				\);
			\item\label{step:General:K2}
				\keyfont{or} it is \emph{Fejérian} with respect to \(\fix T\):
				\begin{equation}\label{eq:General:K2}
					\|x_{k+1}-z\|^2
				{}\leq{}
					\|x_k-z\|^2
					{}-{}
					\sigma\|Rx_k\|^2
				\quad
					\forall z\in\fix T.
				\end{equation}
		\end{enumerate}
	\item\label{step:General:k++}
		Set \(k\gets k+1\) and go to step \ref{step:General:Initial}.
\end{enumerate}
\end{algorithm*}

The general framework prescribes three kinds of updates.
\begin{enumerate}[{label=\(K_{\arabic*}\))},{ref=\(K_{\arabic*}\)}]
	\setcounter{enumi}{-1}
	\item\label{K0}
	\textbf{\emph{Blind} updates.}
		Inspired from \cite{chen1999proximal}, whenever the residual \(\|Rx_k\|\) at iteration \(k\) has \emph{sufficiently} decreased with respect to past iterates we allow for an \emph{uncontrolled} update.
		For an efficient implementation such guess should be somehow reasonable and not completely a ``\emph{blind}'' guess; however, for the sake of global convergence the proposed scheme is robust to any choice.

	\item\label{K1}
		\textbf{\emph{Educated} updates.}
		To encourage favorable updates, similarly to what has been proposed in \cite[§5.3.1]{izmailov2014newton} and \cite[§8.3.2]{facchinei2003finite} an \emph{educated guess} \(x_{k+1}\) is accepted whenever the candidate residual is \emph{sufficiently} smaller than the current.
		
	\item\label{K2}
		\textbf{\emph{Safeguard} (Fejérian) updates.}
		This last kind of updates is similar to \ref{K1} as it is also based on the goodness of \(x_{k+1}\) with respect to \(x_k\).
		The difference is that instead of checking the residual, what needs be \emph{sufficiently} decreased is the distance from each point in \(\fix T\).
		This is meant in a Fejérian fashion as in \Cref{defin:Fejer}.
		
\end{enumerate}
\emph{Blind} \ref{K0}- and \emph{educated} \ref{K1}-updates are somehow complementary: the former is enabled when enough progress has been made in the past, whereas the latter when the candidate update yields a sufficient improvement.
\emph{Progress} and \emph{improvement} are meant in terms of a linear decrease of (the norm of) the residual; at iteration \(k\), \ref{K0} is enabled if \(\|Rx_k\|\leq c_0\|Rx_{\bar k}\|\), where \(c_0\in[0,1)\) is a user-defined constant and \(\bar k\) is the last blind iteration before \(k\); \ref{K1} is enabled if \(\|Rx_{k+1}\|\leq c_1\|Rx_k\|\) where \(c_1\in[0,1)\) is another user-defined constant and \(x_{k+1}\) is the candidate next iterate.
To ensure global convergence, \emph{educated} updates are authorized only if the current residual \(\|Rx_k\|\) is not larger than \(\|Rx_{\tilde k+1}\|\) (up to a linearly decreasing error \(q^{\tilde k}\)); here \(\tilde k\) denotes the last \ref{K1}-update before \(k\).

While \emph{blind} \ref{K0}- and \emph{educated} \ref{K1}-updates are in charge of the asymptotic behavior, what makes the algorithm convergent are \emph{safeguard} \ref{K2}-iterations.

		\subsection{Global weak convergence}
			To establish a notation, we partition the set of iteration indices \(K\subseteq\N\) as \(K_0\cup K_1\cup K_2\).
Namely, relative to \Cref{alg:General}, \(K_0\) \(K_1\) and \(K_2\) denote the sets of indices \(k\) passing the test at steps \ref{step:General:K0}, \ref{step:General:K1} and \ref{step:General:K2}, respectively.
Furthermore, we index the sets \(K_0\) and \(K_1\) of \emph{blind} and \emph{educated} updates as
\begin{equation}\label{eq:K0K1}
	K_0
{}={}
	\set{k_1,k_2,\cdots},
\qquad
	K_1
{}={}
	\set{k_1',k_2',\cdots}.
\end{equation}

To rule out trivialities, throughout the paper we work under the assumption that a solution is not found in a finite number of steps, so that the residual of each iterate is always nonzero.
As long as it is well defined, the algorithm therefore produces an infinite number of iterates.

\begin{thm}[{Global convergence of the general framework \Cref{alg:General}}]\label{thm:General:Global}
Consider the iterates generated by \Cref{alg:General} and suppose that for all \(k\) it is always possible to find a point \(x_{k+1}\) complying with the requirements of either step \ref{step:General:K0}, \ref{step:General:K1} or \ref{step:General:K2}, and further satisfying
\begin{equation}\label{eq:D}
	\|x_{k+1}-x_k\|
{}\leq{}
	D\|Rx_k\|
\quad
	\forall k\in K_0\cup K_1
\end{equation}
for some constant \(D\geq 0\).
Then,
\begin{enumerate}
	\item\label{thm:General:QF}
		\(\seq{x_k}\) is quasi-Fejér monotone with respect to \(\fix T\);
	\item\label{thm:General:ResStr}
		\(Rx_k\to 0\) with \(\seq{\|Rx_k\|}\in\ell^2\);
	\item\label{thm:General:xWeak}
		\(\seq{x_k}\) converges weakly to a point \(x_\star\in\fix T\);
	\item\label{thm:General:K0}
		if \(c_0>0\) the number of \emph{blind updates} at step \ref{step:General:K0} is infinite.
\end{enumerate}
\begin{proof}
See \Cref{proof:thm:General:Global}.
\end{proof}
\end{thm}

		\subsection{Local linear convergence}
			\label{sec:GeneralLinear}
			More can be said about the convergence rates if the mapping \(R\) possesses \DEF{metric subregularity}.
Differently from (bounded) linear regularity \cite{bauschke2015linear}, metric subregularity is a local property and as such it is more general.
For a (possibly multivalued) operator \(R\), metric subregularity at \(\bar x\) is equivalent to calmness of \(R^{-1}\) at \(R\bar x\) \cite[Thm 3.2]{dontchev2004regularity}, and is a weaker condition than metric regularity and Aubin property.
We refer the reader to \cite[§9]{rockafellar2011variational} for an extensive discussion.
\begin{defin}[Metric subregularity at zeros]\label{defin:Subreg}
Let \(\func R\HH\HH\) and \(\bar x\in\zer R\).
\(R\) is \DEF{metrically subregular} at \(\bar x\) if there exist \(\varepsilon,\gamma>0\) such that
\begin{equation}\label{eq:metricsubregularity}
	\dist(x,\zer R)
{}\leq{}
	\gamma
	\|Rx\|
\qquad
	\forall x\in\ball{\bar x}{\varepsilon}.
\end{equation}
\(\gamma\) and \(\varepsilon\) are (one) \DEF{modulus} and (one) \DEF{radius} of subregularity of \(R\) at \(\bar x\), respectively.
\end{defin}
In finite-dimensional spaces, if \(R\) is differentiable at \(\bar x\in\zer R\) and \(\bar x\) is isolated in \(\zer R\) (\eg if it is the unique zero), then metric subregularity is equivalent to nonsingularity of \(JR\bar x\).
Metric subregularity is however a much weaker property than nonsingularity of the Jacobian, firstly because it does not assume differentiability, and secondly because it can cope with `wide' regions of zeros; for instance, any piecewise linear mapping is globally metrically subregular \cite{robinson1981some}.

If the residual \(R=\id-T\) of the \avg{} operator \(T\) is metrically subregular at \(\bar x\in\zer R=\fix T\) with modulus \(\gamma\) and radius \(\varepsilon\), then
\begin{equation}\label{eq:distR}
	\tfrac 1\gamma
	\dist(x,\fix T)
{}\leq{}
	\|Rx\|
{}\leq{}
	2\alpha
	\dist(x,\fix T)
\ifieee\else
	\qquad
	\forall x\in\ball{\bar x}{\varepsilon}.
\fi
\end{equation}
\ifieee
	for all \(x\in\ball{\bar x}{\varepsilon}\).
\fi
Consequently, if \(\|Rx_k\|\to 0\) for some sequence \(\seq{x_k}\subset\HH\), so does \(\dist(x_k,\fix T)\) with the same asymptotic rate of convergence, and viceversa.
Metric subregularity is the key property under which the residual in the classical KM scheme achieves linear convergence; in the next result we show that this asymptotic behavior is preserved in the general framework of \Cref{alg:General}.

\begin{thm}[{Linear convergence of the general framework \Cref{alg:General}}]\label{thm:General:Linear}
Suppose that the hypotheses of \Cref{thm:General:Global} hold, and suppose further that \(\seq{x_k}\) converges strongly to a point \(x_\star\) (this being true if \(\HH\) is finite dimensional) at which \(R\) is metrically subregular.

Then, \(\seq{x_k}\) and \(\seq{Rx_k}\) are \(R\)-linearly convergent.
\begin{proof}
See \Cref{proof:thm:General:Linear}.
\end{proof}
\end{thm}

		\subsection{Main idea}
			Being interested in solving the nonlinear equation \eqref{eq:Problem}, one could think of implementing one of the many existing \emph{fast} methods for nonlinear equations that achieve fast asymptotic rates, such as Newton-type schemes.
At each iteration, such schemes compute an update direction \(d_k\) and prescribe steps of the form \(x_{k+1}=x_k+\tau_kd_k\), where \(\tau_k>0\) is a stepsize that needs be sufficiently small in order for the method to enjoy global convergence; on the other hand, fast asymptotic rates are ensured if \(\tau_k=1\) is eventually always accepted.
The stepsize is a crucial feature of fast methods, and a feasible \(\tau_k\) is usually backtracked with a line search on a smooth merit function.
Unfortunately, in meaningful applications of the problem at hand arising from fixed-point theory the residual mapping \(R\) is nonsmooth, and the typical merit function \(x\mapsto\|Rx\|^2\) does not meet the necessary smoothness requirement.

What we propose in this paper is a hybrid scheme that allows for the employment of any (fast) method for solving nonlinear equations, with global convergence guarantees that do not require smoothness, but which is based only on the nonexpansiveness of \(T\).
Once \emph{fast} directions \(d_k\) are selected,
\Cref{alg:General} can be specialized as follows:
\begin{enumerate}
	\item
		\emph{blind} updates as in step \ref{step:General:K0} shall be of the form \(x_{k+1}=x_k+d_k\);
	\item
		\emph{educated} updates as in step \ref{step:General:K1} shall be of the form \(x_{k+1}=x_k+\tau_kd_k\), with \(\tau_k\) small enough so as to ensure the acceptance condition \(\|Rx_{k+1}\|\leq c_1\|Rx_k\|\);
	\item
		\emph{safeguard} updates as in step \ref{step:General:K2} shall be employed as \emph{last resort} both for globalization purposes and for well definedness of the scheme.
\end{enumerate}
Ideally, the scheme should eventually reduce to the local scheme \(x_{k+1}=x_k+d_k\) when \emph{good} directions \(d_k\) are used.

In \Cref{sec:GenMann} we address the problem of providing explicit \emph{safeguard} updates that comply with the quasi-Fejér monotonicity requirement of step \ref{step:General:K2}.
Because of the arbitrarity of the other two updates, once we succeed in this task \Cref{alg:General} will be of practical implementation.
In \Cref{sec:SuperMann} we will then discuss specific \ref{K0}- and \ref{K1}-updates to be used at steps \ref{step:General:K0} and \ref{step:General:K1} that ensure global \emph{and} fast convergence, yet maintaining the simplicity of fixed-point iterations of \(T\) (evaluations of \(T\) and direct linear algebra).

	\section{Generalized Mann Iterations}
		\label{sec:GenMann}

		\subsection{The classical \KM{} scheme}
			Starting from a point \(x_0\in\HH\), the classical \KM{} scheme (KM) performs the following updates
\begin{equation}\label{eq:KM}
	x_{k+1}
{}={}
	T_{\lambda_k}x_k
{}={}
	(1-\lambda_k)x_k + \lambda_kTx_k
\end{equation}
and converges weakly to a fixed point of \(T\) provided that \(\lambda_k\in[0,\nicefrac1\alpha]\) and
\(
	\seq{\lambda_k(\nicefrac1\alpha-\lambda_k)}
	\notin
	\ell^1
\)
\cite[Thm. 5.14]{bauschke2011convex}.
The key property of KM iterations is Fejér monotonicity:
\[
	\|x_{k+1}-z\|^2
{}\leq{}
	\|x_k-z\|^2
	{}-{}
	\lambda_k(\nicefrac1\alpha-\lambda_k)
	\|Rx_k\|^2
\quad
	\forall z\in\fix T.
\]
In particular, in \Cref{alg:General} KM iterations can be used as \emph{safeguard} updates at step \ref{step:General:K2}.
The drawback of such a selection is that it completely discards the hypothetical fast update direction \(d_k\) that \emph{blind} and \emph{educated} updates try to enforce.
This is particularly penalizing when the local method for computing the directions \(d_k\) is a \emph{quasi-Newton} scheme; such methods are indeed very sensitive to past iterations, and discarding directions is neither theoretically sound nor beneficial in practice.

\begin{figure}[tb]
	\setlength\fboxsep{0pt}%
	\begin{minipage}[t][][t]{\linewidth}%
		\vspace*{0pt}%
		\centering
		{{%
		\pgfkeys{/pgf/images/include external/.code={\includegraphics[width=0.5\linewidth]{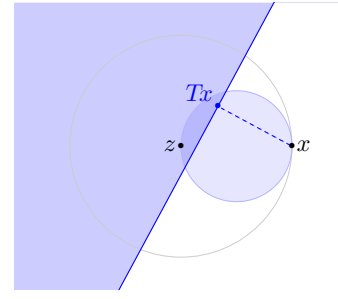}}}%
		\tikzsetnextfilename{Mann}%
		\def\radius{2}
\def\xmax{3}
\begin{tikzpicture}%
	\pgfplotsset{%
		mystyle/.style={%
			ticks=none,
			axis lines=none,
			xmin={-\xmax},
			xmax={\xmax},
			ymin={-\xmax},
			ymax={\xmax},
			enlargelimits=false,
			scale only axis,
		},%
	}%
	\tikzset{%
		textnode/.style={%
			scale=1.75,
			inner sep=0pt,
			outer sep=2pt,
		},%
		nosep/.style={%
			inner sep=0pt,
			outer sep=0pt,
		}%
	}%
	\begin{axis}[mystyle]%
		\draw[draw=none,name path=top] (axis cs:-\xmax,\xmax) -- (axis cs:\xmax,\xmax) {};
		\draw[draw=none,name path=bottom] (axis cs:-\xmax,-\xmax) -- (axis cs:\xmax,-\xmax) {};
		\draw[draw=none,name path=left] (axis cs:-\xmax,-\xmax) -- (axis cs:-\xmax,\xmax) {};
		\draw[draw=none,name path=right] (axis cs:\xmax,-\xmax) -- (axis cs:\xmax,\xmax) {};
		\node[nosep] (z) at (axis cs:0,0) {};
		\node[nosep] (x) at (axis cs:\radius,0) {};
		\node[nosep] (Tx) at (axis cs:\radius/3,{1.25*\radius/3}) {};
		\node (circle_z) at (z) [
			draw,
			very thin,
			black!20,
			circle through=(x),
		] {};
		\node (xOz) at ($(x)!0.5!(z)$) {};
		\node (circle_Tx) at (xOz) [
			draw,
			thin,
			blue!30,
			fill=blue,
			fill opacity=0.1,
			circle through=(x),
		] {};
		\draw[blue,thick,name path=C_x] ($(Tx)!20!90:(x)$) -- ($(Tx)!20!-90:(x)$);
		\draw[blue!50,very thin,dashed] (x) -- (Tx);
		\addplot[fill=blue,fill opacity=0.2] fill between[of = C_x and top];
		\node at (z) {\textbullet};
		\node at (x) {\textbullet};
		\node[blue] at (Tx) {\textbullet};
		
		\node[textnode] at (z) [anchor=east] {\(z\)};
		\node[textnode] at (x) [anchor=west] {\(x\)};
		\node[textnode, blue] at (Tx) [anchor=south east] {\(T\!x\)};
		\draw[
			blue,
			very thin,
			dashed,
			shorten >=2pt,
			shorten <=1pt,
		] (x) -- (Tx);
	\end{axis}%
\end{tikzpicture}%
%
	}}%
	\end{minipage}
	\begin{minipage}[t][][t]{\linewidth}%
		\vspace*{0pt}%
		\caption[Mann iteration as projection on \(C_x\)]{\emph{Mann iteration of a FNE operator \(T\) as projection on \(C_x\)} (the blue half-space, as defined in \eqref{eq:C_x} for \(\alpha=\nicefrac12\)).
			The outer circle is the set of all possible images of a nonexpansive operator, given that \(z\) is a fixed point.
			The inner circle corresponds to the possible images of \emph{firmly} nonexpansive operators.
			Notice that \(C_x\) separates \(x\) from \(z\) as long as \(Tx\) is contained in the small circle, which characterizes firm nonexpansiveness.%
		}%
		\label{fig:Mann}%
	\end{minipage}%
\end{figure}

In this section we provide alternative \emph{safeguard} updates that while ensuring the desirable Fejér monotonicity are also amenable to taking into account arbitrary directions.
\added{%
	The key idea lies in intepreting KM iterations as projections onto suitable half-spaces (see \cref{fig:Mann}), and then exploiting known properties of projections.
	These facts are shown in the next result.
	To this end, let us remark that the projection \(\proj_C\) onto a nonempty closed and convex set \(C\) is FNE \cite[Prop. 4.8]{bauschke2011convex}, and that consequently its \(\lambda\)-averaging \(\proj_{C,\lambda}\) is \avg[\nicefrac\lambda2] for any \(\lambda\in(0,2]\), as it follows from \eqref{eq:AvgAvg}.
}%
\begin{prop}[KM iterations as projections]\label{lem:KMProjection}
For \(x\in\HH\), define
\begin{equation}\label{eq:C_x}
\ifieee
	\mathtight
\fi
	C_x
{}={}
	C_x^{T,\alpha}
{\coloneqq{}}
	\set{z\in\HH}[
		\|Rx\|^2 - 2\alpha\innprod{Rx}{x-z}
	{}\leq{}
		0
	].
\end{equation}
Then,
\begin{enumerate}
	\item\label{item:xC_x}
		\(x\in C_x\) iff \(x\in\fix T\);
	\item\label{item:fixC_x}
		\(
			\fix T
		{}={}
			\bigcap_{x\in\HH}C_x
		\);
	\item\label{item:KMProjection}
		for any \(\lambda\in[0,\nicefrac 1\alpha]\) it holds that
		\(\displaystyle
			T_\lambda x
		{}={}
			\proj_{C_x,2\alpha\lambda}x
		{}={}
			(1-2\alpha\lambda)x
			{}+{}
			2\alpha\lambda\proj_{C_x}x
		\).
\end{enumerate}
\begin{proof}
The set \(C_x\) can be equivalently expressed as
\[
	C_x
{}={}
	\set{z\in\HH}[
		\innprod{x-T_{\nicefrac{1}{2\alpha}}x}{z-T_{\nicefrac{1}{2\alpha}}x}
	{}\leq{}
		0
	].
\]
\ref{item:xC_x} is of immediate verification, and \ref{item:fixC_x} then follows from \cite[Cor. 4.16]{bauschke2011convex} combined with \eqref{eq:AvgAvg}.

We now show \ref{item:KMProjection}.
If \(Rx=0\), then \(x\in\fix T\) and \(C_x=\HH\), and the claim is trivial.
Otherwise, notice that
\begin{equation}\label{eq:Cx2}
	C_x
{}={}
	\set{z\in\HH}[
		\innprod{Rx}{z}
	{}\leq{}
		\innprod{Rx}{x-\tfrac{1}{2\alpha}Rx}
	]
\end{equation}
and the claim can be readily verified using the formula for the projection on a halfspace
\(
	H_{v,\beta}
{}\coloneqq{}
	\set{
		z\in\HH
	}[
		\innprod vz\leq\beta
	]
\)
\begin{equation}\label{eq:HalfSpaceProjection}
	\proj_{H_{v,\beta}}x
{}={}
	x
	{}-{}
	\frac{
		[\innprod vx-\beta]_+
	}{
		\|v\|^2
	}
	v
\end{equation}
defined for \(v\in\HH\setminus\set{0}\) and \(\beta\in\R\) \cite[Ex. 28.16(iii)]{bauschke2011convex}.
\end{proof}
\end{prop}

		\subsection{Generalized Mann projections}
			Though particularly attractive for its simplicity and global convergence properties, the KM scheme \eqref{eq:KM} finds its main drawback in its convergence rate, being \(Q\)-linear at best and highly sensitive to ill conditioning of the problem.
In response to these issues, \Cref{alg:General} allows for the integration of fast local methods still ensuring global convergence properties.
The efficiency of the resulting scheme, which will be proven later on, is based on an ad hoc selection of \emph{safeguard} updates for step \ref{step:General:K2} which is based on the following generalization of \Cref{lem:KMProjection}.
\begin{prop}\label{lem:DriftedKM}
Suppose that \(x,w\in\HH\) are such that
\begin{equation}\label{eq:SepH}
	\rho
{}\coloneqq{}
	\|Rw\|^2 - 2\alpha\innprod{Rw}{w-x}
{}>{}
	0.
\end{equation}
For \(\lambda\in[0,\nicefrac1\alpha]\) let
\begin{equation}\label{eq:DriftedKM}
	x^+
{}\coloneqq{}
	x-\lambda\frac{\rho}{\|Rw\|^2}Rw.
\end{equation}
Then, the following hold:
\begin{enumerate}
	\item\label{prop:GenMannproj}
		\(x^+=\proj_{C_w,2\alpha\lambda}x\) where \(C_w=C_w^{T,\alpha}\) as in \eqref{eq:C_x};
	\item\label{prop:DriftedKMFejer}
		\(
			\|x^+-z\|^2
		{}\leq{}
			\|x-z\|^2
			{}-{}
			\lambda(\nicefrac1\alpha-\lambda)
			\frac{\rho^2}{\|Rw\|^2}
		~~
			\forall z\in\fix T\).
\end{enumerate}
\begin{proof}
\ref{prop:GenMannproj} easily follows from \eqref{eq:Cx2} and \eqref{eq:HalfSpaceProjection}, since by condition \eqref{eq:SepH} the positive part in the formula may be omitted.
In turn, \ref{prop:DriftedKMFejer} follows from \cite[Prop. 4.25(iii)]{bauschke2011convex} by observing that \(\proj_{C_w,2\alpha\lambda}\) is \avg[\alpha\lambda]
due to \cite[Prop.s 4.8]{bauschke2011convex} and \eqref{eq:AvgAvg}, and that \(\fix T\subseteq C_w\) as shown in \cref{item:fixC_x}.
\end{proof}
\end{prop}
Notice that condition \eqref{eq:SepH} is equivalent to \(x\notin C_w\).
Therefore, \Cref{prop:DriftedKMFejer} states that whenever a point \(x\) lies outside the half-space \(C_w\) for some \(w\in\HH\), since \(\fix T\subseteq C_w\) (cf. \cref{lem:KMProjection}) the projection onto \(C_w\) moves closer to \(\fix T\).
This means that after moving from \(x\) along a candidate direction \(d\) to the point \(w=x+d\), even though \(w\) might be farther from \(\fix T\) the point \(x^+=\proj_{C_w}x\) is not.
We may then use this projection as a safeguard step to prevent from diverging from the set of fixed points.
Based on this, we define a \DEF{generalized KM update along a direction \(d\)}.

\begin{defin}[GKM update]\label{defin:GKM}
A \DEF{generalized KM update (GKM) at \(x\) along \(d\)} for the \avg{} operator \(\func T\HH\HH\) with relaxation \(\lambda\in[0,\nicefrac1\alpha]\) is
\[
	x^+
{}\coloneqq{}
	\begin{cases}[l@{~~}l]
		x & \text{if } w\in\fix T \\
		x-\lambda\frac{[\rho]_+}{\|Rw\|^2}Rw & \text{othwerwise,}
	\end{cases}
\]
where \(w=x+d\) and
\(
	\rho
{}\coloneqq{}
	\|Rw\|^2 - 2\alpha\innprod{Rw}{w-x}
\).
In particular, \(d=0\) yields the classical KM update \(x^+=T_\lambda x\).
\end{defin}

		\subsection{Line search for GKM}
\ifieee
	\begin{figure*}
\else
	\begin{figure}[tb]
\fi
	\setlength\fboxsep{0pt}%
	\begin{subfigure}{0.33\textwidth}%
		\vspace*{0pt}%
		{{%
		\pgfkeys{/pgf/images/include external/.code={\includegraphics[width=\linewidth]{#width=\linewidth}}}%
		\tikzsetnextfilename{SuperMann1}%
		\def\radius{2}
\def\xmax{3}
\begin{tikzpicture}%
	\pgfplotsset{%
		mystyle/.style={%
			ticks=none,
			axis lines=none,
			xmin={-\xmax},
			xmax={\xmax},
			ymin={-\xmax},
			ymax={\xmax},
			enlargelimits=false,
			scale only axis,
		},%
	}%
	\tikzset{%
		textnode/.style={%
			scale=1.75,
			inner sep=0pt,
			outer sep=2pt,
		},%
		nosep/.style={%
			inner sep=0pt,
			outer sep=0pt,
		}%
	}%
	\begin{axis}[mystyle]
		\draw[draw=none,name path=top] (axis cs:-\xmax,\xmax) -- (axis cs:\xmax,\xmax) {};
		\draw[draw=none,name path=bottom] (axis cs:-\xmax,-\xmax) -- (axis cs:\xmax,-\xmax) {};
		\draw[draw=none,name path=left] (axis cs:-\xmax,-\xmax) -- (axis cs:-\xmax,\xmax) {};
		\draw[draw=none,name path=right] (axis cs:\xmax,-\xmax) -- (axis cs:\xmax,\xmax) {};
		\node[nosep] (z) at (axis cs:0,0) {};
		\node[nosep] (x) at (axis cs:\radius,0) {};
		\node[nosep] (Tx) at (axis cs:\radius/3,{1.25*\radius/3}) {};
		\node[nosep] (w) at (axis cs:{0.7*\radius},{-0.9*\radius}) {};
		\node (circle_z) at (z) [
			draw,
			very thin,
			black!20,
			circle through=(x),
		] {};
		\node (xOz) at ($(x)!0.5!(z)$) {};
		\node (wOz) at ($(w)!0.5!(z)$) {};
		\node at (wOz) [
			draw,
			name path=circle_wOz,
			circle through=(z),
			orange!50,
			fill=orange!20,
			fill opacity=0.3,
			ultra thin,
		] {};
		\node (xOw) at ($(x)!0.5!(w)$) {};
		\node (circle_xOw) at (xOw) [
			draw,
			dashed,
			very thin,
			black!20,
			circle through=(x),
		] {};
		\node (OTw) at ($(xOw)+(Tx)-(x)$) {};
		\draw[
			dashed,
			orange!50,
			fill=orange!20,
			fill opacity=0.3,
			ultra thin,
		] \getDistance{x}{w} (OTw) circle (\n1);
		\node at (z) {\textbullet};
		\node at (x) {\textbullet};
		\node[blue] at (Tx) {\textbullet};
		\node[orange] at (w) {\textbullet};
		
		\node[textnode] at (z) [anchor=east] {\(z\)};
		\node[textnode] at (x) [anchor=west] {\(x\)};
		\node[textnode, blue] at (Tx) [anchor=south east] {\(T\!x\)};
		\node[textnode] at (w) [anchor=north west] {\(w\)};
		\node[textnode,black!50,scale=0.8] at (axis cs:\radius,-\radius/2) [anchor=west] {\(B_{x,w}\)};
		\draw[
			-{Stealth[scale=1.5]},
			ultra thin,
			black!50!white,
			shorten >=2pt,
			shorten <=1pt,
		] (x) -- node[textnode, pos=0.5, black, anchor=east] {\(d\)} (w);
		\clip \getDistance{w}{z} (wOz) circle (\n1);
		\fill[
			orange,
			fill opacity=0.1,
		] \getDistance{w}{x} (OTw) circle (\n1);
	\end{axis}%
\end{tikzpicture}%
%
	}}%
		\caption{}%
		\label{fig:SuperMann1}%
	\end{subfigure}%
	\hfill
	\begin{subfigure}{0.33\textwidth}%
		\vspace*{0pt}%
		{{%
		\pgfkeys{/pgf/images/include external/.code={\includegraphics[width=\linewidth]{#width=\linewidth}}}%
		\tikzsetnextfilename{SuperMann2}%
		\def\radius{2}
\def\xmax{3}
\begin{tikzpicture}%
	\pgfplotsset{%
		mystyle/.style={%
			ticks=none,
			axis lines=none,
			xmin={-\xmax},
			xmax={\xmax},
			ymin={-\xmax},
			ymax={\xmax},
			enlargelimits=false,
			scale only axis,
		},%
	}%
	\tikzset{%
		textnode/.style={%
			scale=1.75,
			inner sep=0pt,
			outer sep=2pt,
		},%
		nosep/.style={%
			inner sep=0pt,
			outer sep=0pt,
		},%
	}%
	\begin{axis}[mystyle]
		\draw[draw=none,name path=top] (axis cs:-\xmax,\xmax) -- (axis cs:\xmax,\xmax) {};
		\draw[draw=none,name path=bottom] (axis cs:-\xmax,-\xmax) -- (axis cs:\xmax,-\xmax) {};
		\draw[draw=none,name path=left] (axis cs:-\xmax,-\xmax) -- (axis cs:-\xmax,\xmax) {};
		\draw[draw=none,name path=right] (axis cs:\xmax,-\xmax) -- (axis cs:\xmax,\xmax) {};
		\node[nosep] (z) at (axis cs:0,0) {};
		\node[nosep] (x) at (axis cs:\radius,0) {};
		\node[nosep] (Tx) at (axis cs:\radius/3,{1.25*\radius/3}) {};
		\node[nosep] (w) at (axis cs:{0.7*\radius},{-0.9*\radius}) {};
		\node[nosep] (Tw) at ($(x)!0.5!(w)+(Tx)-(x) + (-30:{transformdirectionx(0.8)})$) {};
		\node (circle_z) at (z) [
			draw,
			very thin,
			black!20,
			circle through=(x),
		] {};
		\node (wOz) at ($(w)!0.5!(z)$) {};
		\node at (wOz) [
			draw=none,
			name path=circle_wOz,
			circle through=(z),
		] {};
		\node (xOw) at ($(x)!0.5!(w)$) {};
		\node (OTw) at ($(xOw)+(Tx)-(x)$) {};
		\node (circle_xOw) at (xOw) [
			draw,
			dashed,
			very thin,
			black!20,
			circle through=(x),
		] {};
		\draw[
			draw=none,
			dashed,
		] \getDistance{x}{w} (OTw) circle (\n1);
		\draw[
			orange,
			thick,
			name path=C_w,
		] ($(Tw)!5!90:(w)$) -- ($(Tw)!5!-90:(w)$);
		\draw[
			orange!50,
			very thin,
			dashed,
		] (w) -- (Tw);
		\addplot[
			fill=orange,
			fill opacity=0.1,
		] fill between[
			of = C_w and top,
		];
		\node at (z) {\textbullet};
		\node at (x) {\textbullet};
		\node[blue] at (Tx) {\textbullet};
		\node[orange] at (w) {\textbullet};
		\node[black] at (Tw) {\textbullet};
		
		\node[textnode] at (z) [anchor=east] {\(z\)};
		\node[textnode] at (x) [anchor=west] {\(x\)};
		\node[textnode, blue] at (Tx) [anchor=south east] {\(T\!x\)};
		\node[textnode] at (w) [anchor=north west] {\(w\)};
		\node[textnode] at (Tw) [anchor=south east] {\(T\!w\)};
		\node[textnode,black!50,scale=0.8] at (axis cs:\radius,-\radius/2) [anchor=west] {\(B_{x,w}\)};
		\draw[
			-{Stealth[scale=1.5]},
			ultra thin,
			black!50!white,
			shorten >=2pt,
			shorten <=1pt,
		] (x) -- node[textnode, pos=0.5, black, anchor=east] {\(d\)} (w);
		\clip \getDistance{w}{z} (wOz) circle (\n1);
		\fill[
			orange,
			fill opacity=0.1,
		] \getDistance{w}{x} (OTw) circle (\n1);
		\clip \getDistance{x}{w} (OTw) circle (\n1);
		\node (circle_xOw) at (xOw) [
			draw,
			dashed,
			thick,
			red,
			fill=red,
			fill opacity=0.3,
			pattern color=red,
			circle through=(x),
		] {};
	\end{axis}%
\end{tikzpicture}%
%
	}}%
		\caption{}%
		\label{fig:SuperMann2}%
	\end{subfigure}%
	\hfill
	\begin{subfigure}{0.33\textwidth}%
		\vspace*{0pt}%
		{{%
		\pgfkeys{/pgf/images/include external/.code={\includegraphics[width=\linewidth]{#width=\linewidth}}}%
		\tikzsetnextfilename{SuperMann3}%
		\def\radius{2}
\def\xmax{3}
\begin{tikzpicture}%
	\pgfplotsset{%
		mystyle/.style={%
			ticks=none,
			axis lines=none,
			xmin={-\xmax},
			xmax={\xmax},
			ymin={-\xmax},
			ymax={\xmax},
			enlargelimits=false,
			scale only axis,
		},%
	}%
	\tikzset{%
		textnode/.style={%
			scale=1.75,
			inner sep=0pt,
			outer sep=2pt,
		},%
		nosep/.style={%
			inner sep=0pt,
			outer sep=0pt,
		},%
	}%
	\begin{axis}[mystyle]
		\draw[draw=none,name path=top] (axis cs:-\xmax,\xmax) -- (axis cs:\xmax,\xmax) {};
		\draw[draw=none,name path=bottom] (axis cs:-\xmax,-\xmax) -- (axis cs:\xmax,-\xmax) {};
		\draw[draw=none,name path=left] (axis cs:-\xmax,-\xmax) -- (axis cs:-\xmax,\xmax) {};
		\draw[draw=none,name path=right] (axis cs:\xmax,-\xmax) -- (axis cs:\xmax,\xmax) {};
		\node[nosep] (z) at (axis cs:0,0) {};
		\node[nosep] (x) at (axis cs:\radius,0) {};
		\node[nosep] (Tx) at (axis cs:\radius/3,{1.25*\radius/3}) {};
		\node[nosep] (w_old) at (axis cs:{0.7*\radius},{-0.9*\radius}) {};
		\node[nosep] (w) at (axis cs:{0.85*\radius},{-0.45*\radius}) {};
		\node[nosep] (Tw) at ($(x)!0.5!(w) + (Tx) - (x) + (-45:{transformdirectionx(0.4)})$) {};
		\node (circle_z) at (z) [
			draw,
			very thin,
			black!20,
			circle through=(x),
		] {};
		\node (wOz) at ($(w)!0.5!(z)$) {};
		\node at (wOz) [
			draw=none,
			name path=circle_wOz,
			circle through=(z),
		] {};
		\node (xOw) at ($(x)!0.5!(w)$) {};
		\node (OTw) at ($(xOw)+(Tx)-(x)$) {};
		\node (circle_xOw) at (xOw) [
			draw,
			dashed,
			very thin,
			black!20,
			circle through=(x),
		] {};
		\draw[
			draw=none,
			dashed,
			ultra thin,
		] \getDistance{x}{w} (OTw) circle (\n1);
		\node (C_w1) at ($(Tw)!5!90:(w)$) {};
		\node (C_w2) at ($(Tw)!5!-90:(w)$) {};
		\draw[
			orange,
			thick,
			name path=C_w,
		] (C_w1) -- (C_w2);
		\draw[
			orange!50,
			very thin,
			dashed,
		] (w) -- (Tw);
		\addplot[
			fill=orange,
			fill opacity=0.075,
		] fill between[
			of = C_w and top,
		];
		\node[nosep] (x_plus) at ($(C_w1)!(x)!(C_w2)$) {};
		\node at (z) {\textbullet};
		\node at (x) {\textbullet};
		\node[blue] at (Tx) {\textbullet};
		\node[orange!50!white] at (w_old) {\textbullet};
		\node[orange] at (w) {\textbullet};
		\node[black] at (Tw) {\textbullet};
		\node[blue] at (x_plus) {\textbullet};
		
		\node[textnode] at (z) [anchor=east] {\(z\)};
		\node[textnode] at (x) [anchor=west] {\(x\)};
		\node[textnode,blue] at (Tx) [anchor=south] {\(T\!x\)};
		\node[textnode,orange,opacity=0.4] at (w_old) [anchor=north west] {\(w\)};
		\node[textnode,orange] at (w) [anchor=north west] {\(w\)};
		\node[textnode] at (Tw) [anchor=west] {\(T\!w\)};
		\node[textnode,blue] at (x_plus) [anchor=south] {\(x^+\)};
		\draw[
			-{Stealth[scale=1.5]},
			ultra thin,
			black!50!white,
			shorten >=2pt,
			shorten <=1pt,
		] (x) -- node[textnode, pos=0.5, black, anchor=west] {\(\tau d\)} (w);
		\draw[
			-{Stealth[scale=1.5]},
			ultra thin,
			black!50!white,
			opacity=0.4,
			shorten >=2pt,
			shorten <=1pt,
		] (x) -- (w_old);
		\draw[
			blue,
			very thin,
			dashed,
		] (x) -- (x_plus);
		\clip \getDistance{w}{z} (wOz) circle (\n1);
		\fill[
			orange,
			fill opacity=0.1,
		] \getDistance{w}{x} (OTw) circle (\n1);
	\end{axis}%
\end{tikzpicture}%
%
	}}%
		\caption{}%
		\label{fig:SuperMann3}%
	\end{subfigure}%
	\caption[SuperMann iteration]{%
		\emph{SuperMann iteration of a FNE operator \(T\) as projection on \(C_w\)}.
		
		\textbf{(a)} the darker orange region represents the area in which \(Tw\) must lie given the points \(x\), \(Tx\) and the fixed point \(z\) as prescribed by firm nonexpansiveness of \(T\).
		
		\textbf{(b)} if \(Tw\) lies (also) in the ball \(B_{x,w}\) as in \eqref{eq:Bxw}, then the half-space \(C_w\) (shaded in orange) separates \(x\) from \(w\), which is to be avoided.
		
		\textbf{(c)} when \(w\) is close enough to \(x\) the feasible region for \(Tw\) has empty intersection with \(B_{x,w}\) and \(C_w\) does not contain \(x\).%
	}%
	\label{fig:SuperMann}%
\ifieee
	\end{figure*}
\else
	\end{figure}
\fi
It is evident from \Cref{defin:GKM} that a GKM update trivializes to \(x^+=x\) if either \(w\in\fix T\) or \(\rho\leq 0\).
Having \(w\in\fix T\) corresponds to having found a solution to problem \eqref{eq:Problem}, and the case deserves no further investigation.
In this section we address the remaining case \(\rho\leq 0\), showing how it can be avoided by simply introducing a suitable line search.
In order to recover the same global convergence properties of the classical KM scheme we need something more than simply imposing \(\rho>0\).
The next result addresses this requirement, showing further that it is achieved for any direction \(d\) by sufficiently small stepsizes.
\begin{thm}\label{thm:LS}
	Let \(x,d\in\HH\) and \(\sigma\in[0,1)\) be fixed, and consider
	\[
		\bar\tau
	{}={}
		\begin{cases}[c@{~~}l]
			1 & \text{if } d=0\\
			\tfrac{1-\sigma}{4\alpha}\tfrac{\|Rx\|}{\|d\|} & \text{otherwise}.
		\end{cases}
	\]
	Then, for all \(\tau\in(0,\bar\tau]\) the point \(w=x+\tau d\) satisfies
	\begin{equation}\label{eq:LS}
		\rho
	{}\coloneqq{}
		\|Rw\|^2 - 2\alpha\innprod{Rw}{w-x}
	{}\geq{}
		\sigma
		\|Rw\|\|Rx\|.
	\end{equation}
	\begin{proof}
	Let a constant \(c\geq 0\) to be determined be such that
	\[
		\tau\|d\|
	{}={}
		\|w-x\|
	{}\leq{}
		c\|Rx\|.
	\]
	Observe that
	\(
		\rho
	{}={}
		4\alpha^2
		\innprod{
			w-T_{\nicefrac{1}{2\alpha}}w
		}{
			x-T_{\nicefrac{1}{2\alpha}}w
		}
	\),
	and recall from \eqref{eq:AvgRes} and \eqref{eq:AvgAvg} that \(T_{\nicefrac{1}{2\alpha}}\) is FNE with residual
	\(
		\id-T_{\nicefrac{1}{2\alpha}}
	{}={}
		\frac{1}{2\alpha}R
	\).
	Then,
	\begin{align*}
		\rho
	{}={} &
		4\alpha^2
		\left(
			\|w-T_{\nicefrac{1}{2\alpha}}w\|^2
			{}+{}
			\innprod{w-T_{\nicefrac{1}{2\alpha}}w}{x-w}
		\right)
	\shortintertext{using Cauchy-Schwartz inequality,}
	{}\geq{} &
		4\alpha^2
		\|w-T_{\nicefrac{1}{2\alpha}}w\|
		\bigl(
			\|w-T_{\nicefrac{1}{2\alpha}}w\|
			{}-{}
			\|x-w\|
		\bigr)
	\shortintertext{the bound on \(\|x-w\|\),}
	{}\geq{} &
		2\alpha
		\|Rw\|
		\bigl(
			\|w-T_{\nicefrac{1}{2\alpha}}w\|
			{}-{}
			2\alpha c\|x-T_{\nicefrac{1}{2\alpha}}x\|
		\bigr)
	\shortintertext{the (reverse) triangular inequality,}
	{}\geq{} &
	\ifieee
			2\alpha
			\|Rw\|
			\Bigl(
				\hphantom{{}-{}}
				(1-2\alpha c)
				\|x-T_{\nicefrac{1}{2\alpha}}x\|
		\\
		&
			\hphantom{2\alpha\|Rw\|\Bigl(}
				-{}
				\|(\id-T_{\nicefrac{1}{2\alpha}})w-(\id-T_{\nicefrac{1}{2\alpha}})x\|
			\Bigr)
	\else
		2\alpha
		\|Rw\|
		\bigl(
			(1-2\alpha c)
			\|x-T_{\nicefrac{1}{2\alpha}}x\|
			{}-{}
			\|(\id-T_{\nicefrac{1}{2\alpha}})w-(\id-T_{\nicefrac{1}{2\alpha}})x\|
		\bigr)
	\fi
	\shortintertext{the nonexpansiveness of \(\id-T_{\nicefrac{1}{2\alpha}}\)}
	{}\geq{} &
		2\alpha
		\|Rw\|
		\bigl(
			\tfrac{1-2\alpha c}{2\alpha}
			\|Rx\|
			{}-{}
			\|w-x\|
		\bigr)
	\shortintertext{and again the bound on \(\|w-x\|\),}
	{}\geq{} &
		(1-4\alpha c)
		\|Rw\|
		\|Rx\|
	\end{align*}
	equating \(\sigma=1-4\alpha c\) the assert follows.
	\end{proof}
\end{thm}
Notice that if \(d=0\), then \(\rho=\|Rx\|^2\geq\sigma\|Rx\|^2\) for any \(\sigma\in[0,1)\), and therefore the line search condition \eqref{eq:LS} is always satisfied; in particular, the classical KM step \(x^+=Tx\) is always accepted regardless of the value of \(\sigma\).

Let us now observe how a GKM projection extends the classical KM depicted in \Cref{fig:Mann} and how the line search works.
In the following we use the notation of \Cref{thm:LS}, and for the sake of simplicity we consider \(\sigma=0\) in \eqref{eq:LS} and a FNE operator \(T\).
Suppose that the fixed point \(z\) and the points \(x\), \(Tx\), and \(w\) are as in \Cref{fig:SuperMann1}; due to firm nonexpansiveness, the image \(Tw\) of \(w\) is somewhere in both the orange circles.
We want to avoid the unfavorable situation depicted in \Cref{fig:SuperMann2}, where the couple \((w,Tw)\) generates a halfspace \(C_w\) that contains \(x\), \ie, such that \(\rho\leq 0\): in fact, with simple algebra it can be seen that \(\rho\leq 0\) iff \(Tw\) belongs to the dashed circle of \Cref{fig:SuperMann2}:
\begin{equation}\label{eq:Bxw}
	B_{x,w}
{}\coloneqq{}
	\set{\bar w}[
		\innprod{w-\bar w}{x-\bar w}
	{}\leq{}
		0
	].
\end{equation}
Since the dashed orange circle (in which \(Tw\) must lie) is simply the translation by a vector \(Tx-x\) of \(B_{x,w}\), both having diameter \(\tau\|d\|\), for sufficiently small \(\tau\) the two have empty intersection, meaning that \(\rho>0\) regardless of where \(Tw\) is.




	\section{The SuperMann scheme}
		\label{sec:SuperMann}
		In this section we introduce the \refSM{}, a special instance of the general framework of \Cref{alg:General} that employs GKM updates as safeguard \ref{K2}-steps.
While the global worst-case convergence properties of \emph{SuperMann} are the same as for the classical KM scheme, its asymptotic behavior is determined by how \emph{blind} \ref{K0}- and \emph{educated} \ref{K1}-updates are selected.
In \Cref{sec:Superlinear} we will characterize the ``quality'' of update directions and the mild requirements under which superlinear convergence rates are attained; in particular, \Cref{sec:Broyden} is dedicated to the analysis of quasi-Newton Broyden's directions.

\added{%
	The scheme follows the same philosophy of the general abstract framework.
	The main idea is globalizing a local method for solving the monotone equation \(Rx=0\), in such a way that when the iterates get close enough to a solution the fast convergence of the local method is automatically triggered.
	Approaching a solution is possible thanks to the generalized KM updates (step \ref{step:SuperMann:K2}), provided enough backtracking is performed, as ensured by \cref{prop:DriftedKMFejer,thm:LS}.
	When a basin of fast (\ie superlinear) attraction for the local method is reached, the (norm of) \(Rx\) will decrease more than linearly, and the condition triggering the \emph{educated} updates of step \ref{step:SuperMann:K1} (which is checked first) will be verified without performing any backtracking.
}%

To discuss its global and local convergence properties we stick to the same notation of the general framework of \Cref{alg:General}, denoting the sets of \emph{blind}, \emph{educated}, and \emph{safeguard} updates as \(K_0\), \(K_1\) and \(K_2\), respectively.

\begin{algorithm*}[tb]
	\algcaption[SuperMann scheme]{%
		\emph{SuperMann} scheme for solving \eqref{eq:Problem}, given an \avg{} operator \(T\) with residual \(R=\id-T\)%
	}%
	\label{alg:SuperMann}%
	\newcommand\keyfont[1]{\textsc{#1}}
\begin{tabularx}{\linewidth}{@{}l@{~~}X@{}}
		\keyfont{Require}%
	&
		\(x_0\in\HH\),~
		\(c_0,c_1,q\in[0,1)\),~
		\(\beta,\sigma\in(0,1)\),~
		\(\lambda\in(0,\nicefrac1\alpha)\).%
	\\
		\keyfont{Initialize}%
	&
		\(\eta_0=r_{\rm safe}=\|Rx_0\|\),~
		\(k=0\)%
\end{tabularx}
\begin{enumerate}[{label=\textbf{\arabic*.}},{ref=\textbf{\arabic*}}]
	\item\label{step:SuperMann:Initial}
		\keyfont{If} \(Rx_k=0\), \keyfont{then stop}.
	\item\label{step:SuperMann:d}
		Choose an update direction \(d_k\in\HH\)
	\item\label{step:SuperMann:K0}
		{\boldmath \((K_0)\)}~
		\keyfont{If} \(\|Rx_k\|\leq c_0\eta_k\),~
		\keyfont{then} set
		\(
			\eta_{k+1}
		{}={}
			\|Rx_k\|
		\),
		proceed with a \emph{blind update} \(x_{k+1}{={}}w_k{\coloneqq{}} x_k{+}d_k\) and go to step \ref{step:SuperMann:k++}.%
	\item
		Set \(\eta_{k+1}=\eta_k\)~ and ~\(\tau_k=1\).%
	\item\label{step:SuperMann:w}
		Let \(w_k=x_k+\tau_kd_k\).
		\begin{enumerate}[{label=\textbf{\theenumi(\alph*)}},{ref=\textbf{\theenumi(\alph*)}}]
			\item\label{step:SuperMann:K1}
				{\boldmath \((K_1)\)}~
				\keyfont{If} the \emph{safe condition} \(\|Rx_k\|\leq r_{\rm safe}\) holds and \(w_k\) is \emph{educated}:
				\[
					\|Rw_k\|
				{}\leq{}
					c_1\|Rx_k\|
				\]
				\keyfont{then} set \(x_{k+1}=w_k\), update
				\(
					r_{\rm safe}
				{}={}
					\|Rw_k\|
					{}+{}
					q^k
				\),
				and go to step \ref{step:SuperMann:k++}.
			\item\label{step:SuperMann:K2}
				{\boldmath \((K_2)\)}~
				\keyfont{If}~
				\(
					\rho_k
				{}\coloneqq{}
					\|Rw_k\|^2
				{}-{}
					2\alpha
					\innprod{Rw_k}{w_k-x_k}
				{}\geq{}
					\sigma\|Rw_k\|
					\mathrlap{
						\|Rx_k\|
					}
				\)
				
				\keyfont{then} set
				\[
					x_{k+1}
				{}={}
					\smash{
						x_k - \lambda\frac{\rho_k}{\|Rw_k\|^2}Rw_k
					}
				\]
				\keyfont{otherwise} set \(\tau_k\gets\beta\tau_k\) and go to step \ref{step:SuperMann:w}.%
		\end{enumerate}
	\item\label{step:SuperMann:k++}
		Set \(k\gets k+1\) and go to step \ref{step:SuperMann:Initial}.
\end{enumerate}
\end{algorithm*}

		\subsection{Global and linear convergence}
			To comply with \eqref{eq:D}, we impose the following requirement on the magnitude of the directions (see also \cref{rem:D}).
\begin{ass}\label{ass:d}
There exists a constant \(D\geq 0\) such that the directions \(\seq{d_k}\) in the \refSM{} satisfy
\begin{equation}\label{eq:SuperMann:dBound}
	\|d_k\|
{}\leq{}
	D\|Rx_k\|
\qquad
	\forall k\in\N.
\end{equation}
\end{ass}

\begin{thm}[{Global and linear convergence of the \refSM[]}]\label{thm:SuperMann}
Consider the iterates generated by the \refSM{} with \(\seq{d_k}\) selected so as to satisfy \Cref{ass:d}.
Then,
\begin{enumerate}
	\item\label{thm:SuperMann:QuasiFejer}
		\(\seq{x_k}\) is quasi-Fejér monotone with respect to \(\fix T\);
	\item\label{thm:SuperMann:Tau_k}
		\(\tau_k=1\) if \(d_k=0\), and
		\(
				\tau_k
			{}\geq{}
			\min\set{
				\beta\frac{1-\sigma}{4\alpha D}
				{},{}
				1
			}
		\)
		otherwise.
	\item\label{thm:SuperMann:ResStr}
		\(Rx_k\to 0\) with \(\seq{\|Rx_k\|}\in\ell^2\);
	\item\label{thm:SuperMann:xWeak}
		\(\seq{x_k}\) converges weakly to a point \(x_\star\in\fix T\);
	\item\label{thm:SuperMann:K0}
		if \(c_0>0\) the number of blind updates at step \ref{step:SuperMann:K0} is infinite.
\end{enumerate}
Moreover, if \(\seq{x_k}\) converges strongly to a point \(x_\star\) (this being true if \(\HH\) is finite dimensional) at which \(R\) is metrically subregular, then
\begin{enumerate}[resume]
	\item\label{thm:SuperMann:Linear}
	\(\seq{x_k}\) and \(\seq{Rx_k}\) are \(R\)-linearly convergent.
\end{enumerate}
\begin{proof}
See \Cref{proof:thm:SuperMann}.
\end{proof}
\end{thm}

		\subsection{Superlinear convergence}
			\label{sec:Superlinear}
			Though global convercence of the \refSM[] is independent of the choice of the directions \(d_k\), its performance and tail convergence surely does.
We characterize the \emph{quality} of the directions \(d_k\) in terms of the following definition.
\begin{defin}[{Superlinear directions for the \refSM[]}]\label{defin:dSuperlinear}
Relative to the sequence \(\seq{x_k}\) generated by the \refSM[], we say that \(\seq{d_k}\subset\HH\) are \DEF{superlinear directions} if the following limit holds
\[
	\lim_{k\to\infty}{
		\frac{
			\|R(x_k+d_k)\|
		}{
			\|Rx_k\|
		}
	}
{}={}
	0.
\]
\end{defin}
\begin{rem}\label{rem:dSuperlinear}
\Cref{defin:dSuperlinear} makes no mention of a limit point \(x_\star\) of the sequence \(\seq{x_k}\), differently from the definition in \cite{facchinei2003finite} which instead requires \(\frac{\|x_k+d_k-x_\star\|}{\|x_k-x_\star\|}\) to be vanishing with no mention of \(R\).
Due to \(2\alpha\)-Lipschitz continuity of \(R\), whenever the directions \(d_k\) are bounded as in \eqref{eq:SuperMann:dBound} we have
\[
	\frac{\|R(x_k+d_k)\|}{\|Rx_k\|}
{}\leq{}
	2\alpha D
	\frac{\|x_k+d_k-x_\star\|}{\|d_k\|}.
\]
Invoking \cite[Lem. 7.5.7]{facchinei2003finite} it follows that \Cref{defin:dSuperlinear} is implied by the one in \cite{facchinei2003finite} and is therefore more general.
\end{rem}
\begin{thm}\label{thm:dSuperlinear}
Consider the iterates generated by the \refSM{} with either \(c_0>0\) or \(c_1>0\), and with \(\seq{d_k}\) being superlinear directions as in \Cref{defin:dSuperlinear}.
Then,
\begin{enumerate}
	\item\label{thm:dSuperlinear:tau}
		eventually, stepsize \(\tau_k=1\) is always accepted and \emph{safeguard} updates \(K_2\) are deactivated (\ie, the scheme reduces to the local method \(x_{k+1}=x_k+d_k\));
	\item\label{thm:dSuperlinear:Rx}
		\(\seq{Rx_k}\) converges \(Q\)-superlinearly;
	\item\label{thm:dSuperlinear:x}
		if the directions \(d_k\) satisfy \Cref{ass:d}, then \(\seq{x_k}\) converges \(R\)-super\-linearly;
	\item\label{thm:dSuperlinear:K0}
		if \(c_0>0\), then the complement of \(K_0\) is finite.
\end{enumerate}
\begin{proof}
See \Cref{proof:thm:dSuperlinear}.
\end{proof}
\end{thm}

\Cref{thm:dSuperlinear} shows that when the directions \(d_k\) are \emph{good}, then eventually the \refSM[] reduces to the local method \(x_{k+1}=x_k+d_k\) and consequently inherits its local convergence properties.
The following result specializes to the choice of semismooth Newton directions.
\begin{cor}[Superlinear convergence for semismooth Newton directions]\label{thm:Semismooth}
Suppose that \(\HH\) is finite dimensional, and that \(R\) is semismooth.
Consider the iterates generated by the \refSM{} with either \(c_0>0\) or \(c_1>0\) and directions \(d_k\) chosen as solutions of
\begin{equation}\label{eq:dSemismooth}
	(G_k+\mu_k\id)d_k=-Rx_k
\quad\text{for some }
	G_k\in\partial Rx_k
\end{equation}
where \(\partial R\) denotes the \DEF{Clarke generalized Jacobian} of \(R\) and \(0\leq \mu_k\to 0\).
Suppose that the sequence \(\seq{x_k}\) converges to a point \(x_\star\) at which all the elements in \(\partial R\) are nonsingular.

Then, \(\seq{d_k}\) are superlinear directions as in \Cref{defin:dSuperlinear}, and in particular all the claims of \Cref{thm:dSuperlinear} hold.
\begin{proof}
Any \(G_k\in\partial R\) is positive semidefinite due to the monotonicity of \(R\), and therefore \(d_k\) as in \eqref{eq:dSemismooth} is well defined for any \(\mu_k>0\).
The bound \eqref{eq:SuperMann:dBound} holds due to \cite[Thm. 7.5.2]{facchinei2003finite}.
Moreover,
\[
	\frac{\|Rx_k+G_kd_k\|}{\|d_k\|}
{}={}
	\mu_k
{}\to{}
	0
\]
as \(k\to\infty\), and the proof follows invoking \cite[Thm. 7.5.8(a)]{facchinei2003finite} and \cref{rem:dSuperlinear}.
\end{proof}
\end{cor}
Notice that since \(\partial R=\id-\partial T\), nonsingularity of the elements in \(\partial R(x_\star)\) is equivalent to having \(\|G\|<1\) for all \(G\in\partial T(x_\star)\), \ie that \(T\) is a local contraction around \(x_\star\).

Despite the favorable properties of semismooth Newton methods, in this paper we are oriented towards choices of directions that (1) are defined for any nonexpansive mapping, regardless of the (generalized) first-order properties, and that (2) require exactly the same oracle information as the original KM scheme.
This motivates the investigation of quasi-Newton directions, whose superlinear behavior is based on the classical Dennis-Moré criterion, which we provide next.
We first recall the notions of \emph{semi-} and \emph{strict differentiability}.
\begin{defin}\label{defin:GenDiff}
We say that \(\func R\HH\HH\) is
\begin{enumerate}
\item
	\DEF{strictly differentiable} at \(\bar x\) if it is differentiable there with \(\jac R{\bar x}\) satisfying
	\begin{equation}\label{eq:StrictDiff}
		\lim_{%
			\substack{(y,z)\to(\bar x,\bar x)\\y\neq z}
		}{
			\frac{
				\left\|
					Ry-Rx-\jac R{\bar x}(y-x)
				\right\|
			}{
				\|y-x\|
			}
		}
	{}={}
		0;
	\end{equation}
\item
	\DEF{semidifferentiable} at \(\bar x\) if there exists a continuous and positively homogeneous function
	\(
		\func{\Bjac R{\bar x}}{\HH}{\HH}
	\),
	called the \DEF{semiderivative} of \(R\) at \(\bar x\), such that
	\[
		Rx
	{}={}
		R\bar x
		{}+{}
		\Bjac R{\bar x}[x-\bar x]
		{}+{}
		o(\|x-\bar x\|);
	\]
\item
	\DEF{calmly semidifferentiable} at \(\bar x\) if there exists a neighborhood \(U_{\bar x}\) of \(\bar x\) in which \(R\) is semidifferentiable and such that for all \(w\in\HH\) with \(\|w\|=1\) the function \(U_{\bar x}\ni x\mapsto\Bjac Rx[w]\) is Lipschitz continuous at \(\bar x\).
\end{enumerate}
\end{defin}
There is a slight ambiguity in the literature, whereas \emph{strict} differentiability is sometimes referred to rather as \emph{strong} differentiability \cite{ip1992local,pang1990newtons}.
We choose to stick the proposed terminology, following \cite{rockafellar2011variational}.
Semidifferentiability is clearly a milder property than differentiability in that the mapping \(\Bjac R{\bar x}\) needs not be linear.
More precisely, since the residual \(R\) of a nonexpansive operator is (globally) Lipschitz continuous, then semidifferentiability is equivalent to directional differentiability \cite[Prop. 3.1.3]{facchinei2003finite} and the semiderivative is sometimes called \emph{\(B\)-derivative} \cite{ip1992local,facchinei2003finite}.
The three concepts in \Cref{defin:GenDiff} are related as \emph{(iii)} \(\Rightarrow\) \emph{(i)} \(\Rightarrow\) \emph{(ii)} \cite[Thm. 2]{pang1990newtons} and neither requires the existence of the (classical) Jacobian around \(\bar x\).

\begin{thm}[Dennis-Moré criterion for superlinear convergence]\label{thm:Local}
Consider the iterates generated by the \refSM{} and suppose that \(\seq{x_k}\) converges strongly to a point \(x_\star\) at which \(R\) is strictly differentiable.
Suppose further that the update directions \(\seq{d_k}\) satisfy \Cref{ass:d} and the Dennis-Moré condition
\begin{equation}\label{eq:DM}
	\lim_{k\to\infty}{
		\frac{
			\|Rx_k+\jac R{x_\star}d_k\|
		}{
			\|d_k\|
		}
	}
{}={}
	0.
\end{equation}
Then, the directions \(d_k\) are superlinear as in \Cref{defin:dSuperlinear}.
In particular, all the claims of \Cref{thm:dSuperlinear} hold.
\begin{proof}
\ifieee
	\begin{align*}
		0
	{{}\overrel{\eqref{eq:DM}}} &
		\!\lim_{k\to\infty}{
			\frac{
				\bigl\|
					Rx_k
					{}+{}
					\jac R{x_\star}d_k
					{}+{}
					R(x_k+d_k)
					{}-{}
					R(x_k+d_k)
				\bigr\|
			}{
				\|d_k\|
			}
		}
	\\
	{=} &
		\!\lim_{k\to\infty}{
			\frac{
				\bigl\|
					R(x_k+d_k)
				\bigr\|
			}{
				\|d_k\|
			}
		}
	~{}\overrel[\geq]{\eqref{eq:SuperMann:dBound}}{}~
		\frac 1D
		\lim_{k\to\infty}{
			\frac{
				\bigl\|
					R(x_k+d_k)
				\bigr\|
			}{
				\|Rx_k\|
			}
		}
	\end{align*}
\else
	\begin{align*}
		0
	{}\overrel{\eqref{eq:DM}}{} &
		\lim_{k\to\infty}{
			\frac{
				\bigl\|
					Rx_k
					{}+{}
					\jac R{x_\star}d_k
					{}+{}
					\bigl(
						R(x_k+d_k)
						{}-{}
						R(x_k+d_k)
					\bigr)
				\bigr\|
			}{
				\|d_k\|
			}
		}
	{}={}
		\lim_{k\to\infty}{
			\frac{
				\bigl\|
					R(x_k+d_k)
				\bigr\|
			}{
				\|d_k\|
			}
		}
	\\
	{}\overrel[\geq]{\eqref{eq:SuperMann:dBound}}{} &
		\frac 1D
		\lim_{k\to\infty}{
			\frac{
				\bigl\|
					R(x_k+d_k)
				\bigr\|
			}{
				\|Rx_k\|
			}
		}
	\end{align*}
\fi
where in the second equality we used strict differentiability of \(R\) at \(x_\star\).
\end{proof}
\end{thm}

		\subsection{A modified Broyden's direction scheme}
			\label{sec:Broyden}
			In practical application the Hilbert space \(\HH\) is finite dimensional, and consequently it can be identified with \(\R^n\).
Then, the computation of quasi-Newton directions \(d_k\) in the \refSM[] amounts to selecting
\begin{subequations}\label{eq:QN}
\begin{equation}\label{eq:dQN}
	d_k=-B_k^{-1}Rx_k,
\end{equation}
where \(B_k\in\R^{n\times n}\) are recursively defined by low-rank updates satisfying a \emph{secant condition}, starting from an invertible matrix \(B_0\).
The most popular quasi-Newton scheme is the 2-rank BFGS formula, which also enforces symmetricity.
As such, BFGS is well performing only when the Jacobian at the solution \(JRx_\star\) possesses this property, a requirement that is not met by the residual \(R\) of generic nonexpansive mappings.

For this reason we consider Broyden's method as a universal alternative.
We adopt Powell's modification \cite{powell1970numerical} to enforce nonsingularity and make \eqref{eq:dQN} well defined: for a fixed parameter \(\bar\vartheta\in(0,1)\), matrices \(B_k\) are recursively defined as
\begin{equation}\label{eq:Broyden}
	\textstyle
	B_{k+1}
{}={}
	B_k
	{}+{}	
	\frac{
		1
	}{
		\|s_k\|_2^2
	}
	\bigl(\tilde y_k-B_ks_k\bigr)
	\trans{s_k}
\end{equation}
where for
\(
	\gamma_k
{}\coloneqq{}
	\frac{\innprod{B_k^{-1}y_k}{s_k}_{\!_2}}{\|s_k\|_{_2}^2}
\)
we have defined
\begin{equation}\label{eq:PowellTheta}
	\begin{cases}
		s_k & {}= w_k-x_k \\
		y_k & {}= Rw_k-Rx_k \\
		\tilde y_k & {}= (1-\vartheta_k)B_ks_k + \vartheta_ky_k
	\end{cases}
\ifieee
	\!\vartheta_k
{}\coloneqq{}
	\begin{cases}[l@{~\text{if }}l]
			1
		&
			|\gamma_k| {\,\geq\,} \bar\vartheta
		\\
			\frac{
				1-\sign(\gamma_k)
				\bar\vartheta
			}{
				1-\gamma_k
			}
		&
			|\gamma_k| {\,<\,} \bar\vartheta
	\end{cases}
\else
	\quad\text{and}\quad
	\vartheta_k
{}\coloneqq{}
	\begin{cases}[l@{~~\text{if }}l]
			1
		&
			|\gamma_k| \geq \bar\vartheta
		\\
			\frac{
				1-\sign(\gamma_k)
				\bar\vartheta
			}{
				1-\gamma_k
			}
		&
			|\gamma_k| < \bar\vartheta
	\end{cases}
\fi
\end{equation}
with the convention \(\sign 0=1\).
Letting \(H_k\coloneqq B_k^{-1}\) and using the Sherman-Morrison identity, the inverse of \(B_k\) is given by
\begin{equation}\label{eq:BroydenInverse}
	\textstyle
	H_{k+1}
{}={}
	H_k
	{}+{}	
	\frac{
		1
	}{
		\innprod{H_k\tilde y_k}{s_k}_{\!_2}
	}
	\bigl(s_k-H_k\tilde y_k\bigr)
	\bigl(\trans{s_k}H_k\bigr).
\end{equation}
Consequently, there is no need to compute and store the matrices \(B_k\) and we can directly operate with their inverses \(H_k\).
\end{subequations}
\begin{thm}[{Superlinear convergence of the \refSM[] with Broyden's directions}]\label{thm:SuperMann:Broyden}
Suppose that \(\HH\) is finite dimensional.
Consider the sequence \(\seq{x_k}\) generated by the \refSM{}, \(\seq{d_k}\) being selected with the modified Broyden's scheme \eqref{eq:QN} for some \(\bar\vartheta\in(0,1)\).

Suppose that \(\seq{H_k}\) remains bounded, and that \(R\) is calmly semidifferentiable and metrically subregular at the limit \(x_\star\) of \(\seq{x_k}\).
Then, \(\seq{d_k}\) satisfies the Dennis-Moré condition \eqref{eq:DM}.
In particular, all the claims of \Cref{thm:Local} hold.
\begin{proof}
See \Cref{proof:thm:SuperMann:Broyden}.
\end{proof}
\end{thm}
\added{%
	\begin{rem}\label{rem:D}%
		It follows from \Cref{thm:SuperMann:xWeak} that the \refSM[] is globally convergent as long as \(\|d_k\|\leq D\|Rx_k\|\) for some constant \(D\).
		To enforce it we may select a (large) constant \(D>0\) and as a possible choice truncate \(d_k\gets D\frac{\|Rx_k\|}{\|d_k\|}d_k\) whenever \(d_k\) does not satisfy \eqref{eq:SuperMann:dBound}.
	\end{rem}
}%
Let us observe that in order to achieve superlinear convergence \emph{the \refSM[] does not require nonsingularity of the Jacobian at the solution}.
This is the standard requirement for asymptotic properties of quasi-Newton schemes, which is needed to show first that the method converges at least linearly.
\cite{artacho2014local} generalizes this property invoking the concepts of (strong) metric (sub)regularity (see also \cite{dontchev2004regularity} for an extensive review on these properties).
However, if \(R\) is strictly differentiable at \(x_\star\), then strong subregularity, regularity and strong regularity are equivalent to injectivity, surjectivity and invertibility of \(\jac R{x_\star}\), respectively, these conditions being all equivalent for mappings \(\HH\to\HH\) with \(\HH\) finite dimensional.
In particular, contrary to the \refSM[] standard approaches require the solution \(x_\star\) at least to be isolated%
\ifieee
	.
\else
	, a property that rules out many interesting applications (cf. \cref{sec:SCS}).
\fi

			\subsubsection*{Restarted (modified) Broyden's scheme}
				\label{sec:RBroyden}
				\added{%
	Broyden's scheme requires storing and operating with \(n\times n\) matrices, where \(n\) is the dimension of the optimization variable, and is consequently feasible in practice only for small problems.
	Alternatively, one can restrict Broyden's update rule \eqref{eq:BroydenInverse} to only the most recent pairs of vectors \((s_i,y_i)\).
	As detailed in \Cref{alg:RBroyden}, this can be done by keeping track of the last vectors \(s_i\) and some auxiliary vectors
	\(
		\tilde s_i
	{}={}
		\tfrac{
			s_i-H_i\tilde y_i
		}{
			\innprod{s_i}{H_i\tilde y_i}_2
		}
	\).
	These are stored in some buffers \(S\) and \(\tilde S\), which are initially empty and can contain up to \(m\) vectors.
	The \emph{memory} \(m\) is a small integer typically between \(3\) and \(20\); when the memory is full, the buffers are emptied and Broyden's scheme is \emph{restarted}.
	The choice of a \emph{restarted} rather than a \emph{limited-memory} variant
	obviates the need of a nested for-loop to account for Powell's modification.
}%
\begin{algorithm}
	\begin{algorithmic}[1]%
	\itemindent=-.5\algorithmicindent%
	\item[]\hspace*{-\algorithmicindent}\textfillwidthof[l]{\textbf{Output: }}{\textbf{Input: }}%
		old buffers \(S,\tilde S\);~
		new pair \((s,y)\);~
		current \(Rx\)%
	\item[]\hspace*{-\algorithmicindent}\textbf{Output: }%
		new buffers \(S,\tilde S\);~
		update direction \(d\)%
	\STATE%
		\(d\gets-Rx\),~
		\(\tilde s\gets y\)%
	\FOR{ \(i=1\ldots \#S\) }%
		\item[]\(
			\tilde s
		{}\gets{}
			\tilde s
			{}+{}
			\innprod{s_i}{\tilde s}_{\!_2}\tilde s_i
		\),~
		\(
			d
		{}\gets{}
			d
			{}+{}
			\innprod{s_i}{d}_{\!_2}\tilde s_i
		\)%
	\ENDFOR{}%
	\STATE%
		compute \(\vartheta\) as in \eqref{eq:PowellTheta} with
		\(
			\gamma
		{}={}
			\tfrac{1}{\|s\|_2^2}\innprod{\tilde s}{s}_{\!_2}
		\)%
	\STATE\(
		\tilde s
	{}\gets{}
		\frac{
			\vartheta
		}{
			(1-\vartheta+\vartheta\gamma)\|s\|_2^2
		}(s-\tilde s)
	\),~
	\(
		d
	{}\gets{}
		d
		{}+{}
		\innprod{s}{d}_{\!_2}\tilde s
	\)%
	\STATE%
		\textbf{if } \(\#S = m\) \textbf{ then }
			\(S,\tilde S\gets [\,]\)
		\textbf{ else }
			\(S\gets [S,s]\), \(\tilde S\gets [\tilde S,\tilde s]\)%
\end{algorithmic}
	\caption{%
		Restarted Broyden's scheme with memory \(m\) using Powell's modification%
	}%
	\label{alg:RBroyden}%
\end{algorithm}
		\subsection{Parameters selection in \emph{SuperMann}}
			\label{sec:Parameters}
			As shown in \Cref{thm:dSuperlinear}, the \refSM[] makes sense as long as either \(c_0>0\) or \(c_1>0\); indeed, \emph{safeguard} \ref{K2}-steps are only needed for globalization, while it is \emph{blind} \ref{K0}- and \emph{educated} \ref{K1}-steps that exploit the quality of the directions \(d_k\).
Evidently, \ref{K1}-updates are more reliable than \ref{K0}-updates in that they take into account the residual of the candidate next point.
As such, it is advisable to select \(c_1\) close to \(1\) and use small values of \(c_0\) if more conservatism and robustness are desired.
To further favor \ref{K1}-updates, the parameter \(q\) used for updating the safeguard \(r_{\rm safe}\) at step \ref{step:SuperMann:K1} may be also chosen very close to \(1\).

As to \emph{safeguard} \ref{K2}-steps, a small value of \(\sigma\) makes condition \eqref{eq:LS} easier to satisfy and results in fewer backtrackings; the averaging factor \(\lambda\) may be chosen equal to 1 whenever possible, \ie if \(\alpha\lneqq 1\) (which is the typical case when, \eg \(T\) comes from splitting schemes in convex optimization), or any close value otherwise.
\added{%
	In the simulations of \Cref{sec:Simulations} we used \(c_0=c_1=q=0.99\), \(\sigma=0.1\), \(\lambda=1\) and \(\beta=\nicefrac 12\).
	For a matter of scaling, we multiplied the summable term \(q^k\) by \(\|Rx^0\|\) in updating the parameter \(r_{\rm safe}\) at step \cref{step:SuperMann:K1}.
	The directions were computed according to the restarted modified Broyden's scheme (\cref{alg:RBroyden}) with memory \(m=20\) and \(\bar\vartheta=0.2\); we applied the truncation rule as in \Cref{rem:D} with \(D=10^4\).
	We also imposed a maximum of \(8\) backtrackings after which a nominal V\~u-Condat iteration would be executed.
}%

		\subsection{Comparisons with other methods}
			\subsubsection{Hybrid global and local phase algorithms}
				\emph{Blind} \(K_0\)-updates in the \refSM[] are inspired from \cite[Alg. 1]{chen1999proximal}, and so is the notation \(K_0=\set{k_0,k_1,\ldots}\).

\emph{Educated} \(K_1\)- and \emph{safeguard} \(K_2\)-updates instead play the role of \emph{inner}- and \emph{outer-phases} in the general algorithmic framework described in \cite[§5.3]{izmailov2014newton} for finding a zero of a candidate merit function \(\varphi\) (e.g. \(\varphi(x)=\frac12\|Rx\|^2\) in our case).
Differently from \cite[Alg. 5.16]{izmailov2014newton} where all previous inner-phase iterations are discarded as soon as the required sufficient decrease is not met, the \refSM[] allows for an alternation of phases that eventually stabilizes on the fast local one, provided the solution is sufficiently regular.
Our scheme is more in the flavor of \cite[Alg. 5.19]{izmailov2014newton}, although with less conservative requirements for triggering \emph{inner} \(K_1\)-updates (\(\varphi(x_{k+1})\) is here compared with \(\varphi(x_k)\), whereas in the cited scheme with the smallest past value).

			\subsubsection{Inexact Newton methods for monotone equations}
				\label{sec:InexactNewton}
				The GKM updates are closely related to the extra-gradient steps described in \cite[Alg. 2.1]{solodov1998globally}.
This work introduces an inexact Newton algorithm for solving systems of continuous monotone equations \(Rx=0\), where \(\id-R\) needs not be nonexpansive.
At a given point \(x\), first a direction \(d\) is computed as (possibly approximate) solution of \(Gd=-Rx\), where \(G\) is some positive definite matrix; then, an intermediate point \(w=x+\tau d\) is retrieved with a line search on \(\tau\) that ensures the condition
\begin{equation}\label{eq:InexactNewton_LS}
	\|Rw\|^2
	{}-{}
	\innprod{Rw}{x-Tw}
{}\leq{}
	-\sigma\tau
	\|d\|^2
\end{equation}
for some \(\sigma>0\); here, we defined \(T\coloneqq\id-R\) to highlight the symmetry with \eqref{eq:C_x}.
Finally, the new iterate is given by
\ifieee
	\(x^+=\proj_{H_w}x\), where
\fi
\begin{equation}
\label{eq:InexactNewton_x}
\ifieee\else
		x^+=\proj_{H_w}x
	\quad\text{where}\quad
\fi
	H_w
{}\coloneqq{}
	\set{z\in\HH}[
		\|Rw\|^2
		{}-{}
		\innprod{Rw}{z-Tw}
	{}\geq{}
		0
	].
\end{equation}

\begin{figure*}
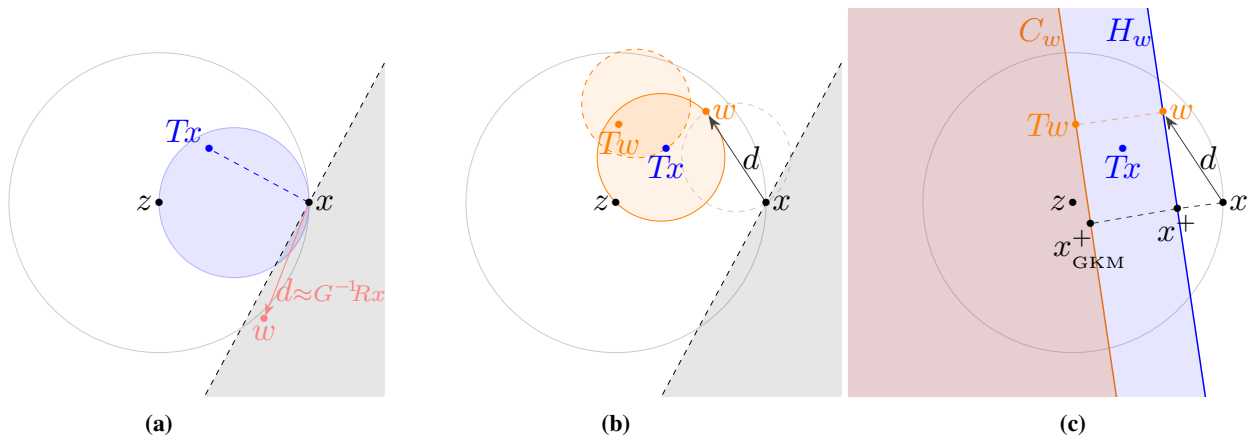

	\setlength\fboxsep{0pt}%
	\begin{subfigure}{0.33\textwidth}%
		\vspace*{0pt}%
		{{%
		\pgfkeys{/pgf/images/include external/.code={\includegraphics[width=\linewidth]{#width=\linewidth}}}%
		\tikzsetnextfilename{InexactNewton1}%
		\def\radius{2}
\def\xmax{3}
\begin{tikzpicture}%
	\pgfplotsset{%
		mystyle/.style={%
			ticks=none,
			axis lines=none,
			xmin={-\xmax},
			xmax={\xmax},
			ymin={-\xmax},
			ymax={\xmax},
			enlargelimits=false,
			scale only axis,
		},%
	}%
	\tikzset{%
		textnode/.style={%
			scale=1.75,
			inner sep=0pt,
			outer sep=2pt,
		},%
		nosep/.style={%
			inner sep=0pt,
			outer sep=0pt,
		},%
		bad/.style={%
			red!50,
			ultra thin,
		},%
		shortarrow/.style={%
			-{Stealth[scale=1.5]},
			shorten >=2pt,
			shorten <=1pt,
		},%
	}%
	\begin{axis}[mystyle]
		\node[nosep] (top_left) at (axis cs:-\xmax,\xmax) {} ;
		\node[nosep] (top_right) at (axis cs:\xmax,\xmax) {} ;
		\node[nosep] (bottom_left) at (axis cs:-\xmax,-\xmax) {} ;
		\node[nosep] (bottom_right) at (axis cs:\xmax,-\xmax) {} ;
		\node[nosep] (z) at (axis cs:0,0) {};
		\node[nosep] (x) at (axis cs:\radius,0) {};
		\node[nosep] (Tx) at (axis cs:\radius/3,{1.25*\radius/3}) {};
		\node[nosep, bad] (w_bad) at (axis cs:{0.7*\radius},{-0.9*\radius}) {};
		\node (circle_z) at (z) [
			draw,
			very thin,
			black!20,
			circle through=(x),
		] {};
		\node (xOz) at ($(x)!0.5!(z)$) {};
		\node (circle_Tx) at (xOz) [
			draw,
			thin,
			blue!30,
			fill=blue,
			fill opacity=0.1,
			circle through=(x),
		] {};
		\node[nosep] (H_x1) at ($(x)!5!90:(Tx)$) {};
		\node[nosep] (H_x2) at ($(x)!5!-90:(Tx)$) {};
		\draw[
			black,
			thin,
			dashed,
		] (H_x1) -- (H_x2);
		\fill[
			black,
			opacity=0.1,
		] (H_x1.center) -- (H_x2.center) -- (bottom_right.center) -- cycle ;
		\draw[
			blue,
			ultra thin,
			dashed,
		] (x) -- (Tx);
		\draw[
			bad,
			shortarrow,
		] (x) -- node[textnode, bad, pos=0.75, anchor=west] {\hspace{-2pt}\(d{\scriptstyle{\approx} G^{{-}\!1}\!\!Rx}\)} (w_bad);
		\node at (z) {\textbullet};
		\node at (x) {\textbullet};
		\node[blue] at (Tx) {\textbullet};
		\node[bad] at (w_bad) {\textbullet};
		
		\node[textnode] at (z) [anchor=east] {\(z\)};
		\node[textnode] at (x) [anchor=west] {\(x\)};
		\node[textnode, blue] at (Tx) [anchor=south east] {\(T\!x\)};
		\node[textnode, bad] at (w_bad) [anchor=north] {\(w\)};
	\end{axis}%
\end{tikzpicture}%
%
	}}
		\caption{}%
		\label{fig:InexactNewton:d}
	\end{subfigure}%
	\hfill
	\begin{subfigure}{0.33\textwidth}%
		\vspace*{0pt}%
		{{%
		\pgfkeys{/pgf/images/include external/.code={\includegraphics[width=\linewidth]{#width=\linewidth}}}%
		\tikzsetnextfilename{InexactNewton2}%
		\def\radius{2}
\def\xmax{3}
\begin{tikzpicture}%
	\pgfplotsset{%
		mystyle/.style={%
			ticks=none,
			axis lines=none,
			xmin={-\xmax},
			xmax={\xmax},
			ymin={-\xmax},
			ymax={\xmax},
			enlargelimits=false,
			scale only axis,
		},%
	}%
	\tikzset{%
		textnode/.style={%
			scale=1.75,
			inner sep=0pt,
			outer sep=2pt,
		},%
		nosep/.style={%
			inner sep=0pt,
			outer sep=0pt,
		},%
		bad/.style={%
			red!50,
			ultra thin,
		},%
		shortarrow/.style={%
			-{Stealth[scale=1.5]},
			shorten >=2pt,
			shorten <=1pt,
		},%
	}%
	\begin{axis}[mystyle]
		\node[nosep] (top_left) at (axis cs:-\xmax,\xmax) {} ;
		\node[nosep] (top_right) at (axis cs:\xmax,\xmax) {} ;
		\node[nosep] (bottom_left) at (axis cs:-\xmax,-\xmax) {} ;
		\node[nosep] (bottom_right) at (axis cs:\xmax,-\xmax) {} ;
		\node[nosep] (z) at (axis cs:0,0) {};
		\node[nosep] (x) at (axis cs:{\radius},{0}) {};
		\node[nosep] (Tx) at (axis cs:{\radius/3},{1.25*\radius/3}) {};
		\node[nosep, orange] (w) at (axis cs:{0.6*\radius},{0.7*\radius}) {};
		\node[nosep, orange] (Tw) at (axis cs:{0.02*\radius},{0.6*\radius}) {};
		\node (circle_z) at (z) [
			draw,
			very thin,
			black!20,
			circle through=(x),
		] {};
		\node (wOz) at ($(w)!0.5!(z)$) {};
		\node at (wOz) [
			draw,
			orange,
			fill=orange,
			fill opacity=0.1,
			name path=circle_wOz,
			circle through=(z),
		] {};
		\node (xOw) at ($(x)!0.5!(w)$) {};
		\node (OTw) at ($(xOw)+(Tx)-(x)$) {};
		\node (circle_xOw) at (xOw) [
			draw,
			dashed,
			very thin,
			black!20,
			circle through=(x),
		] {};
		\draw[
			orange,
			fill=orange,
			fill opacity=0.1,
			dashed,
		] \getDistance{x}{w} (OTw) circle (\n1);
		\node[nosep] (H_x1) at ($(x)!5!90:(Tx)$) {};
		\node[nosep] (H_x2) at ($(x)!5!-90:(Tx)$) {};
		\draw[
			black,
			thin,
			dashed,
		] (H_x1) -- (H_x2);
		\fill[
			black,
			opacity=0.1,
		] (H_x1.center) -- (H_x2.center) -- (bottom_right.center) -- cycle ;
		\draw[
			shortarrow,
			black!75,
		] (x) -- node[textnode, black, pos=0.5, anchor=west] {\(d\)} (w);
		\node at (z) {\textbullet};
		\node at (x) {\textbullet};
		\node[blue] at (Tx) {\textbullet};
		\node[orange] at (w) {\textbullet};
		\node[orange] at (Tw) {\textbullet};
		
		\node[textnode] at (z) [anchor=east] {\(z\)};
		\node[textnode] at (x) [anchor=west] {\(x\)};
		\node[textnode, blue] at (Tx) [anchor=north] {\(T\!x\)};
		\node[textnode, orange] at (w) [anchor=west] {\(w\)};
		\node[textnode, orange] at (Tw) [anchor=north] {\(T\!w\)};
	\end{axis}%
\end{tikzpicture}%
%
	}}%
		\caption{}%
	\end{subfigure}%
	\hfill
	\begin{subfigure}{0.33\textwidth}%
		\vspace*{0pt}%
		{{%
		\pgfkeys{/pgf/images/include external/.code={\includegraphics[width=\linewidth]{#width=\linewidth}}}%
		\tikzsetnextfilename{InexactNewton3}%
		\def\radius{2}
\def\xmax{3}
\begin{tikzpicture}%
	\pgfplotsset{%
		mystyle/.style={%
			ticks=none,
			axis lines=none,
			xmin={-\xmax},
			xmax={\xmax},
			ymin={-\xmax},
			ymax={\xmax},
			enlargelimits=false,
			scale only axis,
		},%
	}%
	\tikzset{%
		textnode/.style={%
			scale=1.75,
			inner sep=0pt,
			outer sep=2pt,
		},%
		nosep/.style={%
			inner sep=0pt,
			outer sep=0pt,
		},%
		bad/.style={%
			red!50,
			ultra thin,
		},%
		shortarrow/.style={%
			-{Stealth[scale=1.5]},
			shorten >=2pt,
			shorten <=1pt,
		},%
	}%
	\begin{axis}[mystyle]
		\node[nosep] (top_left) at (axis cs:-\xmax,\xmax) {} ;
		\node[nosep] (top_right) at (axis cs:\xmax,\xmax) {} ;
		\node[nosep] (bottom_left) at (axis cs:-\xmax,-\xmax) {} ;
		\node[nosep] (bottom_right) at (axis cs:-\xmax,\xmax) {} ;
		\node[nosep] (z) at (axis cs:0,0) {};
		\node[nosep] (x) at (axis cs:{\radius},{0}) {};
		\node[nosep] (Tx) at (axis cs:{\radius/3},{1.25*\radius/3}) {};
		\node[nosep, orange] (w) at (axis cs:{0.6*\radius},{0.7*\radius}) {};
		\node[nosep, orange] (Tw) at (axis cs:{0.02*\radius},{0.6*\radius}) {};
		\node (circle_z) at (z) [
			draw,
			very thin,
			black!20,
			circle through=(x),
		] {};
		\node[nosep] (C_w1) at ($(Tw)!5!90:(w)$) {};
		\node[nosep] (C_w2) at ($(Tw)!5!-90:(w)$) {};
		\draw[
			orange,
			thick,
		] (C_w1) -- (C_w2);
		\draw[
			orange!50,
			very thin,
			dashed,
		] (w) -- (Tw);
		\fill[
			orange,
			opacity=0.2,
		] (C_w1.center) -- (C_w2.center) -- (bottom_left.center) -- (top_left.center) -- cycle ;
		\node[nosep] (GKM) at ($(C_w1)!(x)!(C_w2)$) {};
		\node[textnode, orange, anchor=north east, inner sep=1pt] at (axis cs:{0},{\xmax}) {\(C_w\)} ;
		\node[nosep] (H_w1) at ($(w)!5!90:(Tw)$) {};
		\node[nosep] (H_w2) at ($(w)!5!-90:(Tw)$) {};
		\draw[
			blue,
			thick,
		] (H_w1) -- (H_w2);
		\fill[
			blue,
			opacity=0.1,
		] (H_w1.center) -- (H_w2.center) -- (top_left.center) -- (bottom_left.center) -- cycle ;
		\node[nosep] (IS) at ($(H_w1)!(x)!(H_w2)$) {};
		\node[textnode, blue, anchor=north west, inner sep=1pt] at (axis cs:{0.15*\radius},{\xmax}) {\(H_w\)} ;
		\draw[
			shortarrow,
			black!75,
		] (x) -- node[textnode, black, pos=0.5, anchor=west] {\(d\)} (w);
		\draw[
			ultra thin,
			black!75,
			dashed,
		] (x) -- (IS);
		\draw[
			ultra thin,
			black!75,
			dashed,
		] (IS) -- (GKM);
		\node at (z) {\textbullet};
		\node at (x) {\textbullet};
		\node[blue] at (Tx) {\textbullet};
		\node[orange] at (w) {\textbullet};
		\node[orange] at (Tw) {\textbullet};
		\node[black] at (IS) {\textbullet};
		\node[black] at (GKM) {\textbullet};
		
		\node[textnode] at (z) [anchor=east] {\(z\)};
		\node[textnode] at (x) [anchor=west] {\(x\)};
		\node[textnode, blue] at (Tx) [anchor=north] {\(T\!x\)};
		\node[textnode, orange] at (w) [anchor=west] {\(w\)};
		\node[textnode, orange] at (Tw) [anchor=east] {\(T\!w\)};
		
		\node[textnode, black] at (IS) [anchor=north] {\(x^+\)};
		\node[textnode, black] at (GKM) [anchor=north] {\(x_{_{\rm GKM}}^+\)};
	\end{axis}%
\end{tikzpicture}%
%
	}}%
		\caption{}%
		\label{fig:InexactNewton:x+}%
	\end{subfigure}%
	\caption[Comparison between GKM and \cite{solodov1998globally}]{%
		The positive definiteness of \(G\) prevents the update directions \(d\) in the scheme of \cite{solodov1998globally} to point in the gray-shaded area.
		As a result, differently from the GKM scheme the cited algorithm is not robust to any choice of direction (\eg, it cannot accept the one as in \Cref{fig:SuperMann}).
		In any case, the half-space \(C_w\) onto which \(x\) is projected according to the GKM scheme is properly contained in the half-space \(H_w\) corresponding to the update of \cite{solodov1998globally}; consequently, the GKM update is always closer to any solution.%
	}%
	\label{fig:InexactNewton}%
\end{figure*}

Letting \(C_w\) be the half-space as in \cref{lem:DriftedKM}, so that \(x_{_{\rm GKM}}^+=\proj_{C_w}x\) (for simplicity we set \(\lambda=1\)), for the half-spaces \eqref{eq:InexactNewton_x} it holds that
\[
	\zer R\subseteq C_w\subseteq H_w,
\]
the last inclusion holding as equality iff \(Rw=0\).
This means that in the GKM scheme, the same \(w\) yields an iterate \(x_{_{\rm GKM}}^+\) which is closer to any \(z\in\zer R\) with respect to \(x^+\) (cf. \cref{fig:InexactNewton}).
Notice further that the hyperplanes delimiting the two half-spaces are parallel, with \(\boundary C_w\) passing by \(Tw\) (or \(T_{\nicefrac{1}{2\alpha}}w\) for generic \(\alpha\)'s) and \(\boundary H_w\) by \(w\).

The requirement of positive definiteness of matrix \(G\) in defining the update direction \(d\) is due to the fact that \cite{solodov1998globally} addresses a broader class of monotone operators; we instead exploited at full the nonexpansiveness of \(\id-R\) and as a result have complete freedom in selecting \(d\) (\cref{fig:InexactNewton:d}) and better projections (\cref{fig:InexactNewton:x+}).

Moreover, it can be easily verified that the proposed extra-gradient step does not extend the classical KM iteration unless \(T\) has a very peculiar property, namely that
\(
	Rx
{}={}
	\frac{\innprod{RTx}{Rx}}{\|Rx\|^2}RTx
\)
for every \(x\).
(In particular, for such a \(T\) necessarily \(\|Rx\|=\|RTx\|\) for all \(x\), and consequently there cannot exist \(\alpha\in(0,1)\) for which \(T\) is \avg{}).

			\subsubsection{Line-search for KM}
				The recent work \cite{giselsson2016line} proposes an acceleration of the classical KM scheme for finding a fixed point of an \avg{} operator \(T\) based on a line search on the relaxation parameter.
Namely, instead of the \emph{nominal} update
\(
	\bar x=T_\lambda x
\)
with \(\lambda\in[0,\nicefrac1\alpha]\) as in \eqref{eq:KM}, values \(\lambda'>\nicefrac1\alpha\) are first tested and the update
\(
	x^+=T_{\lambda'} x
\)
is accepted as long as
\(
	\|Rx^+\|\leq c_1\|R\bar x\|
\)
holds for some constant \(c_1\in(0,1)\).

In the setting of the \refSM[], this corresponds to selecting \(d_k=-Rx_k\), discarding \emph{blind} updates (\ie setting \(c_0=0\)), foretracking \emph{educated} updates and using plain KM iterations as safeguard steps.
Convergence can be enhanced and the method is attractive when \(T=S_2\circ S_1\) is the composition of an affine mapping \(S_1\) and a cheap operator \(S_2\), in which case the line search is inexpensive.
However, though preserving the same theoretical convergence guarantees of KM (hence of the \refSM[]), it does not improve its best-case local linear rate.

Although other choices \(d_k\) may also be considered, however fast directions such as Newton-type ones would be discarded and replaced by nominal KM updates every time the candidate point \(x_k+d_k\) does not meet some requirements.
Avoiding this take-it-or-leave-it behavior is exactly the primary goal of GKM iterations, so that candidate good directions are never discarded.


			\subsubsection{Smooth optimization with envelope functions}
				For solving nonsmooth minimization problems in composite form, \cite{patrinos2013proximal,patrinos2014douglas} introduced \emph{forward-backward envelope} (FBE) and \emph{Douglas-Rachford envelope} (DRE) functions.
The original nonsmooth problem is recast into the minimization of continuous (possibly continuously differentiable) real-valued exact penalty functions closely related to FBS and DRS, named \emph{envelopes} due to their kinship with the Moreau \emph{envelope} and the proximal point algorithm.
This paved the way for the employment of fast methods for smooth unconstrained minimization problems \cite{patrinos2013proximal,patrinos2014douglas,stella2017forward}, or for globalizing convergence of fast methods for solving nonlinear equations \cite{themelis2016forward,stella2017simple}.
Although they have the advantage of being suited for nonconvex problems, however their employment is limited to composite operators as described above and they cannot handle, for instance, saddle-point convex-concave optimization problems typically arising from primal-dual splittings such as V\~u-Condat \cite{condat2013primal}.
The \refSM[] instead offers a unifying framework that is based uniquely on evaluations of the nonexpansive mapping \(T\), regardless of their structure.

\ifieee
	\section{Simulations -- Linear optimal control}
		\label{sec:Simulations}
\else
	\section{Simulations}
		\label{sec:Simulations}
		We conclude with some numerical examples to give tangible evidence of the robustifying and enhancing effect that the \refSM[] has on fixed-point iterations.
In all simulations we deactivated blind updates by setting \(c_0=0\), and we selected \(\sigma=10^{-3}\) for safeguard updates and \(c_1=q=1-\sigma\) for educated updates.
Due to problem size we used restarted Broyden's directions with Powell's modification \eqref{eq:PowellTheta} with a memory buffer of \(20\) vectors.

		\subsection{Cone programs}
			\label{sec:SCS}
			We consider cone problems of the form
\begin{equation}\label{eq:CP}
	\minimize_{x\in\R^n}{%
		\innprod cx
	}
\qquad\text{s.t.}\quad
	Ax+s=b,~
	s\in\mathcal K
\end{equation}
where \(\mathcal K\) is a nonempty closed convex cone.
Almost any convex program can be recast as \eqref{eq:CP}, and many convex optimization solvers address problems by first translating them into this form.
The KKT conditions for optimality of the primal-dual couple
\(
	\bigl(
		(x_\star,s_\star),
		(y_\star,r_\star)
	\bigr)
\)
are
\ifieee
	\begin{multline*}
		Ax_\star+s_\star=b,~~
		s_\star\in \mathcal K,~~
		\trans A y_\star+c=r_\star,
	\\
		r_\star=0,~~
		y_\star\in\conj{\mathcal K},~~
		\trans cx_\star+\trans by_\star=0
	\end{multline*}
\else
	\[
		Ax_\star+s_\star=b,~~
		s_\star\in K,~~
		\trans A y_\star+c=r_\star,~~
		r_\star=0,~~
		y_\star\in K^\star,~~
		\trans cx_\star+\trans by_\star=0
	\]
\fi
where \(\conj{\mathcal K}\) is the dual cone of \(\mathcal K\).
A recently developed conic solver for \eqref{eq:CP} is SCS \cite{odonoghue2016conic}, which solves the corresponding so-called \emph{homogeneous self-dual embedding}
\begin{equation}\label{eq:SCS-embedding}
	\text{find}
\quad
	u\in\mathcal C
\qquad
	\text{s.t.}
\quad
	Qu\in\conj{\mathcal C}
\end{equation}
where
\[
	\mathcal C
{}={}
	\R^n\times\conj{\mathcal K}\times\R_+
\quad\text{and}\quad
	Q
{}={}
	\setlength\arraycolsep{2pt}
	\begin{bmatrix}
		0		& A^\top	& c \\
		-A		& 0			& b \\
		-c^\top	& -b^\top	& 0
	\end{bmatrix}.
\]
Problem \eqref{eq:SCS-embedding} can be equivalently reformulated as the variational inequality
\begin{equation}\label{eq:SCS-VI}
	\text{find}~
	u\in\mathcal C
\quad
	\text{s.t.}
\quad
	0\in Qu+N_{\mathcal C}(u).
\end{equation}
Indeed, for all \(u\in\mathcal C\) we have
\ifieee
	\begin{align*}
		N_{\mathcal C}(u)
	{}={} &
		\set{y}[\innprod{v-u}{y}\leq0~~\forall v\in\mathcal C]
	\\
	{}={} &
		\set{u}^\bot
		{}\cap{}
		\set{y}[\innprod{v}{y}\leq0~~\forall v\in\mathcal C]
	{}={}
		\set{u}^\bot
		{}\cap{}
		(-\conj{\mathcal C})
	\end{align*}
\else
	\begin{align*}
		N_{\mathcal C}(u)
	{}={} &
		\set{y}[\innprod{v-u}{y}\leq0~~\forall v\in\mathcal C]
	{}={}
		\set{u}^\bot
		{}\cap{}
		\set{y}[\innprod{v}{y}\leq0~~\forall v\in\mathcal C]
	\\
	{}={} &
		\set{u}^\bot
		{}\cap{}
		(-\conj{\mathcal C})
	\end{align*}
\fi
where the second equality follows by considering, \eg, \(v=\frac 12u\) and \(v=\frac32u\), which both belong to \(\mathcal C\) being it a cone.
From this equivalence and the fact that \(Qu\in\set{u}^\bot\) for any \(u\) due to the skew symmetry of \(Q\), the equivalence of \eqref{eq:SCS-embedding} and \eqref{eq:SCS-VI} is apparent.
This leads to the short and elegant interpretation of SCS as Douglas-Rachford splitting (DRS) applied to the splitting \(N_{\mathcal C}+Q\) in \eqref{eq:SCS-VI}, which, after a well known change of variables and index shifting, reads
\begin{equation}
	\begin{cases}
		\tilde u_{k+1} & {}\approx(I+Q)^{-1}(u_k+v_k)\\
		u_{k+1} & {}=\Pi_{\mathcal C}(\tilde u_{k+1}-v_k)\\
		v_{k+1} & {}=v_k-\tilde u_{k+1}+u_{k+1}.
	\end{cases}
\end{equation}
The ``\(\approx\)'' symbol refers to the fact that \(v_k\) may be retrieved inexactly by means of conjugate gradient method (CG); see \cite{odonoghue2016conic} for a detailed discussion.

Here we consider instead DRS applied to the (equivalent) splitting \(Q+N_{\mathcal C}\) in \eqref{eq:SCS-VI}, namely
\begin{equation}\label{eq:DR}
	\begin{cases}
		v_{k+1} & {}\approx(I+Q)^{-1}\fillwidthof[l]{(u_k+v_k)}{(u_k)}\\
		w_{k+1} & {}=\Pi_{\mathcal C}(2v_{k+1}-u_k)\\
		u_{k+1} & {}=u_k+w_{k+1}-v_{k+1}.
	\end{cases}
\end{equation}
For any initial point \(u_0\), the variable \(v_k\) converges to a solution to \eqref{eq:SCS-VI} \cite[Thm. 25.6(i),(iv)]{bauschke2011convex}.
DRS is a (firmly) nonexpansive operator and as such it can be integrated in the \refSM[] with \(\lambda\in(0,2)\); in these simulations we set \(\lambda=1\).

\renewcommand\floatpagefraction{0.1}
\begin{figure}[t]
		\caption[SCS vs SuperSCS]{%
		Comparison between Splitting Cone Solver \cite{odonoghue2016conic} (blue) and its enhancement with the \refSM[] for solving a cone program \eqref{eq:CP}.
	}%
	\ifieee\else
		\addtocounter{figure}{-1}%
	\fi
	\label{fig:SCS}%
	\begin{subfigure}{\linewidth}
		{{%
		\pgfkeys{/pgf/images/include external/.code={\includegraphics[width=\linewidth]{#width=\linewidth}}}%
		\tikzsetnextfilename{Simulation_SCS_m487_n325_condA1e2}%
		\input{./TeX/Tikz/Simulation_SCS_m487_n325_condA1e2.tex}%
	}}
		\caption[SCS vs SuperSCS for a well-conditioned problem]{%
			On the \(x\)-axis the number of times a linear system is solved, the most expensive operation, needed for computing the resolvent of \(Q\).
			SCS performs quite well, however its super-enhancement converges considerably faster in terms of operations.
		}%
		\label{fig:SCSexact}%
	\end{subfigure}
	\begin{subfigure}{\linewidth}
		{{%
		\pgfkeys{/pgf/images/include external/.code={\includegraphics[width=\linewidth]{#width=\linewidth}}}%
		\tikzsetnextfilename{Simulation_SCS_m487_n325_condA1e2_CG}%
		\input{./TeX/Tikz/Simulation_SCS_m487_n325_condA1e2_CG.tex}%
	}}
		\caption[SCS vs SuperSCS for a well-conditioned problem (solved with CG)]{%
			Comparison with respect to the same problem, but with linear systems solved approximately with CG on a reduced system.
			On the \(x\)-axis the number of times either the operator \(A\) or \(\trans A\) is called, which amounts to the most expensive operations.
			The comparison between SCS and super-SCS is quite identical to that in which the linear system is solved exactly.%
		}%
		\label{fig:SCSinexact}%
	\end{subfigure}
\end{figure}
We run a cone problem \eqref{eq:CP} of size \(m=487\) and \(n=325\), with \(\operatorname{dens}(A)=0.01\) and \(\cond(A)=100\), both by solving exactly the linear systems and by adopting the CG technique.
\(\mathcal C\) is the cartesian product of all the primitive cones implemented in SCS solver: positive orthant, second-order, positive semidefinite, (dual) exponential, and (dual) power cones.
We reported primal residual, dual residual, and duality gap; consistently with SCS' termination criterion, the algorithm is stopped when all these quantities are below some tolerance \cite[§3.5]{odonoghue2016conic}, which we set to \(10^{-6}\).

Notice that if \(u\) solves \eqref{eq:SCS-embedding}, or equivalently \eqref{eq:SCS-VI}, then so does any multiple \(tu\) with \(t>0\).
In particular no isolated solution exists, and therefore whenever the residual \(R\) of the DRS operator is differentiable at a solution \(u_\star\), \(\jac R{u_\star}\) is singular.
Fortunately, the \refSM[] does not necessitate nonsingularity of the Jacobian but merely metric subregularity, the same property that enables linear convergence rate of the original DRS (or equivalently SCS).
In particular, whenever the original SCS scheme is linearly convergent, the SuperMann enhancement is provably superlinear provided that \(R\) is strictly differentiable at the limit point.
However, since \emph{restarted} Broyden directions are implemented instead of the full-memory method, rather than superlinear convergence we expect an ``extremely steep'' linear convergence.

In \Cref{fig:SCS} we can observe how the original SCS scheme (blue) converges at a fair linear rate; however, its super-enhancement greatly outperforms it both when solving linear systems exactly and approximately.

		\subsection{Lasso}
			We consider a lasso problem
\[
	\minimize_{x\in\R^n}{%
		\tfrac12\|Ax-b\|^2
		{}+{}
		\nu\|x\|_1
	}
\]
where \(A\in\R^{m\times n}\), \(b\in\R^m\) and \(\nu>0\).
In \Cref{fig:Lasso} the comparison of forward-backward splitting (or proximal gradient, in blue) and its super-enhanced version (red) on a random problem with \(m=1500\) \(n=5000\) and \(\nu=10^{-2}\).
On the \(x\)-axis the number of matvecs, being them the most expensive operations of FB and hence of super-FB, and on the \(y\)-axis the fixed-point residual.
Though superlinear convergence cannot be observed due to the fact that a limited-memory method is used for computing directions, however an outstanding speedup is noticeable.
\begin{figure}[t]
	\begin{minipage}[t]{0.33\linewidth}\vspace*{0pt}
		{{%
		\pgfkeys{/pgf/images/include external/.code={\includegraphics[width=\linewidth]{#width=\linewidth}}}%
		\tikzsetnextfilename{Simulation_Lasso}%
		\input{./TeX/Tikz/Simulation_Lasso.tex}%
	}}
	\end{minipage}
	\begin{minipage}[t]{0.66\linewidth}\vspace*{0pt}
		\caption[Lasso simulation]{Comparison between FBS and Super-FBS (using modified Broyden limited-memory directions) in a lasso problem.}
		\label{fig:Lasso}
	\end{minipage}
\end{figure}

		\subsection{Constrained linear optimal control}
\fi
\added{%
		For matrices \(A_t\) and \(B_t\) of suitable size, \(t=0,\ldots,N-1\), consider a state-input dynamical system
\begin{subequations}\label{eq:ControlP}
	\begin{equation}\label{eq:Dynamics}
		x_{t+1}
	{}={}
		A_tx_t + B_tu_t,
	\quad
		t=0,\ldots,N-1,
	\end{equation}
	where the \(x_0\in\R^{n_x}\) is given, and the next states \(x_t\in\R^{n_x}\) are determined by the user-defined inputs \(u_\tau\in\R^{n_u}\), \(\tau=0,\ldots,t-1\).
	States \(\bm x=(x_1,\ldots,x_N)\) can be expressed in terms of the inputs \(\bm u=(u_0,\ldots,u_{N-1})\) through a linear operator \(L\in\R^{Nn_x\times Nn_u}\) as \(\bm x=L\bm u+b\) for some constant \(b\in\R^{Nn_x}\).
	The goal is to choose inputs that minimize a cost
	\begin{equation}\label{eq:Cost}
		\ell(\bm u,\bm x)
	{}={}
		\sum_{t=0}^{N-1}{
			\ell_t(u_t,x_t)
		}
		{}+{}
		\ell_N(x_N)
	\end{equation}
	subject to some constraints
	\begin{equation}\label{eq:Constraints}
		x_{t+1}\in\mathcal X_{t+1},
	\quad
		u_t\in\mathcal U_t,
	\quad
		t=0,\ldots,N-1.
	\end{equation}
\end{subequations}

\ifieee
	\subsection{V\~u-Condat splitting}
\else
	\subsubsection{V\~u-Condat splitting}
\fi
The constraint sets in \eqref{eq:Constraints} are typically simple and easy to project onto (boxes, Euclidean balls\ldots).
However, while simple input constraints can be easily handled, due to the coupling enforced by the dynamics \eqref{eq:Dynamics}, expressing \(\mathcal X_{t+1}\) in terms of the optimization variable \(\bm u\) results in much more complicated sets (polyhedra, ellipsoids\ldots).
To avoid this complication we make use of the extremely versatile algorithm that V\~u-Condat three-term splitting offers \cite[Alg. 3.1]{condat2013primal}.
In its general form, the algorithm addresses problems of the form
\begin{equation}\label{eq:ProblemVC}
	\minimize_{x\in\R^n}{
		f(x)+g(x)+h(Lx)
	}
\end{equation}
where \(\func{f}{\R^n}{\R}\) is convex with \(L_f\)-Lipschitz continuous gradient, \(\func{g}{\R^n}{\Rinf}\) and \(\func{h}{\R^m}{\Rinf}\) are convex, and \(L\in\R^{n\times m}\), by iterating the following steps:
\begin{equation}\label{eq:VC}
	\begin{cases}
			x^+
		& {}={}
			\prox_{\fillwidthof[l]{\tau\conj h}{\tau g}}\!\!
			\left(
				x-\tau(\nabla f(x)+\trans Ly)
			\right)
		\\
			y^+
		& {}={}
			\prox_{\tau\conj h}\!\!
			\left(
				y+\sigma L(2x^+-x)\vphantom{\trans L}
			\right).
	\end{cases}
\end{equation}
Here,
\(
	0
{}<{}
	\tau
{}<{}
	\frac{2}{L_f}
\)
and
\(
	0
{}<{}
	\sigma
{}<{}
	\frac{1}{\|L\|^2}
	\bigl(
		\frac1\tau-\frac{L_f}{2}
	\bigr)
\)
are stepsizes, and \(y\in\R^m\) is a Lagrange multiplier.
V\~u-Condat splitting is a primal-dual method that generalizes FBS by allowing an extra nonsmooth term \(h\) and a linear operator \(L\) (by neglecting \(h\) and \(L\) one recovers the proximal gradient iterations of FBS).

The optimal control problem \eqref{eq:ControlP} can be cast into V\~u-Condat splitting form \eqref{eq:ProblemVC} by simply letting
\(
	f(\bm u)
{}={}
	\ell(\bm u,L\bm u)
\),
\(
	g=\indicator_{\mathcal U}
\)
and
\(
	h=\indicator_{\mathcal X}({}\cdot{}+b)
\),\footnote{%
	\(\indicator_C\) denotes the \emph{indicator function} of the nonempty closed convex set \(C\), namely \(\indicator_C(x)=0\) if \(x\in C\) and \(\indicator_C(x)=\infty\) otherwise.%
}
where
\(\mathcal U=\mathcal U_0\times\cdots\times\mathcal U_{N-1}\) and \(\mathcal X=\mathcal X_1\times\cdots\times\mathcal X_N\) (in particular, \(n=Nn_u\) and \(m=Nn_x\)).
Then, \(\prox_{\tau g}=\proj_{\mathcal U}\) and \(\prox_{\sigma\conj h}(y)=y-\sigma\proj_{\mathcal X}(\sigma^{-1}y+b)+b\).
Notice that \(\proj_{\mathcal U}\) and \(\proj_{\mathcal X}\) are fully decoupled as the projection of each input and state onto the corresponding constraint set.
Moreover, the full matrix \(L\) needs not be computed, as both \(L\) and \(\trans L\) can be treated as abstract operators that simulate forward and backward dynamics.

Apparently, the appeal of V\~u-Condat splitting in addressing the optimal control problem lies in the extreme simplicity of its operations and low memory requirements, making it particularly suited for medium-to-large-scale problems in which traditional interior point algorithms fail.
However, like all first-order methods it is extremely sensitive to ill conditioning, which gets worse as the problem size increases.
Fortunately, this splitting fits into the \textit{SuperMann} framework.
The operator \(T\) that maps \((x,y)\) into \((x^+,y^+)\) as in \eqref{eq:VC} is averaged in the Hilbert space \(\HH_P\), where \(\HH_P\) is defined as \(\R^n\times\R^m\) equipped with the scalar product \(\innprod{z}{z'}_P\coloneqq\innprod{z}{Pz'}\), where \(P\coloneqq{}\){\scriptsize\arraycolsep=0.3\arraycolsep\(
	\begin{pmatrix}
		\tau^{-1} I & -\trans L\\
		-L & \sigma^{-1} I
	\end{pmatrix}
\)}
\cite[proof of Thm. 3.1]{condat2013primal}.


\ifieee
	\subsection{Oscillating masses experiment}
\else
	\subsubsection{Oscillating masses experiment}
\fi
We tried this approach on the benchmark problem of controlling a chain of oscillating masses connected by springs and with both ends attached to walls.
The chain is composed of \(2K\) bodies of unit mass subject to a viscous friction of \(0.1\), the springs have elastic constant \(1\) and no damping, and the system is controlled through \(K\) actuators, each being a force acting on a pair of masses, as depicted in \Cref{fig:chain}.
Therefore \(n_x=4K\) (the states are the displacement from the rest position and velocity of each mass) and \(n_u=K\).
The inputs are constrained in \([-2,2]\), while the position and velocity of each mass is constrained in \([-5,5]\).

The continuous-time system was discretized with a sampling time \(T_s = 0.1s\).
We considered quadratic stage costs \(\frac12\trans xQx\) for the states and \(\tfrac12\trans uu\) for the inputs, where \(Q\) is diagonal positive definite with random diagonal entries, and generated a random (feasible) initial state \(x_0\).
Notice that a QP reformulation would require the computation of the full cost matrix, differently from the splitting approach where only the small dynamics matrices \(A\) and \(B\) are needed, as \(L\) and \(\trans L\) can be abstract operators.

\begin{figure}
	{{%
		\pgfkeys{/pgf/images/include external/.code={\includegraphics[width=\linewidth]{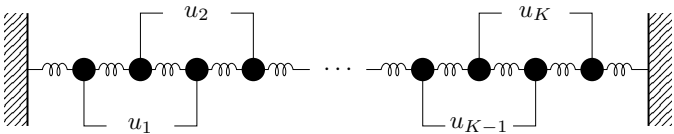}}}%
		\tikzsetnextfilename{chain}%
		\begin{tikzpicture}[every node/.style={outer sep=0pt},scale=0.8]
\tikzstyle{spring}=[decorate,decoration={coil, pre length=6, post length=6, segment length=3}]
\fill[pattern = north east lines] (0,-1) rectangle (-0.4,1);
\draw[thick] (0,-1) -- (0,1);

\def \springlength {1.0}
\def \handlespace {0.6}

\foreach \k in {1}
{

    \def \basex {((\k - 1)*4)*\springlength}
    \draw[spring] ({\basex},0) -- ({\basex+\springlength},0);
    \node[circle,fill=black] (a) at ({\basex+\springlength},0) {};
    \draw[spring] ({\basex+\springlength},0) -- ({\basex+2*\springlength},0);
    \node[circle,fill=black] (a) at ({\basex+2*\springlength},0) {};
    \draw[spring] ({\basex+2*\springlength},0) -- ({\basex+3*\springlength},0);
    \node[circle,fill=black] (a) at ({\basex+3*\springlength},0) {};
    \draw[spring] ({\basex+3*\springlength},0) -- ({\basex+4*\springlength},0);
    \node[circle,fill=black] (a) at ({\basex+4*\springlength},0) {};


    \pgfmathsetmacro\ucount{(\k-1)*2+1}

    \draw ({\basex+\springlength},0) -- ({\basex+\springlength},-1) -- ({\basex+2*\springlength-\handlespace},-1);
    \draw ({\basex+3*\springlength},0) -- ({\basex+3*\springlength},-1) -- ({\basex+2*\springlength+\handlespace},-1);
    \node at ({\basex+2*\springlength},-1) {\(\:u_1\:\)}; 

    \pgfmathsetmacro\ucount{(\k-1)*2+2}

    \draw ({\basex+2*\springlength},0) -- ({\basex+2*\springlength},+1) -- ({\basex+3*\springlength-\handlespace},+1);
    \draw ({\basex+4*\springlength},0) -- ({\basex+4*\springlength},+1) -- ({\basex+3*\springlength+\handlespace},+1);
    \node at ({\basex+3*\springlength},+1) {\(\:u_2\:\)}; 
}

\draw[spring] ({4*\springlength},0) -- ({5*\springlength},0);
\node at ({5.5*\springlength},0) {\(\cdots\)};

\foreach \k in {1}
{

    \def \basex {(\k*4+2)*\springlength}
    \draw[spring] ({\basex},0) -- ({\basex+\springlength},0);
    \node[circle,fill=black] (a) at ({\basex+\springlength},0) {};
    \draw[spring] ({\basex+\springlength},0) -- ({\basex+2*\springlength},0);
    \node[circle,fill=black] (a) at ({\basex+2*\springlength},0) {};
    \draw[spring] ({\basex+2*\springlength},0) -- ({\basex+3*\springlength},0);
    \node[circle,fill=black] (a) at ({\basex+3*\springlength},0) {};
    \draw[spring] ({\basex+3*\springlength},0) -- ({\basex+4*\springlength},0);
    \node[circle,fill=black] (a) at ({\basex+4*\springlength},0) {};


    \pgfmathsetmacro\ucount{(\k-1)*2+1}

    \draw ({\basex+\springlength},0) -- ({\basex+\springlength},-1) -- ({\basex+2*\springlength-\handlespace},-1);
    \draw ({\basex+3*\springlength},0) -- ({\basex+3*\springlength},-1) -- ({\basex+2*\springlength+\handlespace},-1);
    \node at ({\basex+2*\springlength},-1) {\(\:u_{K-1}\:\)}; 

    \pgfmathsetmacro\ucount{(\k-1)*2+2}

    \draw ({\basex+2*\springlength},0) -- ({\basex+2*\springlength},+1) -- ({\basex+3*\springlength-\handlespace},+1);
    \draw ({\basex+4*\springlength},0) -- ({\basex+4*\springlength},+1) -- ({\basex+3*\springlength+\handlespace},+1);
    \node at ({\basex+3*\springlength},+1) {\(\:u_{K}\:\)}; 
}

\def \basex {(2*4+2)*\springlength};

\draw[spring] ({\basex}, 0) -- ({\basex+\springlength}, 0);
\fill[pattern = north east lines] ({\basex+\springlength},-1) rectangle ({\basex+\springlength+0.4},1);
\draw[thick] ({\basex+\springlength},-1) -- ({\basex+\springlength},1);

%
%
%
%
%
%
%
%
%
\end{tikzpicture}%
	}}%
	\caption{Oscillating masses}%
	\label{fig:chain}%
\end{figure}
We simulated different scenarios for all combinations of \(K\in\set{8,16}\) and \(N\in\set{10,20,30,40,50}\).
We compared Vu-Condat splitting (VC) with its `super' enhancement (SuperVC); parameters were set as detailed in \Cref{sec:Parameters}.
\Cref{fig:springmass_K8_N30c} shows a comparison of the convergence rates for one problem instance, while \Cref{table:MPC} offers an overview of the whole experiment: SuperVC is roughly 13 times faster on average and 21 times better in worst-case performance than VC algorithm in reaching the termination criterion \(\|Rx^k\|\leq10^{-4}\|Rx^0\|\).

\begin{table}
	\centering
	\noindent
	\setlength\tabcolsep{2pt}%
\footnotesize
\begin{tabularx}{\linewidth}{@{}r @{\,} | rr | rr | rr | rr | rr |@{}}
	\multicolumn{1}{c}{}%
	&
	\multicolumn{10}{c}{\textbf{Number of calls to \(L\) and \(\trans L\) (\(\times 10^3\))}}%
\\\cline{2-11}%
	\multicolumn{1}{@{}c|}{\(K=8\)\vphantom{\(X^{X^X}\)}}%
	&
	\multicolumn{2}{>{\centering \arraybackslash}X|}{\(N=10\)}&%
	\multicolumn{2}{>{\centering \arraybackslash}X|}{\(N=20\)}&%
	\multicolumn{2}{>{\centering \arraybackslash}X|}{\(N=30\)}&%
	\multicolumn{2}{>{\centering \arraybackslash}X|}{\(N=40\)}&%
	\multicolumn{2}{>{\centering \arraybackslash}X|@{}}{\(N=50\)}%
\\
	& \multicolumn{1}{c}{avg} & \multicolumn{1}{c|}{max}
	& \multicolumn{1}{c}{avg} & \multicolumn{1}{c|}{max}
	& \multicolumn{1}{c}{avg} & \multicolumn{1}{c|}{max}
	& \multicolumn{1}{c}{avg} & \multicolumn{1}{c|}{max}
	& \multicolumn{1}{c}{avg} & \multicolumn{1}{c|}{max}
\\\cline{2-11}\vphantom{\(X^{X^X}\)}%
	VC      &   19.0 &  337.1       &   15.0 &  174.4       &   25.0 &   400+       &   21.0 &  136.5       &   16.0 &   61.9 \\
	SVC     &    1.0 &    5.5       &    1.0 &    4.3       &    2.0 &   19.3       &    2.0 &   10.9       &    2.0 &    6.6 \\
\cline{2-11}
	\multicolumn{11}{@{}c@{}}{}
	\\[-5pt]
\cline{2-11}%
	\multicolumn{1}{@{}c|}{\(K=16\)\vphantom{\(X^{X^X}\)}}%
	&
	\multicolumn{2}{c|}{\(N=10\)}&%
	\multicolumn{2}{c|}{\(N=20\)}&%
	\multicolumn{2}{c|}{\(N=30\)}&%
	\multicolumn{2}{c|}{\(N=40\)}&%
	\multicolumn{2}{c|@{}}{\(N=50\)}%
\\
	& \multicolumn{1}{c}{avg} & \multicolumn{1}{c|}{max}
	& \multicolumn{1}{c}{avg} & \multicolumn{1}{c|}{max}
	& \multicolumn{1}{c}{avg} & \multicolumn{1}{c|}{max}
	& \multicolumn{1}{c}{avg} & \multicolumn{1}{c|}{max}
	& \multicolumn{1}{c}{avg} & \multicolumn{1}{c|}{max}
\\\cline{2-11}\vphantom{\(X^{X^X}\)}%
	VC      &   62.0 &   400+       &   30.0 &  344.9       &   30.0 &   400+       &   65.0 &   400+       &   29.0 &  318.6 \\
	SVC     &    4.0 &   39.5       &    2.0 &   11.6       &    3.0 &   46.6       &    8.0 &   58.1       &    3.0 &   26.1 \\
\cline{2-11}
\end{tabularx}
	\caption{%
		Comparison between V\~u-Condat algorithm (VC) and its ``super'' enhancement (SuperVC) in solving the oscillating masses problem with \(\|Rx^k\|\leq 10^{-4}\|Rx^0\|\) as termination criterion.
		Average and worst performances among 25 simulations with randomly generated starting point \(x_0\) for each combination of \(K\in\set{8,16}\) and \(N\in\set{10,20,30,40,50}\).
		The tables compare the number of calls to the operators \(L\) and \(\trans L\), which are the expensive operations (the rest are projections on boxes).
		In four problems V\~u-Condat exceeded \(4\cdot 10^5\) many calls (corresponding to \(10^5\) iterations) and was stopped prematurely.%
	}%
	\label{table:MPC}%
\end{table}
\begin{figure}
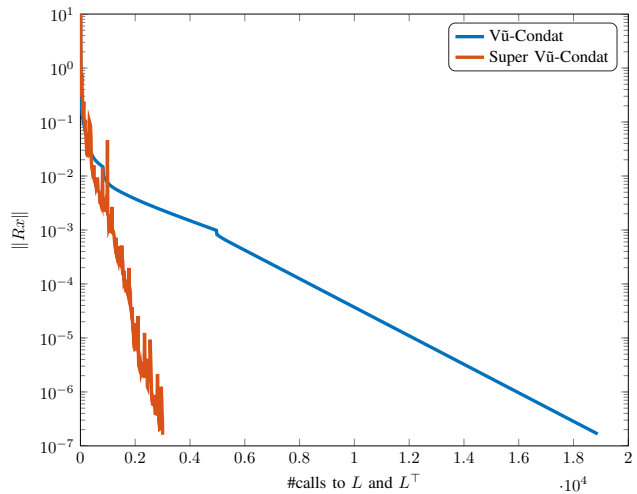

	\centering
	{{%
		\pgfkeys{/pgf/images/include external/.code={\includegraphics[width=.949\linewidth]{#width=.949\linewidth}}}%
		\tikzsetnextfilename{springmass_K8_N30c}%
		\input{./TeX/Tikz/springmass_K8_N30c.tex}%
	}}
	\caption{%
		Random simulation of the spring-mass control problem with \(K=8\) and \(N=30\).%
	}%
	\label{fig:springmass_K8_N30c}%
\end{figure}

}%

	\section{Conclusions}
		\label{sec:Conclusions}
		We proposed the \refSM, a novel algorithm for finding fixed points of a nonexpansive operator \(T\) that generalizes and greatly improves the classical \KM{} (KM) scheme, enjoying the same favorable properties:
global convergence with worst-case sublinear rate,
cheap iterations based solely on evaluations of \(T\), and
easy codability.
The \refSM[] is an extremely versatile algorithm, its flexibility being twofold: on one hand it works for any nonexpansive operator \(T\) by requiring only the oracle \(x\mapsto Tx\); on the other hand it allows for the integration of any fast local method for solving nonlinear equations, leaving much freedom for trading-off cheap iterations or faster convergence.
The remarkable performance of the method is supported both in practice with promising simulations and in theory where the employment of quasi-Newton directions is shown to yield asymptotic superlinear convergence rates provided a condition analogous to the famous result by Dennis and Moré is satisfied.
Most importantly, superlinear convergence does not require nonsingularity of the Jacobian of the residual at the solution but merely metric subregularity, and as such can be achieved even when the solution is not isolated.

We encourage the employment of the \refSM[] to improve and robustify convex splitting algorithms; in particular, we strongly believe that its integration in generic solvers which are based on fixed-point iterations of nonexpansive operators such as SCS \cite{odonoghue2016conic} would be extremely beneficial.


	\ifams
		\bibliographystyle{plain}
	\fi
	\ifieee
		\bibliographystyle{IEEEtran}
		\let\appendix\appendices
		\let\endappendix\endappendices
	\fi

\bibliography{TeX/Bibliography.bib}


	\begin{appendix}
		\ifieee
			\crefalias{section}{appsec}%
		\fi
		\renewcommand\appendixproofnameref[1]{}%
		\proofsection{sec:General}
			\begin{appendixproof}{thm:General:Global}
\begin{proofitemize}
	\item\ref{thm:General:QF}:
		we start observing that because of \eqref{eq:D} and the triangular inequality, for all \(k\in K_0\cup K_1\) we have
		\begin{align}
		\label{eq:K+1z}
			\|x_{k+1}-z\|
		{}\leq{} &
			\|x_k-z\|
			{}+{}
			D\|Rx_k\|
		\qquad
			\forall z\in\fix T
		\shortintertext{and since \(R\) is \(2\alpha\)-Lipschitz continuous we also have that}
		\label{eq:Rk+1}
			\|Rx_{k+1}\|
		{}\leq{} &
		\ifieee\else
			\|Rx_k\|
			{}+{}
			\|Rx_{k+1}-Rx_k\|
		{}\leq{}
		\fi
			(1+2\alpha D)
			\|Rx_k\|.
		\end{align}
		Combining \cite[Prop. 3.2(i)]{combettes2001quasi} with \eqref{eq:General:K2} and \eqref{eq:K+1z}, it follows that in order to prove quasi-Fejér monotonicity it suffices to show that the sequence \(\seq{\|Rx_k\|}[k\in K_0\cup K_1]\) is summable.
		Let \(K_0\) and \(K_1\) be indexed as in \eqref{eq:K0K1}.
		Since \(\eta_k\) is kept constant whenever \(k\notin K_0\),
		\begin{equation}\label{eq:K0Linear}
		\ifieee
				\eta_{k_\ell}
			{}={}
				\|Rx_{k_{\ell-1}}\|
			{}\leq{}
				c_0\eta_{k_{\ell-1}}
			\cdots\leq{}
				c_0^{\ell-1}\eta_{k_1}
			{}={}
				c_0^{\ell-1}\eta_0.
		\else
				\eta_{k_\ell}
			{}={}
				\|Rx_{k_{\ell-1}}\|
			{}\leq{}
				c_0\eta_{k_{\ell-1}}
			{}\leq\cdots\leq{}
				c_0^{\ell-1}\eta_{k_1}
			{}={}
				c_0^{\ell-1}\eta_0
			\qquad
				\forall k_\ell\in K_0.
		\fi
		\end{equation}
		In particular, \(\seq{\|Rx_{k_\ell}\|}[k_\ell\in K_0]\) is summable (regardless of whether \(K_0\) is finite or not).
		
		As for \(k_\ell'\in K_1\), the safeguard parameter \(r_{\rm safe}\) ensures that
		\ifieee
			\begin{align*}
				\|Rx_{k'_\ell}\|
			{}\leq{} &
				\|Rx_{k'_{\ell-1}+1}\|+q^{k'_{\ell-1}}
				{}\leq{}
				 c_1\|Rx_{k'_{\ell-1}}\|+q^{k'_{\ell-1}}
			\\
			{}\leq{} &
				 c_1\|Rx_{k'_{\ell-1}}\|+q^{\ell-1}
			\end{align*}
		\else
			\[
				\|Rx_{k'_\ell}\|
			{}\leq{}
				\|Rx_{k'_{\ell-1}+1}\|+q^{k'_{\ell-1}}
			{}\leq{}
				 c_1\|Rx_{k'_{\ell-1}}\|+q^{k'_{\ell-1}}
			{}\leq{}
				 c_1\|Rx_{k'_{\ell-1}}\|+q^{\ell-1}
			\quad
				\forall k_\ell'\in K_1.
			\]
		\fi
		Iterating the inequality, for any \(\rho\in(0,1)\) such that \(\rho>\max\set{c_1,q}\) we have
		\begin{equation}\label{eq:K1linear}
			\|Rx_{k'_\ell}\|
		{}\leq{}
			\rho^{\ell-1}
			\|Rx_{k'_1}\|
			{}+{}
			\sum_{i=1}^{\ell-1} c_1^{i-1}\rho^{\ell-i}
		{}\leq{}
			C\rho^\ell
		\end{equation}
		where
		\(
			C
		{}\coloneqq{}
			\tfrac 1\rho
			\left(
				\|Rx_{k'_1}\|
				{}+{}
				\sum_{i\in\N}\left(\nicefrac{c_1}{\rho}\right)^i
			\right)
		{}<{}
			\infty
		\).
		In particular, also \(\seq{\|Rx_k\|}[k\in K_1]\) is summable.

	\item\ref{thm:General:ResStr}:
		due to quasi-Fejér monotonicity, for all \(z\in\fix T\) there exists \(\seq{\varepsilon_k(z)}\in\ell_1^+\) such that
		\[
			\|x_{k+1}-z\|^2
		{}\leq{}
			\|x_k-z\|^2
			{}+{}
			\varepsilon_k(z).
		\]
		Combining this with \eqref{eq:General:K2} and telescoping the inequalities, we obtain that for all \(z\in\fix T\)
		\begin{equation}\label{eq:sigmaConstant}
			\|x_0-z\|^2
		{}\geq{}
			\sigma
			\sum_{k\in K_2}\|Rx_k\|^2
			~{}-{}~~
			\sum_{\mathclap{k\in K_0\cup K_1}}\varepsilon_k(z).
		\end{equation}
		Since the sequence \(\seq{\varepsilon_k(z)}[k\in K_0\cup K_1]\) is summable, then so is \(\seq{\|Rx_k\|^2}[k\in K_2]\).
		In turn, since \(\seq{\|Rx_k\|}[k\in K_0\cup K_1]\) is also summable it follows that the whole sequence of residuals is square-summable.

	\item\ref{thm:General:xWeak}:
		follows by combining \ref{thm:General:ResStr} with \cref{prop:WeakFix}.

	\item\ref{thm:General:K0}:
		trivially follows from the already proven point \ref{thm:General:ResStr}, together with the observation that since \(\eta_k\) is kept constant whenever \(k\notin K_0\), the condition \(\|Rx_k\|\leq c_0\eta_k\) will be satisfied infinitely often if \(c_0>0\).
\qedhere
\end{proofitemize}
\end{appendixproof}

We now state two lemmas which will be needed in the proof of \Cref{thm:General:Linear}.

\begin{lem}[Asymptotic properties of \(K_0\) and \(K_1\)]\label{prop:K0K1}
Suppose the hypotheses of \Cref{thm:General:Global} hold and let \(\seq{x_k}\) be the sequence generated by \Cref{alg:General}.
Then,
\begin{enumerate}
	\item\label{prop:K0linear}
		\(\seq{\|Rx_k\|}[k\in K_0]\) is \(Q\)-linearly convergent;
	\item\label{prop:K1linear}
		\(\seq{\|Rx_k\|}[k\in K_1]\) is \(R\)-linearly convergent;
	\item\label{prop:K1<2K0}
		if \(c_0>0\) then for some \(\varrho\in(0,1]\) and \(\beta\in\R\)
		\[
			\ell_0(k)\geq\varrho\,\ell_1(k)-\beta
		\quad\forall k\in\N,
		\]
		where
		\(
			\ell_j(k)
		{}\coloneqq{}
			\#\set{k'\in K_j}[k'\leq k]
		\),
		\(j=0,1,2\), is the number of times \(K_j\) was visited up to iteration \(k\).
\end{enumerate}
\begin{proof}
\begin{proofitemize}
	\item\ref{prop:K0linear} and \ref{prop:K1linear}:
		already shown in \eqref{eq:K0Linear} and \eqref{eq:K1linear}.
	\item\ref{prop:K1<2K0}:
		if \(c_1=0\), then \(K_1=\emptyset\) and the claim trivially holds with \(\varrho=1\) and \(\beta=0\).
		Otherwise, from \eqref{eq:K1linear} and due to the definition of \(\ell_1(k)\) there exist \(C>0\) and \(\rho\in(0,1)\) such that
		\[
			\|Rx_k\|
		{}\leq{}
			C\rho^{\ell_1(k)}
		\qquad
			\forall k\in K_1.
		\]
		If \(k\in K_1\), then \(\|Rx_k\|\) didn't pass the test at step \ref{step:General:K0}, therefore
		\[
			C\rho^{\ell_1(k)}
		{}\geq{}
			\|Rx_k\|
		{}\geq{}
			\eta_k
		{}\overrel{\eqref{eq:K0Linear}}{}
			\|Rx_0\|c_0^{\ell_0(k)}.
		\]
		The proof now follows by simply taking the logarithm on the outer inequality.
	\qedhere
\end{proofitemize}
\end{proof}
\end{lem}

\begin{lem}\label{prop:Rlinear}
Let \(\seq{u_k}\subset[0,+\infty)\) be a sequence, and let \(K_1,K_2\subseteq\N\) be such that \(\N=K_1\cup K_2\).
Let \(K_1\) be indexed as \(K_1=\set{k_0',k_1'\ldots}\), and suppose that there exist \(a,b>0\) and \(\rho\in(0,1)\) such that
\[
	\begin{cases}[r @{{}\leq{}} l @{\quad} l]
		u_{k+1} & au_k & \text{for all } k\in\N,\\
		u_{k_\ell'} & b\rho^\ell & \text{for all } k_\ell'\in K_1,\\
		u_{k+1} & \rho u_k & \text{for all } k\in K_2.
	\end{cases}
\]
Then, there exists \(\sigma\in(0,1)\) such that \(u_k\leq ab\sigma^k\).
\begin{proof}
To exclude trivialities we assume that \(K_1\) and \(K_2\) are both infinite.
To arrive to a contradiction, for all \(\sigma\in(0,1)\) let \(k=k(\sigma)\) be the minimum such that \(u_k>ab\sigma^k\).
Let \(\sigma\geq\rho\) be fixed.
If \(k-1\in K_2\), then
\[
	\rho u_{k-1}
{}\geq{}
	u_k
{}>{}
	ab\sigma^k
{}\geq{}
	ab\rho\sigma^{k-1}
\]
and therefore \(u_{k-1}>ab\sigma^{k-1}\) which contradicts minimality of \(k\).
It follows that necessarily \(k-1\in K_1\), hence \(k-1=k_\ell'\in K_1\) for some \(\ell\in\N\).
For all \(n\in\N\), let \(k_{\ell_n}'=k(\rho^{\nicefrac 1n})-1\), \ie the minimum such that
\(
	u_{k_{\ell_n}'+1}
{}>{}
	ab\rho^{\frac{k_{\ell_n}'+1}{n}}
\).
Combining with the property of \(K_1\) we obtain
\begin{equation}\label{eq:ell_n}
	ab\rho^{\frac{k_{\ell_n}'+1}{n}}
{}<{}
	u_{k_{\ell_n}'+1}
{}\leq{}
	au_{k_{\ell_n}'}
{}\leq{}
	ab\rho^{\ell_n}
\end{equation}
and in particular
\(
	\ell_n
{}\leq{}
	\frac{k_{\ell_n}'}{n}
\).
This means that up to \(k=k_{\ell_n}'\) there are at most \(\nicefrac kn\) elements in \(K_1\), and consequently at least \(k-\nicefrac kn\) in \(K_2\).
Therefore,
\[
	\mathloose
	b\rho^{\frac{k+1}{n}}
{}\overrel[<]{\eqref{eq:ell_n}}{}
	u_k
{}\leq{}
	a^{\nicefrac kn}\rho^{k-\nicefrac kn}u_0.
\]
Taking the \(k\)-th square root on the outer inequality we are left with
\[
	\left(\nicefrac 1\rho\right)^{1-\nicefrac 2n-\nicefrac{1}{nk}}
{}<{}
	\left(\nicefrac{u_0}{b}\right)^{\nicefrac 1k}
	a^{\nicefrac 1n}.
\]
Letting \(n\to+\infty\), so that also \(k\to+\infty\), we arrive to the contradiction \(\rho\geq 1\).
\end{proof}
\end{lem}

\begin{appendixproof}{thm:General:Linear}

\noindent
Letting \(e_k\coloneqq\dist(x_k,\fix T)\), because of \eqref{eq:Rk+1} and \eqref{eq:distR} there exists \(B>1\) such that
\begin{equation}\label{eq:k+1bounds}
	\|Rx_{k+1}\|\leq B\|Rx_k\|
	\quad\text{and}\quad
	e_{k+1}\leq Be_k
\quad
	\forall k\in\N.
\end{equation}
Suppose that \(R\) is metrically subregular at \(x_\star\) with radius \(\varepsilon>0\) and modulus \(\gamma>0\); since \(x_k\to x_\star\), up to an index shifting without loss of generality we may assume that \(\seq{x_k}\subset\ball{x_\star}{\varepsilon}\).
Let \(z_k=\proj_{\fix T}x_k\), so that \(e_k=\|x_k-z_k\|\); combining \eqref{eq:General:K2} and \eqref{eq:distR} we obtain
\ifieee
	that for all \(k\in K_2\)
\fi
\begin{equation}\label{eq:K2linear}
	e_{k+1}^2
{}\leq{}
	\|x_{k+1}-z_k\|^2
{}\leq{}
	\|x_k-z_k\|^2
	{}-{}
	\sigma\|Rx_k\|^2
{}\leq{}
	\rho^2e_k^2
\ifieee\else
	\quad
	\forall k\in K_2
\fi
\end{equation}
where
\(
	\rho
{}\coloneqq{}
	\sqrt{
		1-\nicefrac{\sigma}{\gamma^2}
	}
{}\in{}
	(0,1)
\).
By possibly enlarging \(\rho\) we assume \(\rho\geq\max\set{c_0, c_1}\).

If \(c_0=0\), then \(K_0=\emptyset\) and using \cref{prop:K1linear} and \eqref{eq:k+1bounds} we may invoke \cref{prop:Rlinear} to infer \(R\)-linear convergence of the sequence \(\seq{e_k}\) and conclude the proof.

Therefore, let us suppose that \(c_0>0\), so that by \cref{thm:General:K0} the set \(K_0\) contains infinite many indices.
We now show that there exists \(n\in\N\) such that every \(n\) consecutive indices at least one is in \(K_0\).
Let \(k\in K_0\) be fixed and suppose that \(k+1\ldots k+n+1\notin K_0\).
\begin{itemize}
	\item
		If \(c_1=0\) then \(K_1=\emptyset\) and all such indices belong to \(K_2\).
		Then,
  		\ifieee
	   		\begin{align*}
				\|Rx_{k+n+1}\|
			{}\overrel[\leq]{\eqref{eq:distR}}{} &
				2\alpha e_{k+n+1}
			{}\overrel[\leq]{\eqref{eq:K2linear}}{}
				2\alpha\rho^n
				e_{k+1}
			{}\overrel[\leq]{\eqref{eq:k+1bounds}}{}
				2\alpha B\rho^n
				e_k
			\\
			{}\overrel[\leq]{\eqref{eq:distR}}{} &
				2\alpha\gamma B\rho^n
				\|Rx_k\|.
			\end{align*}
		\else
	   		\[
				\addtolength\displayindent{\labelindent+\labelsep}
				\addtolength\displaywidth{-\labelindent-\labelsep}
				\|Rx_{k+n+1}\|
			{}\overrel[\leq]{\eqref{eq:distR}}{}
				2\alpha e_{k+n+1}
			{}\overrel[\leq]{\eqref{eq:K2linear}}{}
				2\alpha\rho^n
				e_{k+1}
			{}\overrel[\leq]{\eqref{eq:k+1bounds}}{}
				2\alpha B\rho^n
				e_k
			{}\overrel[\leq]{\eqref{eq:distR}}{}
				2\alpha\gamma B\rho^n
				\|Rx_k\|.
			\]
		\fi
		Since \(k+n+1\notin K_0\), then \(\|Rx_{k+n+1}\|\) failed the test at step \ref{step:SuperMann:K0} and therefore
		\[
   		\ifieee\else
			\addtolength\displayindent{\labelindent+\labelsep}
			\addtolength\displaywidth{-\labelindent-\labelsep}
		\fi
			c_0\|Rx_k\|
		{}={}
			c_0\eta_{k+n+1}
		{}<{}
			\|Rx_{k+n+1}\|
		{}\leq{}
			2\alpha\gamma B\rho^n
			\|Rx_k\|
		\]
		which proves that \(n\) cannot be arbitrarily large.
	\item
		If instead \(c_1>0\), let \(n_1\) be the number of indices among \(k+1\ldots k+n\) that belong to \(K_1\), and \(n_2=n-n_1\) those belonging to \(K_2\).
		Then, from iteration \(k+1\) to \(k+n+1\) the distance from the fixed set has reduced \(n_2\) times (at least) by a factor \(\rho\) and, due to \eqref{eq:k+1bounds}, increased at most by a factor \(B\) the remaining \(n_1\) times:
		\ifieee
			\begin{align*}
				\|Rx_{k+n+1}\|
			{}\overrel[\leq]{\eqref{eq:distR}}{} &
				2\alpha e_{k+n+1}
			{}\leq{}
				2\alpha\rho^{n_2}
				B^{n_1}
				e_{k+1}
			\\
			{}\overrel[\leq]{\eqref{eq:k+1bounds}}{} &
				2\alpha\rho^{n_2}
				B^{n_1+1}
				e_k
			{}\overrel[\leq]{\eqref{eq:distR}}{}
				2\alpha\gamma
				\rho^{n_2}
				B^{n_1+1}
				\|Rx_k\|.
			\end{align*}
		\else
			\begin{align*}
				\addtolength\displayindent{\labelindent+\labelsep}
				\addtolength\displaywidth{-\labelindent-\labelsep}
				\|Rx_{k+n+1}\|
			{}\overrel[\leq]{\eqref{eq:distR}}{} &
				2\alpha e_{k+n+1}
			{}\leq{}
				2\alpha\rho^{n_2}
				B^{n_1}
				e_{k+1}
			{}\overrel[\leq]{\eqref{eq:k+1bounds}}{}
				2\alpha\rho^{n_2}
				B^{n_1+1}
				e_k
			\\
			{}\overrel[\leq]{\eqref{eq:distR}}{} &
				2\alpha\gamma
				\rho^{n_2}
				B^{n_1+1}
				\|Rx_k\|.
			\end{align*}
		\fi
		Again, since \(k+n+1\notin K_0\) we have
		\(
			c_0\|Rx_k\|
		{}<{}
			2\alpha\gamma
			\rho^{n_2}
			B^{n_1+1}
			\|Rx_k\|
		\),
		and therefore
		\[
			\addtolength\displayindent{\labelindent+\labelsep}
			\addtolength\displaywidth{-\labelindent-\labelsep}
			n_1
		{}>{}
			\tfrac{
				\ln{\nicefrac{c_0}{2\alpha\gamma}}
			}{
				\ln B
			}
			{}-{}
			1
			{}+{}
			\tfrac{
				\ln{\nicefrac 1\rho}
			}{
				\ln B
			}
			n_2.
		\]
		In particular, for large \(n\) the number \(n_1\) of indices in \(K_1\) grows proportionally with respect to \(n\), and from \cref{prop:K1<2K0} we conclude once again that \(n\) cannot be arbitrarily large (since the number of visits to \(K_0\) does not change from \(k+1\) to \(k+n\)).
\end{itemize}

So far we proved that there exists \(n\in\N\) such that every \(n\) indices at least one belongs to \(K_0\).
In particular, indexing \(K_0=\set{k_0,k_1\cdots}\) we have that \(k_\ell\leq n\ell\), hence for all \(k_\ell\in K_0\)
\begin{equation}\label{eq:K0linearGlobal}
	\|Rx_{k_\ell}\|
{}\leq{}
	c_0^\ell\|Rx_0\|
{}\leq{}
	\bigl(c_0^{\nicefrac 1n}\bigr)^{k_\ell}
	\|Rx_0\|
\ifieee\else
	\qquad
	\forall k_\ell\in K_0
\fi
	.
\end{equation}
Moreover, any \(k\in\N\) is at most \(n-1\) indices away from the nearest previous index \(k_\ell\in K_0\); combined with \eqref{eq:K0linearGlobal} and invoking \eqref{eq:k+1bounds} we obtain
\[
	\|Rx_k\|
{}\leq{}
	B^{n-1}
	\|Rx_0\|
	\bigl(c_0^{\nicefrac 1n}\bigr)^{k_\ell}
{}\leq{}
	B^{n-1}
	\|Rx_0\|
	\bigl(c_0^{\nicefrac 1n}\bigr)^k
\]
proving the sought \(R\)-linear convergence of \(\seq{\|Rx_k\|}\).
It follows that for some \(b>0\) and \(r\in(0,1)\) we have \(\|Rx_k\|\leq br^k\) for all \(k\in\N\); then,
\ifieee
	\begin{align*}
		\|x_k-x_\star\|
	{}\leq{} &
		\sum_{j\geq k}\|x_{j+1}-x_j\|
	{}\leq{}
		D\sum_{j\geq k}\|Rx_j\|
	{}\leq{}
		bD\sum_{j\geq k}r^j
	\\
	{}={} &
		\frac{bD}{1-r}r^k
	\end{align*}
\else
	\[
		\|x_k-x_\star\|
	{}\leq{}
		\sum_{j\geq k}\|x_{j+1}-x_j\|
	{}\overrel[\leq]{\eqref{eq:D}}{}
		D\sum_{j\geq k}\|Rx_j\|
	{}\leq{}
		bD\sum_{j\geq k}r^j
	{}={}
		\frac{bD}{1-r}r^k
	\]
\fi
where in the second inequality we used the bound \eqref{eq:D}, which also holds for \(k\in K_2\) (up to possibly enlarging \(D\)) due to the fact that for \(k\in K_2\) under metric subregularity we have
\[
	\|x_{k+1}-x_k\|
{}\leq{}
	\|x_{k+1}-z_k\|
	{}+{}
	\|x_k-z_k\|
{}\leq{}
	2e_k
{}\overrel[\leq]{\eqref{eq:distR}}{}
	2\gamma\|Rx_k\|.
\]
This shows that \(\seq{x_k}\) is \(R\)-linearly convergent too.
\end{appendixproof}

		\proofsection{sec:SuperMann}
			\begin{appendixproof}{thm:SuperMann}

\noindent
Because of \cref{thm:LS} we know that for any direction \(d_k\) a feasible stepsize \(\tau_k\) complying with the requirements of step \ref{step:SuperMann:K2} will eventually be found, lower bounded as in \ref{thm:SuperMann:Tau_k} due to \cref{thm:LS} and \cref{ass:d}.
In particular, the scheme is well defined.
Moreover, from \cref{prop:DriftedKMFejer} we have that there exists a constant \(\underline\sigma>0\) such that
\[
	\|x_{k+1}-z\|^2
{}\leq{}
	\|x_k-z\|^2
	{}-{}
	\underline\sigma\|Rx_k\|^2
\ifieee\else
	\quad
	\text{for all \(k\in K_2\) and \(z\in\fix T\).}
\fi
\]
\ifieee
	for all \(k\in K_2\) and \(z\in\fix T\).
\fi
It follows that the \refSM[] is a special case of \cref{alg:General} and the proof entirely follows from \cref{{thm:General:Global},{thm:General:Linear}}.
\end{appendixproof}

\begin{appendixproof}{thm:dSuperlinear}
\begin{proofitemize}
	\item\ref{thm:dSuperlinear:tau}:
		let \(w_k^0\coloneqq x_k+d_k\).
		Superlinear convergence of \(\seq{d_k}\) then reads
		\(
			\frac{
				\|Rw_k^0\|
			}{
				\|Rx_k\|
			}
		{}\to{}
			0
		\).
		In particular, if \(c_1>0\) then there exists \(\bar k\in\N\) such that
		\(
			\|Rw_k^0\|
		{}\leq{}
			 c_1
			\|Rx_k\|
		\)
		for all \(k\geq\bar k\), \ie the point
		\(
			w_k^0
		{}={}
			x_k+d_k
		\)
		will always pass condition at step \ref{step:SuperMann:K1} resulting in
		\(
			x_{k+1}
		{}={}
			w_k^0
		{}={}
			x_k+d_k
		\)
		for all \(k\geq\bar k\).
		
		Similarly, if \(c_0>0\) then \(K_0\) is infinite as shown in \cref{thm:SuperMann:K0}; moreover, for \(\ell\in\N\)
		\[
			\frac{
				\|Rx_{k_\ell+1}\|
			}{
				\eta_{k_\ell+1}
			}
		{}={}
			\frac{
				\|Rx_{k_\ell+1}\|
			}{
				\|Rx_{k_\ell}\|
			}
		{}={}
			\frac{
				\|R(x_{k_\ell}+d_{k_\ell})\|
			}{
				\|Rx_{k_\ell}\|
			}
		{}\to{}
			0
		\ifieee\else
			\quad
				\text{as }~
				\ell\to\infty
		\fi
		\]
		\ifieee
			as \(\ell\to\infty\),
		\fi
		and therefore the ratio eventually is always smaller than \(c_0\), resulting in \(k_\ell+1\in K_0\) for \(\ell\) large enough.
		Consequently, the sequence will eventually reduce to \(x_{k+1}=x_k+d_k\).
	\item\ref{thm:dSuperlinear:Rx} and \ref{thm:dSuperlinear:x}:
		\(Q\)-superlinear convergence of the sequence \(\seq{Rx_k}\) follows from the fact that \(x_{k+1}=x_k+d_k\) for \(k\geq\bar k\).
		In particular, \(\seq{\|Rx_k\|}\) is summable and there exists a sequence \(\seq{\delta_k}\to 0\) such that \(\|Rx_{k+1}\|\leq\delta_k\|Rx_k\|\) for all \(k\).
		If \(\|d_k\|\leq D\|Rx_k\|\) for some \(D>0\), then
		\[
			\sum_{k\geq\bar k}\|x_{k+1}-x_k\|
		{}\leq{}
			D\sum_{k\geq\bar k}\|Rx_k\|
		{}<{}
			\infty
		\]
		which implies that \(\seq{x_k}\) is a Cauchy sequence, and hence converges to a point, be it \(x_\star\).
		Moreover, by possibly enlarging \(D\) so as to account for the iterates \(k<\bar k\), we have
		\ifieee
			\begin{align*}
				\|x_k-x_\star\|
			{}\leq{} &
				\sum_{j\geq k}{
					\|x_{j+1}-x_j\|
				}
			{}\leq{}
				D
				\sum_{j\geq k}{
					\|Rx_j\|
				}
			\\
			{}\leq{} &
				D\delta_0\delta_1\cdots\delta_{k-1}
				\sum_{j\in\N}{
					\|Rx_j\|
				}
			{}\eqqcolon{}
				\Delta_k.
			\end{align*}
		\else
			\[
				\|x_k-x_\star\|
			{}\leq{}
				\sum_{j\geq k}{
					\|x_{j+1}-x_j\|
				}
			{}\leq{}
				D
				\sum_{j\geq k}{
					\|Rx_j\|
				}
			{}\leq{}
				D\delta_0\delta_1\cdots\delta_{k-1}
				\sum_{j\in\N}{
					\|Rx_j\|
				}
			{}\eqqcolon{}
				\Delta_k.
			\]
		\fi
		This shows that \(\seq{x_k}\) is \(R\)-superlinearly convergent, since \(\nicefrac{\Delta_{k+1}}{\Delta_k}=\delta_k\to 0\).
	\item\ref{thm:dSuperlinear:K0}:
		already shown in the proof of \ref{thm:dSuperlinear:tau}.
	\qedhere
\end{proofitemize}
\end{appendixproof}

\begin{appendixproof}{thm:SuperMann:Broyden}
\added{%
	Let \(G_\star=JRx_\star\in\R^{n\times n}\) and let \(\|{}\cdot{}\|\) denote the Euclidean norm.
}%
From \cite[Lem. 2.2]{ip1992local} we have that there exist a constant \(L\) and a neighborhood \(U_{x_\star}\) of \(x_\star\) such that
\ifieee
	\begin{align*}
		\frac{
			\|y_k-G_\star s_k\|
		}{
			\|s_k\|
		}
	{}={} &
		\frac{
			\|Rw_k-Rx_k-G_\star(w_k-x_k)\|
		}{
			\|w_k-x_k\|
		}
	\\
	{}\leq{} &
		L
		\max{
			\set{
				\|x_k-x_\star\|
				{},{}
				\|w_k-x_\star\|
			}
		}.
	\end{align*}
\else
	\[
		\frac{
			\|y_k-G_\star s_k\|
		}{
			\|s_k\|
		}
	{}={}
		\frac{
			\|Rw_k-Rx_k-G_\star(w_k-x_k)\|
		}{
			\|w_k-x_k\|
		}
	{}\leq{}
		L
		\max{
			\set{
				\|x_k-x_\star\|
				{},{}
				\|w_k-x_\star\|
			}
		}.
	\]
\fi
Because of \eqref{eq:SuperMann:dBound}, the fact that \(\tau_k\leq 1\), and the triangular inequality we have
\(
	\|w_k-x_\star\|
{}\leq{}
	\|x_k-x_\star\|+D\|Rx_k\|
\)
and consequently
\[
	\sum_{k\in\N}{
		\frac{
			\|y_k-G_\star s_k\|
		}{
			\|s_k\|
		}
	}
{}\leq{}
	L
	\sum_{k\in\N}{
		\bigl(
			\|x_k-x_\star\|
			{}+{}
			D\|Rx_k\|
		\bigr)
	}
{}<{}
	\infty
\]
where the last inequality follows from \cref{thm:SuperMann:Linear}.

Let \(E_k=B_k-G_\star\) and let \(\|{}\cdot{}\|_F\) denote the Frobenius norm.
With a simple modification of the proofs of \cite[Thm. 4.1]{ip1992local} and \cite[Lem. 4.4]{artacho2014local} that takes into account the scalar \(\vartheta_k\in[\bar\vartheta,2-\bar\vartheta]\) we obtain
\ifieee
	\begin{align*}
		\left\|E_{k+1}\right\|_F
	{}\leq{} &
		\left\|
			E_k
			\left(
				\id-\vartheta_k\tfrac{s_k\trans{s_k}}{\|s_k\|^2}
			\right)
		\right\|_F
		{}+{}
		\vartheta_k
		\frac{\|y_k-G_\star s_k\|}{\|s_k\|}
	\\
	{}\leq{} &
		\hphantom{{}+{}}
		\left\|E_k\right\|_F
		{}-{}
		\frac{\bar\vartheta(2-\bar\vartheta)}{2\|E_k\|_F}
		\frac{\|E_ks_k\|^2}{\|s_k\|^2}
	\\
	&
		{}+{}
		(2-\bar\vartheta)
		\frac{\|y_k-G_\star s_k\|}{\|s_k\|}.
	\end{align*}
\else
	\begin{align*}
		\left\|E_{k+1}\right\|_F
	{}\leq{} &
		\left\|
			E_k
			\left(
				\id-\vartheta_k\tfrac{s_k\trans{s_k}}{\|s_k\|^2}
			\right)
		\right\|_F
		{}+{}
		\vartheta_k
		\frac{\|y_k-G_\star s_k\|}{\|s_k\|}
	\\
	{}\leq{} &
		\left\|E_k\right\|_F
		{}-{}
		\frac{\bar\vartheta(2-\bar\vartheta)}{2\|E_k\|_F}
		\frac{\|E_ks_k\|^2}{\|s_k\|^2}
	\end{align*}
\fi
The last term on the right-hand side, be it \(\sigma_k\), is summable and therefore the sequence \(\seq{E_k}\) is bounded.
Let \(\bar E\coloneqq\sup\seq{\|E_k\|_F}\), then
\begin{align*}
	\|E_{k+1}\|_F
	{}-{}
	\|E_k\|_F
{}\leq{} &
	\sigma_k
	{}-{}
	\frac{\bar\vartheta(2-\bar\vartheta)}{2\bar E}
	\left(
		\frac{\|(B_k-G_\star)s_k\|}{\|s_k\|}
	\right)^{\!2}.
\end{align*}
Telescoping the above inequality, summability of \(\sigma_k\) ensures that of
\(
	\frac{\|(B_k-G_\star)s_k\|^2}{\|s_k\|^2}
\)
proving in particular the claimed Dennis-Moré condition \eqref{eq:DM}.
\end{appendixproof}

	\end{appendix}

~~~~~~~~~~~~~~~~~~~~~~~~~~~~~~~~~~~~~~~~~~~~~~~~~~~~~~~~~~~~~~~~~~~~~~~~~~~ %

\ifieee
	\hspace*{0pt}
	\begin{IEEEbiography}[{%
		\includegraphics[width=1in,height=1.25in,clip,keepaspectratio]{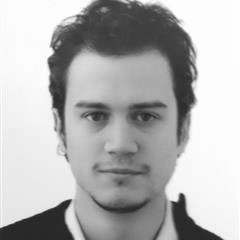}%
	}]{%
		Andreas Themelis received both Bachelor and Master degrees in Mathematics from the University of Florence, Italy, in 2010 and 2013, respectively.
		He is currently pursuing a joint Ph.D at the IMT School for Advanced Studies, Lucca (Italy) and the Department of Electrical Engineering (ESAT) of KU Leuven (Belgium).
		His research focuses on (non)convex nonsmooth optimization, with particular interest in splitting algorithms and monotone operators theory.
	}%
	\end{IEEEbiography}
	\begin{IEEEbiography}[{%
		\includegraphics[width=1in,height=1.25in,clip,keepaspectratio]{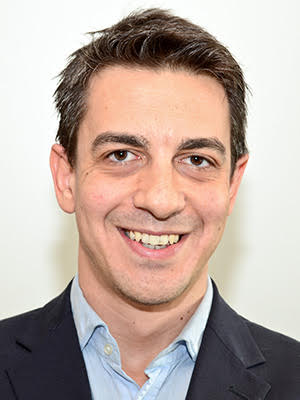}%
	}]{%
		Panagiotis (Panos) Patrinos is currently assistant professor at the Department of Electrical Engineering (ESAT) of KU Leuven, Belgium.
		He received the M.Eng. in Chemical Engineering, M.Sc. in Applied Mathematics and Ph.D. in Control and Optimization from National Technical University of Athens, Greece.
		After his Ph.D. he held postdoctoral positions at the University of Trento and IMT School of Advanced Studies Lucca, Italy, where he became an assistant professor in 2012.
		During fall/winter 2014 he held a visiting assistant professor position in the department of electrical engineering at Stanford University.
		His current research interests are in the theory and algorithms of optimization and predictive control with a focus on large-scale, distributed, stochastic and embedded optimization with a wide range of application areas including smart grids, water networks, aerospace, and machine learning.%
	}%
	\end{IEEEbiography}
\fi

\end{document}